\documentclass{article}

\usepackage {amsmath, amsthm, amssymb}
\usepackage [top=1.5in, right=1in, bottom=1.5in, left=1in] {geometry}
\usepackage {tikz}
\usetikzlibrary{calc}
\usepackage{hyperref}
\hypersetup{
    colorlinks = true,
    citecolor = blue,
}

\usepackage {libertine}
\usepackage [T1]{fontenc}

\theoremstyle {plain}
\newtheorem {thm}  {Theorem}
\newtheorem* {thmA} {Theorem A}
\newtheorem* {thmB} {Theorem B}
\newtheorem* {thmC} {Theorem C}
\newtheorem {prop} {Proposition}
\newtheorem {lem} {Lemma}
\newtheorem {cor} {Corollary}
\newtheorem* {prob} {Problem}

\theoremstyle {definition}
\newtheorem {ex} {Example}
\newtheorem {defn} {Definition}
\newtheorem {rmk} {Remark}

\newcommand {\n}{\mathbf{n}}
\newcommand {\m}{\mathbf{m}}

\title {Mixed Dimer Models for Euler and Catalan Numbers}
\author {Andrew Claussen and Nicholas Ovenhouse}
\date{}

\begin {document}

\maketitle

\begin {abstract}
    We study the enumeration of mixed dimer covers on skew Young diagrams of ribbon shape (also called border strips or snake graphs).
    For the two extreme cases of straight and zigzag shapes, we show that the number of mixed dimer covers are given
    by the Euler and Catalan numbers. We also give $q$-analogs by showing that the rank generating functions of the partial orders on
    mixed dimer covers agree with certain $q$-Euler and $q$-Catalan numbers. These $q$-analogs are a consequence of an isomorphism between
    the partial order on mixed dimer covers and the so-called \emph{middle order} on certain classes of permutations.
\end {abstract}

\setcounter{tocdepth}{1}
\tableofcontents

\section {Introduction and Main Results}

The Catalan numbers and Euler numbers are two very classical combinatorial sequences. 
The Catalan numbers $C_n = \frac{1}{n+1}\binom{2n}{n}$  count lattice paths, trees,
triangulations of polygons, and much more \cite{stanley_catalan}. The Euler numbers $E_n$ count alternating permutations, labeled trees,
standard Young tableaux of certain skew shapes, and other things \cite{stanley_alt}. Our main results are new interpretations of
the Euler and Catalan numbers as the number of \emph{mixed dimer covers} of certain planar graphs.

A mixed dimer cover is a generalization of a perfect matching (perfect matchings are also called \emph{dimer covers}). Given
non-negative integers assigned to the vertices of a graph, a mixed dimer cover is a multiset of edges so that every vertex has the prescribed degree.
If this degree is the same integer $n$ at all vertices, then we call the configurations \emph{$n$-dimer covers}.
It is similar to the notion of \emph{$\mathbf{d}$-factor} from \cite{propp_02}, the difference being that we have multisets of edges,
allowing each edge to be taken with some multiplicity. 

We consider certain skew Young diagrams called \emph{border strips} or \emph{ribbons}.
They are also frequently called \emph{snake graphs} (especially in the cluster algebra literature). 
The weighted enumeration of perfect matchings on snake graphs is related to many interesting parts of mathematics, including 
cluster algebras \cite{propp_20} \cite{ms_09} \cite{msw_13}, continued fractions \cite{snake_cluster}, Jones polynomials \cite{jones_poly},
Teichm\"{u}ller theory \cite{mw}, and $q$-deformed rational numbers \cite{mgo}.

There has recently been some interest in studying higher dimer covers on snake graphs.
In \cite{moz} it was shown that the coordinates on the higher Teichm\"{u}ller space for the orthosymplectic supergroup $\mathrm{OSp}(2|1)$ can be expressed
as weighted sums of double dimer covers on snake graphs. Then \cite{mosz} found enumerative formulas for $n$-dimer covers on snake graphs
for any $n \geq 1$. Also in \cite{musiker_wright}, type $D$ cluster variables were expressed in terms of mixed dimer covers on graphs which are
similar to snake graphs.
In some sense, this paper is a continuation and generalization of \cite{mosz} and \cite{bosz_24} to the case of mixed dimer covers.

\begin {figure}[h]
\centering
\begin {tikzpicture}
    \draw (0,0) -- (4,0) -- (4,1) -- (0,1) -- cycle;
    \draw (1,0) -- (1,1);
    \draw (2,0) -- (2,1);
    \draw (3,0) -- (3,1);

    \draw (0,0) node[below] {$1$};
    \draw (1,0) node[below] {$2$};
    \draw (2,0) node[below] {$3$};
    \draw (3,0) node[below] {$4$};
    \draw (4,0) node[below] {$5$};

    \draw (0,1) node[above] {$1$};
    \draw (1,1) node[above] {$2$};
    \draw (2,1) node[above] {$3$};
    \draw (3,1) node[above] {$4$};
    \draw (4,1) node[above] {$5$};

    \begin {scope}[shift={(6,-1)}]
        \draw (0,0) -- (1,0) -- (1,1) -- (2,1) -- (2,2) -- (3,2) -- (3,3) -- (1,3) -- (1,2) -- (0,2) -- cycle;
        \draw (0,1) -- (1,1) -- (1,2) -- (2,2) -- (2,3);

        \draw (0,0) node[below] {$1$};
        \draw (1,0) node[below] {$1$};
        \draw (0,1) node[left] {$2$};
        \draw (0,2) node[left] {$2$};
        \draw (1,1) node[below right] {$3$};
        \draw (2,1) node[below] {$3$};
        \draw (1,2) node[above left] {$4$};
        \draw (1,3) node[left] {$4$};
        \draw (2,2) node[below right] {$5$};
        \draw (3,2) node[below] {$5$};
        \draw (2,3) node[above] {$6$};
        \draw (3,3) node[above] {$6$};
    \end {scope}
\end {tikzpicture}
\caption {Examples of graphs $\mathcal{G}_4^s$ (left) and $\mathcal{G}_5^z$ (right), with their standard vertex labelings.}
\label {fig:snakes}
\end {figure}
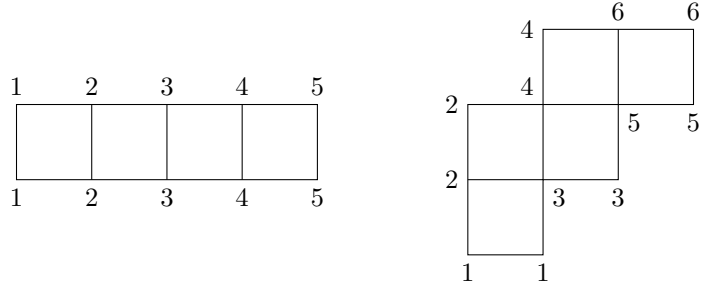

In the early sections of this paper, we will focus on the two simplest families of snake graphs.
The first family we consider are the straight snake graphs $\mathcal{G}_n^s$ (see Figure \ref{fig:snakes} (left)), 
which are $2 \times (n+1)$ grid graphs, formed by $n$ squares attached horizontally in a row.
The second family are the zig-zag snake graphs $\mathcal{G}_n^z$, whose squares are attached in a sequence going 
up, right, up, right, etc (see Figure \ref{fig:snakes} (right)).

These graphs have very simple and well-known formulas for their number of perfect matchings. The number of matchings of $\mathcal{G}_n^s$ are the Fibonacci numbers,
and the number of matchings of the zig-zag $\mathcal{G}_n^z$ is $n+1$.
Our main results are that with their \emph{standard labelings} (the vertex labels in Figure \ref{fig:snakes}), the number of mixed dimer covers on these graphs are
precisely the Euler and Catalan numbers.

\begin {thmA}
    Let $\mathcal{G}_n^s$ and $\mathcal{G}_n^z$ have the vertex labelings as depicted in Figure \ref{fig:snakes}. Then
    \begin {enumerate}
        \item[(a)] The number of mixed dimer covers of $\mathcal{G}_n^s$ is equal to the Euler number $E_{n+2}$.
        \item[(b)] The number of mixed dimer covers of $\mathcal{G}_n^z$ is equal to the Catalan number $C_{n+1}$.
    \end {enumerate}
\end {thmA}

Not only do we have these counting formulas, but we give explicit bijections between mixed dimer covers and certain families of permutations
counted by the Euler and Catalan numbers (namely alternating and 132-avoiding permutations). These bijections actually restrict to certain refinements of the Euler and Catalan numbers. The \emph{Entringer numbers} $E_{n,k}$ count 
alternating permutations $\sigma \in S_n$ with $\sigma(1) = k$. One of the many things counted by the Catalan numbers are $132$-avoiding
permutations. The \emph{generalized Catalan numbers} (also called \emph{ballot numbers}) $C_{n,k}$
count the $132$-avoiding permutations $\sigma \in S_n$ with $\sigma(1) = k$. 
We then have the following more refined version of Theorem A:

\begin {thmB}
    Let $\mathcal{G}_n^s$ and $\mathcal{G}_n^z$ have the vertex labelings as depicted in Figure \ref{fig:snakes}. Then
    \begin {enumerate}
        \item[(a)] The number of mixed dimer covers of $\mathcal{G}_n^s$ with $k$ dimers on the right-most vertical edge is equal to the Entringer number $E_{n+2,k+1}$.
        \item[(b)] The number of mixed dimer covers of $\mathcal{G}_n^z$ with $k$ dimers on the final edge\footnote{See section \ref{sec:catalan} for the 
                definition of \emph{final edge}.} is equal to the ballot number $C_{n+1,k+1}$. \\
    \end {enumerate}
\end {thmB}

The strategy of the proofs of Theorem A and Theorem B are to encode certain recurrences for the Euler/Entringer and and Catalan/ballot numbers using matrix multiplication, and then to show that the corresponding subsets of mixed dimer covers satisfy the same recurrences. We then use the recurrences to give explicit bijections. Under these bijections, one can read off the inversions of the associated permutation from certain edge multiplicities in the mixed dimer covers. 

We then consider a natural partial order on mixed dimer covers, and the associated rank generating functions, to obtain $q$-analogs of the above results.
In \cite{propp_02}, Propp gave a partial order on the set of perfect matchings (more generally $\mathbf{d}$-factors) of a planar graph. It is easy to see that
this partial order also makes sense for mixed dimer covers. We show that the rank generating functions 
agree with well-known $q$-analogs of the Euler and Catalan numbers. Said in another way,
the bijections from Theorem A and Theorem B are order-preserving with respect to certain partial orders on the permutations counted by Euler and Catalan numbers.

\begin {thmC}
    Let $\mathcal{G}_n^s$ and $\mathcal{G}_n^z$ be as before.
    \begin {enumerate}
        \item[(a)] The bijection from mixed dimer covers of $\mathcal{G}_n^s$ to alternating permutations in $S_{n+2}$ is order-preserving with respect to the Bruhat order.
        \item[(b)] The bijection from mixed dimer covers of $\mathcal{G}_n^z$ to $132$-avoiding permutations is order-preserving with respect to the Bruhat order.
    \end {enumerate}
\end {thmC}

More specifically, the bijections in Theorem C are poset isomorphims, with respect to the \emph{middle order} 
on permutations \cite{bft}, which is a coarsening of the Bruhat order. In the later sections of the paper, we discuss more general graphs
(of which the straight and zigzag snakes are special cases), and give an analog of Theorem C. For each graph in this more general family,
we explain how there is some associated subset of permutations such that the middle order on this subset is isomorphic to the lattice of mixed
dimer covers on the associated graph.

The rest of the paper is organized as follows.
In Section 2, we review the definitions of mixed dimer covers, and explain how for straight snake graphs we can count them
using certain matrix products. We then apply this to the \emph{standard labeling} of the vertices, and prove part $(a)$ of Theorems A and B, and derive an explicit bijection between mixed dimer covers and alternating permutations using the matrix recurrence. Lastly, we remark how one can get similar results for the Genocchi numbers using a slight modification.

In Section 3, we consider the case of zigzag shapes $\mathcal{G}^z_n$, and prove the analogous part $(b)$ of Theorems A and B for Catalan and ballot numbers. We also derive an explicit bijection between mixed dimer covers and 132-avoiding permutations.

In Section 4, we consider partial orders on the set of mixed dimer covers (generalizing the usual order on perfect matchings), and we show that the bijections from the previous sections are poset isomorphisms with respect to a certain order on permutations called the \emph{middle order} (because it is between the weak order and the Bruhat order). This gives an expression for certain $q$-analogs of Euler and Catalan numbers as rank generating functions for mixed dimer covers. As an application of our previous matrix recurrence formulas, we give special edge weights which allow the calculation of $q$-Euler and $q$-Catalan numbers using matrices.

In Section 5, we survey some other combinatorial models which are in bijection with mixed dimer covers. Most of the material and results in this section are not new, but we include them here for completeness,
and to showcase the relationship with our mixed dimer model.

Lastly, in Section 6, we discuss the generalization to arbitrary snake graphs (not just straight and zigzag shapes). We are not able to generalize all of the theorems from earlier sections, but we give some results, and pose some interesting questions for future study.

\section {Euler and Entringer Numbers from Straight Snakes}

\subsection {Mixed Dimer Covers on Straight Snakes} \label{sec:straight_snakes}

\begin {defn}
    Let $\mathcal{G} = (V,E)$ be a finite graph, and let $\n \colon V \to \Bbb{N}$ be a function which assigns a non-negative integer
    to each vertex. Define an \emph{$\n$-dimer cover} of $\mathcal{G}$ to be a multiset of edges such that each vertex $v$ is incident to 
    $\n(v)$ edges. Equivalently, it is a function $\m \colon E \to \Bbb{N}$ on edges such that for each $v \in V$, we have $\sum_{e \sim v} \m(e) = \n(v)$.
    We denote the set of $\n$-dimer covers of $\mathcal{G}$ by $\Omega_{\n}(\mathcal{G})$.
\end {defn}

\begin {rmk}
    When $\n$ is the constant function $\n(v) \equiv 1$, an $\n$-dimer cover is the same thing as a perfect matching.
    When $\n(v) \equiv n$, it is commonly called an \emph{$n$-dimer cover}. When $\n$ is not constant,
    we will call elements of $\Omega_\n ( \mathcal{G} )$ \emph{mixed dimer covers}.
\end {rmk}

\begin {figure}
\centering
\begin {tikzpicture}
    \draw (0,0) grid (3,1);
    \draw (0,0) node[below] {$v_{1,0}$};
    \draw (1,0) node[below] {$v_{1,1}$};
    \draw (2,0) node[below] {$v_{1,2}$};
    \draw (3,0) node[below] {$v_{1,3}$};

    \draw (0,1) node[above] {$v_{2,0}$};
    \draw (1,1) node[above] {$v_{2,1}$};
    \draw (2,1) node[above] {$v_{2,2}$};
    \draw (3,1) node[above] {$v_{2,3}$};
\end {tikzpicture}
\caption {The graph $\mathcal{G}^s_3$}
\label {fig:Gn_graph}
\end {figure}
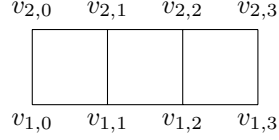

In this section, we will primarily be concerned with the following special case. Let $\mathcal{G}_n^s$ be the $2 \times (n+1)$ grid graph,
with vertices $v_{1,0},\dots,v_{1,n}$ in the bottom row and $v_{2,0},\dots,v_{2,n}$ in the top row (see Figure \ref{fig:Gn_graph}).
Given an integer sequence $m = (m_0,m_1,\dots,m_n) \in \Bbb{N}^{n+1}$, let $\n_m$ be the vertex labeling 
such that $\n_m(v_{1,k}) = \n_m(v_{2,k}) = m_k$ for all $k$.
We will present an enumerative formula for $\left| \Omega_{\n_m}(\mathcal{G}^s_n) \right|$ as a particular entry in a matrix product,
generalizing a formula from \cite{mosz} for the case $m = (n,n,n,\dots,n)$.

\begin {defn}
    Let $m = (m_0,m_1,\dots,m_n)$ be a sequence of vertex multiplicities for $\mathcal{G}^s_n$. 
    Define $\Omega^{ij}_{\n_m}(\mathcal{G}^s_n) := \Omega_{\n_{m'}}(\mathcal{G}^s_n)$,
    where $m' = (m_0+1-i, m_2, m_3, \dots, m_{n-1}, m_n+1-j)$. In other words, $\Omega^{11}_{\n_m}(\mathcal{G}^s_n) = \Omega_{\n_m}(\mathcal{G}^s_n)$,
    and incrementing $i$ (resp $j$) decreases $m_0$ (resp $m_n$). 
\end {defn}

\begin {defn}
    Let ${}^{ij}\Omega_{\n_m}(\mathcal{G}^s_n) \subseteq \Omega_{\n_m}(\mathcal{G}^s_n)$ denote the subset of mixed dimer covers 
    where the left-most vertical edge occurs with multiplicity $m_0+1-i$
    and the right-most vertical edge has multiplicity $m_n+1-j$. We also use the notations ${}^{i \bullet}\Omega_{\n_m}(\mathcal{G}^s_n)$ 
    and ${}^{\bullet j}\Omega_{\n_m}(\mathcal{G}^s_n)$ when only one of the two conditions is imposed.
\end {defn}

\begin {rmk} \label{rmk:straight_snake_bijection}
    Let $m_0,m_1,\dots,m_{n+1}$ be a sequence of vertex multiplicities, and let $m = (m_0,\dots,m_n)$ and $m' = (m_0,\dots,m_{n+1})$.
    Then for $k \leq m_{n+1}+1$, there is a bijection ${}^{\bullet k}\Omega_{\n_{m'}}(\mathcal{G}^s_{n+1}) \longrightarrow \Omega^{1 k}_{\n_m}(\mathcal{G}^s_n)$
    obtained by deleting the right-most two vertices of $\mathcal{G}^s_{n+1}$ (and their incident edges), and keeping all remaining edge multiplicities the same.
\end {rmk}

\begin {defn}
    Let $R_{a,b}$ be the $(a+1) \times (b+1)$ matrix whose $i,j$-entry is $1$ if $i+j \leq b+2$ and $0$ otherwise.
\end {defn}

\begin {ex}
    Some examples for different values of $a,b$ are given below:
    \[
        R_{2,2} = \begin{pmatrix} 1&1&1 \\ 1&1&0 \\ 1&0&0 \end{pmatrix}, \quad \quad 
        R_{2,4} = \begin{pmatrix} 1&1&1&1&1 \\ 1&1&1&1&0 \\ 1&1&1&0&0 \end{pmatrix}, \quad \quad 
        R_{3,1} = \begin{pmatrix} 1&1 \\ 1&0 \\ 0&0 \\ 0&0 \end{pmatrix}
    \]
\end {ex}

\begin {thm} \label{thm:matrix_formula}
    For an integer sequence $m = (m_0,m_1,\dots,m_n)$, define 
    the matrix $M_m$ by
    \[ M_m := R_{m_0,m_0}R_{m_0,m_1}R_{m_1,m_2} \cdots R_{m_{n-1},m_n} \]
    Then the $(i,j)$-entry of $M_m$ is equal to $\left| \Omega^{ij}_{\n_{m}}(\mathcal{G}^s_n) \right|$. 
    In particular, the $(1,1)$-entry is $\left| \Omega_{\n_m}(\mathcal{G}^s_n) \right|$, the total number of mixed dimer covers on $\mathcal{G}^s_n$ with vertex multiplicities $m$.
\end {thm}
\begin {proof}
    For the base case, where there is a single square, the matrix is
    \[ M = R_{m_0,m_0}R_{m_0,m_1} \]
    It is easy to check that this matrix has the following form. The last row and last column consist entirely of $1$'s.
    Removing the last row and column, we obtain a smaller matrix, and the last row and last column of this smaller matrix
    consists entirely of $2$'s. This process repeats all the way. In other words, the matrix looks as follows (if $m_0 < m_1$):
    \[
        M = \begin{pmatrix}
            m_0+1  & \cdots & m_0 + 1 & m_0    & \cdots & \cdots & 3      & 2      & 1 \\
            m_0    & \cdots & m_0     & m_0    & \cdots & \cdots & 3      & 2      & 1 \\
            \vdots & \vdots & \vdots  & \vdots & \ddots & \ddots & \vdots & \vdots & \vdots \\
            3      & \cdots & \cdots  & \cdots & \cdots & 3 & 3  & 2      & 1 \\
            2      & \cdots & \cdots  & \cdots & \cdots & 2 & 2  & 2      & 1 \\
            1      & \cdots & \cdots  & \cdots & \cdots & 1 & 1  & 1      & 1
        \end {pmatrix}
    \]
    If $m_0 > m_1$, then the matrix has the same form, but the upper-left corner has $m_1$ rather than $m_0$.
    We wish to compare the $(i,j)$-entry with $\left| \Omega^{ij}_{\n_m}(\mathcal{G}_1) \right|$. That is, we must count how many mixed dimer covers
    there are on $\mathcal{G}_1$ (which is a single square) in which the left two vertices have degree $m_0+1-i$ and the right two vertices
    have degree $m_1+1-j$. Consider first the case when $m_0-i < m_1 -j$ (i.e. $(i,j)$ is below the ``diagonal'' going through the bottom-right entry of the matrix). 
    Then there are exactly $m_0+2-i$ possible mixed dimer covers.
    Specifically, the left edge can have multiplicity $0,1,2,\dots,m_0+1-i$, and the rest is determined. 
    Since the part of the matrix below the diagonal has constant rows, we clearly see that
    the $(i,j)$ entry is $m_0+2-i$ in this case. The case when $m_0-i > m_1-j$ is similar.

    Now suppose that the result holds for the matrix product $M = R_{m_0,m_0}R_{m_0,m_1} \cdots R_{m_{n-1},m_n}$. We want to multiply on the right
    by $R_{m_n,m_{n+1}}$ to obtain $M' = M R_{m_n,m_{n+1}}$, and see that the entries of this matrix count mixed dimer covers on the larger graph.
    We will do so by checking that this matrix multiplication encodes a recurrence for the mixed dimer covers.

    Let $\mu$ be such that $j = m_{n+1}-m_n+\mu$, and let $\mu_+ := \mathrm{max}(\mu,0)$. Then from the structure of the $R_{m_n,m_{n+1}}$ matrix,
    we have $M'_{ij} = \sum_{k=1}^{m_n+1-\mu_+} M_{ik}$. We want to see that this is also the number of mixed dimer covers on $\mathcal{G}_{n+1}^s$ with
    vertex multiplicities $(m_0+1-i,m_1,m_2,\dots,m_n,m_{n+1}+1-j)$. Suppose the right-most vertical edge has multiplicity $\ell$. We can see that
    it must be in the range $m_{n+1}+1-j-m_n \leq \ell \leq m_{n+1}+1-j$. This lower bound could be negative, so we should also be careful to remember
    that $\ell \geq 0$.  
    By Remark \ref{rmk:straight_snake_bijection}, removing the last square gives a 
    bijection with $\Omega^{i,k}_{\n_m}(\mathcal{G}_n^s)$,
    where $k = m_{n+1}+2-j-\ell$. By induction, the number of such mixed dimer covers is $M_{ik}$. We have shown that $M_{ik}$ is
    the number of mixed dimer covers in $\Omega^{ij}_{\n_m}(\mathcal{G}^s_{n+1}) \cap {}^{\bullet k}\Omega_{\n_m}(\mathcal{G}^s_{n+1})$.
    Hence $\left| \Omega^{ij}_{\n_m}(\mathcal{G}^s_{n+1}) \right| = \sum_k M_{ik}$, summed over all admissible $k$.

    It remains to check that $k$ has the correct range of values.
    As $\ell$ ranges over $m_{n+1}+1-j-m_n \leq \ell \leq m_{n+1}+1-j$,
    we equivalently have $1 \leq k \leq m_n + 1$. The condition $\ell \geq 0$ also means $k \leq m_{n+1}+2-j = m_n+2-\mu$. If $\mu \leq 1$,
    this adds no extra restrictions, but if $\mu > 1$, this gives a smaller upper bound for $k$. So we have $1 \leq k \leq m_n+1-\mu_+$,
    and we recover the matrix multiplication formula.
\end {proof}

\bigskip

\begin {ex} \label{ex:matrix_product}
    For $(m_0,m_1) = (2,3)$, we have
    \[ M = R_{2,2}R_{2,3} = \begin{pmatrix} 3&3&2&1 \\ 2&2&2&1 \\ 1&1&1&1 \end{pmatrix} \]
    The upper-left entry is $3$, corresponding to the 3 mixed dimer covers shown below. 
    
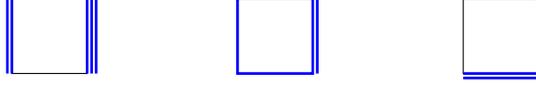
\begin {center}
    \begin {tikzpicture}
        \draw (0,0) -- (1,0) -- (1,1) -- (0,1) -- cycle;
        \draw[blue, line width=1] (0,0) -- (0,1);
        \draw[blue, line width=1] (-0.06,0) -- (-0.06,1);
        \draw[blue, line width=1] (1,0) -- (1,1);
        \draw[blue, line width=1] (1.06,0) -- (1.06,1);
        \draw[blue, line width=1] (1.12,0) -- (1.12,1);

        \begin {scope} [shift={(3,0)}]
            \draw (0,0) -- (1,0) -- (1,1) -- (0,1) -- cycle;
            \draw[blue, line width=1] (0,0) -- (1,0) -- (1,1) -- (0,1) -- cycle;
            \draw[blue, line width=1] (1.06,0) -- (1.06,1);
        \end {scope}

        \begin {scope} [shift={(6,0)}]
            \draw (0,0) -- (1,0) -- (1,1) -- (0,1) -- cycle;
            \draw[blue, line width=1] (0,0) -- (1,0) -- (1,1) -- (0,1);
            \draw[blue, line width=1] (0,-0.06) -- (1,-0.06);
            \draw[blue, line width=1] (0,1.06) -- (1,1.06);
        \end {scope}
    \end {tikzpicture}
    \end {center}
\end {ex}

\begin {ex}
    If $\n(v) \equiv 1$ for all vertices, then we are counting ordinary perfect matchings. 
    The matrix formula in this case becomes $M = R_{1,1}^{n+1} = \left(\begin{smallmatrix} 1&1\\1&0\end{smallmatrix}\right)^{n+1}$,
    a well-known way of generating the Fibonacci numbers.
\end {ex}

\bigskip

We can also upgrade Theorem \ref{thm:matrix_formula} to the weighted case. Putting arbitrary positive weights on the edges
of the graph $\mathcal{G}^s_n$, we can define generalizations of the $R_{a,b}$ matrices (depending on edge parameters) so that the entries of the matrix
product are the weighted sums of all corresponding mixed dimer covers, in which case Theorem \ref{thm:matrix_formula} will be recovered upon
setting all edge weights equal to $1$.

For a parameter $t$, let $U_n(t)$ be the $(n+1) \times (n+1)$ matrix whose $(i,j)$-entry is 0 if $i+j > n+2$ and otherwise $t^{n+2-i-j}$,
and let $T_{a,b}(t)$ be the $(a+1) \times (b+1)$ matrix with $(i,i)$-entry given by $t^{i-1}$ and all off-diagonal entries are zero. For example,
\[
    U_3(t) = \begin{pmatrix} t^3 & t^2 & t & 1 \\ t^2 & t & 1 & 0 \\ t & 1 & 0 & 0 \\ 1 & 0 & 0 & 0 \end{pmatrix}
    \quad \quad \text{and} \quad \quad
    T_{2,4}(t) = \begin{pmatrix} 1 & 0 & 0 & 0 & 0 \\ 0 & t & 0 & 0 & 0 \\ 0 & 0 & t^2 & 0 & 0 \end{pmatrix}
\]

Label the edges of $\mathcal{G}_n^s$ as follows.
The vertical edges (from left-to-right) will be labeled $a_0, a_1, \dots, a_n$. The top horizontal edges (from left-to-right) are $b_1,b_2,\dots,b_n$,
and the bottom horizontal edges $c_1,\dots,c_n$.
Then the analog of Theorem \ref{thm:matrix_formula} is given by the following. We omit the proof, as it is a straightforward generalization of the argument
from Theorem \ref{thm:matrix_formula}.

\begin {thm} \label{thm:weighted_matrix_product}
    For a vertex labeling $m = (m_0,\dots,m_n)$ on $\mathcal{G}^s_n$, and an edge weighting as described above (with edge weights $a_i,b_i,c_i$),
    define the matrix
    \[ X := U_{m_0}(a_0) \prod_{k=1}^n T_{m_{k-1},m_k}(b_kc_k) U_{m_k}(a_k) \]
    Then $X_{ij}$ is the weighted sum over mixed dimer covers in $\Omega^{ij}_{\n_m}(\mathcal{G}^s_n)$, where the weight of a mixed dimer
    is the product of its edge weights (counted with multiplicity).
\end {thm}

For instance, the weighted version of Example \ref{ex:matrix_product} is
\[ 
    X = \begin{pmatrix} a_0^2 & a_0 & 1 \\ a_0 & 1 & 0 \\ 1 & 0 & 0 \end{pmatrix}
    \begin{pmatrix} 1 & 0 & 0 & 0 \\ 0 & b_1c_1 & 0 & 0 \\ 0 & 0 & b_1^2c_1^2 & 0 \end{pmatrix}
    \begin{pmatrix} a_1^3 & a_1^2 & a_1 & 1 \\ a_1^2 & a_1 & 1 & 0 \\ a_1 & 1 & 0 & 0 \\ 1 & 0 & 0 & 0 \end{pmatrix}
\]
The $(1,1)$-entry of this matrix is $a_0^2a_1^3 + a_0a_1^2b_1c_1 + a_1b_1^2c_1^2$, which is indeed the weighted sum of the three
mixed dimer covers pictured above.

\bigskip

\subsection {Euler and Entringer Numbers} \label{sec:euler_straight_snakes}

\bigskip

\begin {defn}
    A permutation $\sigma \in S_n$ is called \emph{alternating} if $\sigma(1) > \sigma(2) < \sigma(3) > \sigma(4) < \cdots$.
    Let $\mathrm{Alt}_n \subset S_n$ be the set of alternating permutations, and $\mathrm{Alt}_{n,k}$ the set of alternating permtuations
    with $\sigma(1) = k$. The \emph{Euler numbers} are defined as $E_n := \left| \mathrm{Alt}_n \right|$, and the \emph{Entringer numbers} are
    $E_{n,k} := \left| \mathrm{Alt}_{n,k} \right|$.
\end {defn}

\begin {rmk} \label{rmk:ralt}
    If instead we have $\sigma(1) < \sigma(2) > \sigma(3) < \cdots$, then we call the permutation \emph{reverse alternating}, and write $\mathrm{RAlt}_n$ for the set of such permutations.
\end {rmk}

\bigskip

It was shown by Andr\'{e} \cite{andre} that the exponential generating function for the sequence $E_n$ is 
\[ \sum_{n \geq 0} E_n \frac{x^n}{n!} = 1 + x + \frac{x^2}{2!} + 2\frac{x^3}{3!} + 5\frac{x^4}{4!} + 16\frac{x^5}{5!} \dots = \text{sec}(x) + \text{tan}(x). \]
For this reason, the odd-indexed Euler numbers are sometimes called \textit{tangent numbers} and the even-indexed Euler numbers are called \textit{secant numbers}.

\begin {table}
\centering
\[
    \begin{array}{r|cccccc}
        n \setminus k & 1 & 2 & 3 & 4 & 5 & 6 \\ \hline
        1             & 1 &   &   &   &   &   \\
        2             & 0 & 1 &   &   &   &   \\
        3             & 0 & 1 & 1 &   &   &   \\
        4             & 0 & 1 & 2 & 2 &   &   \\
        5             & 0 & 2 & 4 & 5 & 5 &   \\
        6             & 0 & 5 & 10 & 14 & 16 & 16
    \end {array}
\]
\caption {The Entringer triangle}
\label {table:entringer}
\end{table}

The refinement $E_{n,k}$ was given initially by Seidel \cite{seidel} in the form of the triangular array shown in Table \ref{table:entringer}, which
we call the \emph{Entringer triangle}.
Seidel showed that the Entringer numbers satisfy a ``boustrophedon'' recurrence relation (Lemma \ref{lem:entringer}(a) below).  
Later, Entringer showed that the numbers $E_{n,k}$ count the number of alternating permutations with first entry $k$ \cite{entringer}. 

Alternating permutations have since been studied in the context of Morse theory \cite{arnold}, pattern avoidance \cite{lewis}, increasing trees \cite{kpp} \cite{stanley_alt}, and more. See \cite{stanley_EC1} and \cite{stanley_alt} for more information and history on alternating permutations and the Euler numbers.

The identities in the following lemma are well-known, and can be found in many of the references given above.

\begin {lem} \label{lem:entringer}
    The Entringer and Euler numbers satisfy the following recursions:
    \begin {enumerate}
        \item[(a)] $\displaystyle E_{n,k} = E_{n,k-1} + E_{n-1,n-k+1}$.
        \item[(b)] $\displaystyle E_{n,k} = \sum_{i > n-k} E_{n-1,i}$.
        \item[(c)] $\displaystyle E_{n+1,n+1} = E_{n}$.
    \end {enumerate}
\end {lem}
\begin {proof}
    We will just prove (b). Part (a) easily follows from (b), and (c) follows by taking $k=n$.

    Let $\sigma \in \mathrm{Alt}_{n,k}$. Then $\sigma(2)$ can be any $i < k$. Removing the $k$ from the beginning
    gives a reverse alternating permutation on the set $[n] \setminus \{k\}$. There is the obvious (order-preserving)
    bijection $[n] \setminus \{k\} \to [n-1]$ obtained by subtracting one from all $i > k$, and so we obtain a reverse alternating permutation
    on $[n-1]$. Taking the complementary permutation (i.e. $i \mapsto n-i$) gives an alternating permutation which begins with $n-\sigma(2)$.
    Since $\sigma(2) < k$, then $n-\sigma(2) > n-k$. Summing over all possibilities for $\sigma(2)$ gives part (b).
\end {proof}

Our main result is a new interpretation of the Euler and Entringer numbers in terms of mixed dimer covers. From now on, we will always consider
the \emph{standard labeling} of $\mathcal{G}_n^s$, given by $m = (m_0,m_1,\dots,m_n) = (1,2,3,\dots,n+1)$, and we will write $\Omega(\mathcal{G}_n^s)$ instead of
$\Omega_{\n_m}(\mathcal{G}_n^s)$, omitting the subscript.

\begin {thm} \label{thm:euler_bijection}
    With the standard labeling of $\mathcal{G}_n^s$, we have:
    \begin{enumerate}
        \item[(a)] $\left|\Omega(\mathcal{G}^s_{n-2})\right| = E_n$. 
        \item[(b)] $\left| {}^{\bullet k}\Omega(\mathcal{G}^s_{n-2}) \right| = \left| \Omega^{1 k}(\mathcal{G}^s_{n-3}) \right| = E_{n,n+1-k}$
    \end{enumerate}
\end {thm}
\begin {proof}
    The first equality $\left| {}^{\bullet k}\Omega(\mathcal{G}^s_{n-2}) \right| = \left| \Omega^{1 k}(\mathcal{G}^s_{n-3}) \right|$ in part $(b)$
    is simply an application of Remark \ref{rmk:straight_snake_bijection}. 
    Let $D_{n,k} := \left| {}^{\bullet,n+1-k}\Omega(\mathcal{G}^s_{n-2}) \right| = \left| \Omega^{1,n+1-k}(\mathcal{G}^s_{n-3}) \right|$ be the number of mixed dimer covers
    of $\mathcal{G}^s_{n-2}$ for which the right-most vertical edge has multiplicity $k-1$, and also the number of mixed dimer covers of $\mathcal{G}^s_{n-3}$ with
    vertex labeling $(1,2,\dots,n-3,k-2)$.
    To prove part $(b)$, we must show that $D_{n,k} = E_{n,k}$, which we will do by analyzing the matrix product formula from Theorem \ref{thm:matrix_formula}
    and seeing that it encodes the same recurrence from Lemma \ref{lem:entringer} (b).

    By Theorem \ref{thm:matrix_formula}, we have that $D_{n,k}$ is equal to the $(1,n+1-k)$-entry
    of the matrix $M_{n-2} = R_{1,1}R_{1,2}R_{2,3}\cdots R_{n-3,n-2}$.
    Since column $n+1-k$ of $R_{n-3,n-2}$ is given by $(1,1,\dots,1,0,0,\dots,0)^\top$, with the first $k-1$ entries equal to $1$, we
    get the recurrence 
    \[ D_{n,k} = \sum_{j=1}^{k-1}(M_{n-3})_{1,j} = \sum_{j=1}^{k-1} D_{n-1,n-j} = \sum_{j=n+1-k}^{n-1}D_{n-1,j} \]
    This is identical to the recurrence in Lemma \ref{lem:entringer}(b). It is simple to check that $E_{n,k}$ and $D_{n,k}$ have the same initial conditions.
    This proves part $(b)$ of the theorem, and part $(a)$ follows by taking the union/sum over $k$.
\end {proof}

\subsection {An Explicit Bijection} \label{sec:bijection}

The equality in Theorem \ref{thm:euler_bijection} followed from the fact that the sizes of certain subsets of mixed dimer covers
satisfy the same recurrence as the Entringer numbers. Also recall from the proof of Lemma \ref{lem:entringer}(b) the explicit bijection
$\Phi \colon \mathrm{Alt}_{n,k} \to \bigcup_{i > n-k} \mathrm{Alt}_{n-1,i}$ which realizes this recurrence: for $\sigma \in \mathrm{Alt}_{n,k}$, remove
the first entry to get $(\sigma_2,\sigma_3,\dots,\sigma_n)$, and then standardize (i.e. subtract $1$ from all $\sigma_i > \sigma_1$) to obtain
a permtuation $\bar{\sigma} \in S_{n-1}$. This will be reverse alternating, so finally multiply on the left by $w_0$ 
(i.e. change all $\bar{\sigma}_i$ to $n-\bar{\sigma}_i$) to obtain another alternating permutation $\tau$ (with $\tau_1 > n-k$).
We will carefully compare this bijective recurrence with a similar one on mixed dimer covers
to give an explicit bijection realizing the equality in Theorem \ref{thm:euler_bijection}.

The mapping above can be described nicely in terms of Lehmer codes. 
Recall that the Lehmer code of a permutation is
\[ L(\sigma) = (L_1,\dots,L_n) \quad \text{ where } L_i := \#\{j>i ~|~ \sigma(j) < \sigma(i) \}. \]
Associating a permutation with its Lehmer code gives a bijection $S_n \to [0,n-1] \times [0,n-2] \times \cdots \times [0,0]$.
In the first step, $\sigma \mapsto \bar{\sigma}$, no inversions are created or destroyed (other than those involving $\sigma_1$), and so if $L(\sigma) = (L_1,\dots,L_n)$ and $L(\bar{\sigma}) = (\bar{L}_1,\dots,\bar{L}_{n-1})$ are the Lehmer codes, we clearly have $\bar{L}_i = L_{i+1}$. In the last step, mapping $\bar{\sigma} \mapsto \tau = w_0 \bar{\sigma}$,
has the effect that every inversion becomes a non-inversion and vice versa. That is, $\bar{\sigma}_i > \bar{\sigma}_j$ if and only if $\tau_i < \tau_j$.
The Lehmer codes $L(\bar{\sigma}) = (\bar{L}_1,\dots,\bar{L}_{n-1})$ and $L(\tau) = (L'_1,\dots,L'_{n-1})$ are therefore related
by $L'_i = n-1-i-\bar{L}_i$. To summarize, the mapping $\mathrm{Alt}_{n,k} \to \bigcup_{i > n-k} \mathrm{Alt}_{n-1,i}$ translates
to the mapping $(L_1,L_2,\dots,L_n) \mapsto (n-2-L_2, n-3-L_3,\dots,1-L_{n-1},0)$ on Lehmer codes.

\begin {figure}[h]
\centering
\begin {tikzpicture}[scale=2]
    \draw (0,0) -- (3,0) -- (3,1) -- (0,1);
    \foreach \x in {1,2} {\draw (\x,0) -- (\x,1);}
    \draw (0,0.5) node {$\cdots$};

    \draw[blue] (3,1.4) node {\small $n-1$};
    \draw[blue] (2,1.4) node {\small $n-2$};
    \draw[blue] (1,1.4) node {\small $n-3$};

    \draw (3,0.5) node[right] {$\ell_1$};
    \draw (2,0.5) node[right] {$\ell_2$};
    \draw (1,0.5) node[right] {$\ell_3$};

    \draw (2.5,1) node[above] {\small $n-1-\ell_1$};
    \draw (1.5,1) node[above] {\small $\ell_1-\ell_2-1$};

    \draw[-latex] (3.5,0.5) -- (4.5,0.5);

    \begin {scope}[shift={(5,0)}]
        \draw (0,0) -- (2,0) -- (2,1) -- (0,1);
        \draw (1,0) -- (1,1);
        \draw (0,0.5) node {$\cdots$};
    
        \draw[blue] (2,1.4) node {\small $n-2$};
        \draw[blue] (1,1.4) node {\small $n-3$};
    
        \draw (2,0.5) node[right] {$n-1-\ell_1+\ell_2$};
        \draw (1,0.5) node[right] {$\ell_3$};
    
        \draw (1.5,1) node[above] {\small $\ell_1-\ell_2-1$};
    \end {scope}
\end {tikzpicture}
\caption {The bijection $\Psi \colon \mathcal{D}_{n,k} \to \bigcup_{i > n-k} \mathcal{D}_{n-1,i}$. Vertex labels are in blue,
and edge labels are the multiplicities in a mixed dimer cover.}
\label{fig:dimer_recurrence_bijection}
\end {figure}
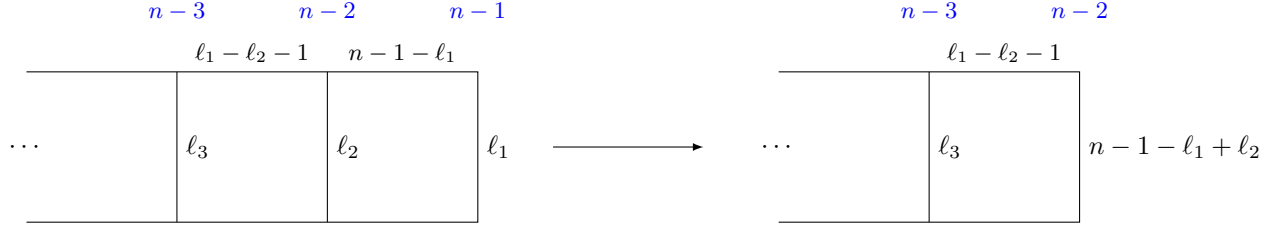

Now let us discuss a similar bijective recurrence for mixed dimer covers. 
Let $\mathcal{D}_{n,k} := {}^{\bullet, n+1-k}\Omega(\mathcal{G}^s_{n-2})$ be the set of mixed dimer covers of $\mathcal{G}^s_{n-2}$ 
with the right-most vertical edge having multiplicity $k-1$.
We will define a bijection $\Psi \colon \mathcal{D}_{n,k} \to \bigcup_{i > n-k} \mathcal{D}_{n-1,i}$. This bijection is illustrated in
Figure \ref{fig:dimer_recurrence_bijection}, and we will explain it now. Let $\ell_1,\ell_2,\dots$ be the multiplicities of the vertical
edges of $\mathcal{G}^s_{n-2}$, ordered from right-to-left. Since the right-most vertical edge has multiplicity $\ell_1$,
then the horizontal edges on the last square must have multiplicity $n-1-\ell_1$. We define the mapping $\Psi$ by sending this mixed dimer cover
to the one on $\mathcal{G}^s_{n-3}$ which is identical to the original mixed dimer cover (on all but the last square) except that its right-most vertical
edge has multiplicity $n-1-\ell_1+\ell_2$ (rather than just $\ell_2$). In this way, we obtain a mixed dimer cover whose right-most vertical edge
has multiplicity at least $n-1-\ell_1$, and thus the map $\Psi$ is well-defined.

We now wish to define a bijection $F \colon \mathrm{Alt}_{n,k} \to \mathcal{D}_{n,k}$ which commutes with the recursive bijections discussed above.
That is, we would like to have $F \circ \Phi = \Psi \circ F$. This condition will uniquely determine the bijection. Note that an element 
$\omega \in \mathcal{D}_{n,k}$ will have $\ell_1 = k-1$ by definition,
and also permutations $\sigma \in \mathrm{Alt}_{n,k}$ have $\sigma_1 = k$ and $L_1 = k-1$. We therefore have that $\ell_1 = L_1$.
As discussed above, $L(\Phi(\sigma)) = (n-2-L_2,\dots)$. On the other hand, the right-most vertical edge of $\Psi(\omega)$ has multiplicity
$n-1-\ell_1+\ell_2$. Requiring that $F \circ \Phi = \Psi \circ F$ therefore imposes the condition $n-2-L_2 = n-1-\ell_1+\ell_2$. Simplifying,
and using the fact that $\ell_1 = L_1$, we get $\ell_2 = L_1 - L_2 - 1$, and thus the multiplicity of the second vertical edge is determined.

\begin {figure}
\centering
\begin {tikzpicture}[scale=1.5]
\draw (0,0) grid (5,1);
\draw (5,0.5) node[right] {$L_1$};
\draw (4,0.5) node[right] {$\Delta_2$};
\draw (3,0.5) node[right] {$\delta_3$};
\draw (2,0.5) node[right] {$\Delta_4$};
\draw (1,0.5) node[right] {$\delta_5$};
\draw (4.5,1) node[above] {$6-L_1$};
\draw (4.5,0) node[below] {$6-L_1$};
\draw (3.5,1) node[above] {$L_2$};
\draw (3.5,0) node[below] {$L_2$};
\draw (2.5,1) node[above] {$4-L_3$};
\draw (2.5,0) node[below] {$4-L_3$};
\draw (1.5,1) node[above] {$L_4$};
\draw (1.5,0) node[below] {$L_4$};
\draw (0.5,1) node[above] {$2-L_5$};
\draw (0.5,0) node[below] {$2-L_5$};
\draw (0,0.5) node[left] {$L_5-1$};
\end {tikzpicture}
\caption {The bijection $F \colon \mathrm{Alt}_7 \to \Omega(\mathcal{G}^s_{5})$. The multiplicities of the vertical edges are given by $\delta_i = L_i-L_{i-1}$ 
and $\Delta_i = L_{i-1}-L_i-1$, where $L(\sigma) = (L_1,L_2,\dots,L_7)$ is the Lehmer code of an alternating permutation. 
Note that the first and last vertical edges also follow this pattern if we define $L_0 := 0$.}
\label{fig:lehmer}
\end {figure}
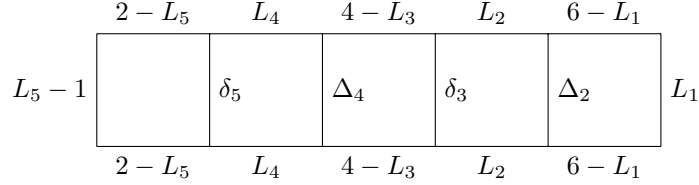

We may continue this process to express all the $\ell_i$ in terms of the $L_i$. Doing so gives
\[ \ell_{2i+1} = L_{2i+1}-L_{2i}, \quad \ell_{2i} = L_{2i-1}-L_{2i}-1 \]
Knowing the multiplicities of these vertical edges uniquely determines the rest of the mixed dimer cover.
The end result in the special case of $n=7$ is picture in Figure \ref{fig:lehmer}.

\begin {rmk}
    The inverse bijection $F^{-1} \colon \Omega(\mathcal{G}^s_{n-2}) \to \mathrm{Alt}_n$ is given as follows. 
    If $\omega \in \Omega(\mathcal{G}^s_{n-2})$ has multiplicities $\ell_i$ on the vertical edges, then the Lehmer code $L(\sigma) = (L_1,\dots,L_n)$
    of $F^{-1}(\omega)$ is given by 
    \[ L_i = \left( \sum_{j=1}^i (-1)^{j+1} \ell_j \right) - \left\lfloor \frac{i}{2} \right\rfloor \]
\end {rmk}

\bigskip

\begin {ex} \label{ex:alt_mixed_bijection}
    There are $E_4=5$ alternating permutations in $S_4$: $2143$, $3142$, $3241$, $4132$, $4231$.
    We also see that the Entringer numbers are $E_{4,1} = 0$, $E_{4,2} = 1$, $E_{4,3} = 2$, and $E_{4,4} = 2$.
    Figure \ref{fig:euler_dimer_example} shows the 5 mixed dimer covers of $\mathcal{G}^s_2$, along with the associated permutations
    and their Lehmer codes.

    \begin {figure}
    \centering
    \begin {tikzpicture}
        \foreach \x in {0, 3, 6, 9, 12} {
            \draw (\x,0) -- (\x+2,0) -- (\x+2,1) -- (\x,1) -- cycle;
            \draw (\x+1,0) -- (\x+1,1);
        }

        \begin {scope}[shift={(0,0)}]
            \draw[red, line width=1.5] (0,0) -- (0,1);
            \draw[red, line width=1.5] (1,0) -- (2,0);
            \draw[red, line width=1.5] (1,-0.1) -- (2,-0.1);
            \draw[red, line width=1.5] (1,1) -- (2,1);
            \draw[red, line width=1.5] (1,1.1) -- (2,1.1);
            \draw[red, line width=1.5] (2,0) -- (2,1);

            \draw (1,-1) node {$L = (1,0,1,0)$};
            \draw (1,-1.5) node {$\sigma = 2143$};
        \end {scope}

        \begin {scope}[shift={(3,0)}]
            \draw[red, line width=1.5] (0,0) -- (0,1);
            \draw[red, line width=1.5] (1,0) -- (2,0);
            \draw[red, line width=1.5] (1,1) -- (2,1);
            \draw[red, line width=1.5] (1,0) -- (1,1);
            \draw[red, line width=1.5] (2,0) -- (2,1);
            \draw[red, line width=1.5] (2.1,0) -- (2.1,1);
            
            \draw (1,-1) node {$L = (2,0,1,0)$};
            \draw (1,-1.5) node {$\sigma = 3142$};
        \end {scope}

        \begin {scope}[shift={(6,0)}]
            \draw[red, line width=1.5] (0,0) -- (2,0) -- (2,1) -- (0,1);
            \draw[red, line width=1.5] (2.1,0) -- (2.1,1);

            \draw (1,-1) node {$L = (2,1,1,0)$};
            \draw (1,-1.5) node {$\sigma = 3241$};
        \end {scope}

        \begin {scope}[shift={(9,0)}]
            \draw[red, line width=1.5] (0,0) -- (1,0) -- (1,1) -- (0,1);
            \draw[red, line width=1.5] (2,0) -- (2,1);
            \draw[red, line width=1.5] (2.1,0) -- (2.1,1);
            \draw[red, line width=1.5] (2.2,0) -- (2.2,1);

            \draw (1,-1) node {$L = (3,1,1,0)$};
            \draw (1,-1.5) node {$\sigma = 4231$};
        \end {scope}

        \begin {scope}[shift={(12,0)}]
            \draw[red, line width=1.5] (0,0) -- (0,1);
            \draw[red, line width=1.5] (1,0) -- (1,1);
            \draw[red, line width=1.5] (1.1,0) -- (1.1,1);
            \draw[red, line width=1.5] (2,0) -- (2,1);
            \draw[red, line width=1.5] (2.1,0) -- (2.1,1);
            \draw[red, line width=1.5] (2.2,0) -- (2.2,1);

            \draw (1,-1) node {$L = (3,0,1,0)$};
            \draw (1,-1.5) node {$\sigma = 4132$};
        \end {scope}
    \end {tikzpicture}
    \caption {The 5 elements of $\Omega(\mathcal{G}^s_2)$.}
    \label{fig:euler_dimer_example}
    \end {figure}
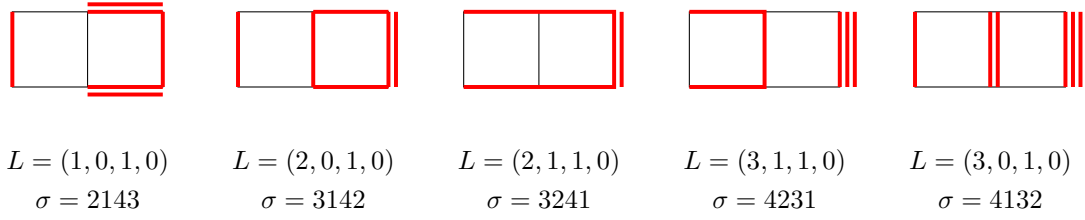
\end {ex}

\begin {ex}
    For an illustration of the bijection for $E_5$, see Figure \ref{fig:alt5_poset}.
\end {ex}

\begin {rmk}
    Let $M = R_{11}R_{12}R_{23} \cdots R_{n-3,n-2}$ be the matrix product from Theorem \ref{thm:matrix_formula} for the
    standard labeling $m = (1,2,3,\dots,n-2)$ of $\mathcal{G}^s_{n-3}$. Then Theorem \ref{thm:euler_bijection}$(b)$ says that the entries in the first row
    of $M$ are the Entringer numbers $E_{n,n}$, $E_{n,n-1}$, $\dots$, $E_{n,3}$, $E_{n,2}$. Using the bijection described in this section,
    we can also give a combinatorial description of the entries in the second row of $M$. By Theorem \ref{thm:matrix_formula}, $M_{2,k}$
    counts mixed dimer covers on $\mathcal{G}^s_{n-3}$ with vertex multiplicities $m = (0,2,3,\dots,n-3,n-1-k)$. Similarly to Remark \ref{rmk:straight_snake_bijection},
    there is a bijection $\Omega^{2,k}(\mathcal{G}^s_{n-3}) \to {}^{1,k}\Omega(\mathcal{G}^s_{n-2})$.
    By the Lehmer code description given in this section, this means $L_{n-2} = 2$ if $n$ is odd, and $L_{n-2} = 0$ if $n$ is even.
    In other words, $M_{2,k}$ counts the subset of $\mathrm{Alt}_{n,n+1-k}$ such that $(-1)^n(\sigma_n-\sigma_{n-2}) > 0$.
\end {rmk}

\subsection {Genocchi Numbers from Straight Snakes}

Consider the \emph{Seidel triangle} in Table \ref{table:genocchi} and denote the entries by $g_{n,k}$. They are defined by the recurrence
\[ g_{n,k} = \sum_{i > \lfloor \frac{n}{2} \rfloor - k} g_{n-1,i} \]

\begin {table}
\centering
\begin {tabular} {c|cccc}
    $n/k$ & 1 & 2 & 3 & 4 \\ \hline
    1     & 1 &   &   &   \\
    2     & 1 &   &   &   \\
    3     & 1 & 1 &   &   \\
    4     & 1 & 2 &   &   \\
    5     & 2 & 3 & 3 &   \\
    6     & 3 & 6 & 8 &   \\
    7     & 8 & 14 & 17 & 17 \\
    8     & 17 & 34 & 48 & 56 
\end {tabular}
\caption {The Seidel triangle for the Genocchi numbers.}
\label {table:genocchi}
\end {table}

\begin {rmk}
    Our triangle differs slightly from what usually appears in the literature (see e.g. \cite{es_00}). In particular,
    we have reversed the order of the even rows. This indexing will be more convenient for our purposes.
\end {rmk}

\begin {defn}
    The entries $G_n := g_{n, \lceil n/2 \rceil}$ at the ends of the rows of the triangle are called the \emph{Genocchi numbers}. 
    More specifically, $G_{2n-1} = g_{2n-1,n}$ are the \emph{Genocchi numbers of the first kind}, and $G_{2n} = g_{2n,n}$ are
    the \emph{Genocchi numbers of the second kind} (also called \emph{median Genocchi numbers}).
\end {defn}

The Genocchi numbers are very closely related to the Euler numbers. Indeed, the recurrence above is very similar to the one for the Entringer numbers.
Moreover, the Genocchi numbers of the first kind and the odd Euler numbers (i.e. the \emph{tangent numbers}) are related
by $G_{2n-1} = \frac{n}{4^{n-1}} E_{2n-1}$.

Several combinatorial interpretations for these numbers are known. They count \emph{Dumont permutations} \cite{dumont}, 
\emph{collapsed permutations} \cite{aht_22}, alternating pistols \cite{dv_80}, and certain restricted semistandard tableaux of skew shape \cite{mz_22}.
It was also shown in \cite{dv_80} that they count alternating permutations whose Lehmer codes have all even entries. We will give
another combinatorial interpretation of the Genocchi numbers (and indeed all entries of the triangle $g_{n,k}$) in terms of mixed dimer covers.
The methods are the same as those used earlier for the Euler and Entringer numbers, and we omit the proofs.

We will consider the straight snake graphs $G^s_n$ with vertex multiplicities $m = (1,1,2,2,\dots,k,k)$ and $m = (1,1,2,2,\dots,k,k,k+1)$.
We then have the following results, analogous to Theorem \ref{thm:euler_bijection}.

\begin {thm}
    Let $m = (1,1,2,2,\dots,k,k)$ or $m = (1,1,2,2,\dots,k,k,k+1)$. Then
    \begin {enumerate}
        \item[(a)] $\left| \Omega_{\n_m}(\mathcal{G}^s_{n-3}) \right| = G_n$
        \item[(b)] $\left| {}^{\bullet k}\Omega_{\n_m}(\mathcal{G}^s_{n-2}) \right| = \left| \Omega^{1,k}_{\n_m}(\mathcal{G}^s_{n-3}) \right| = g_{n,\lceil n/2 \rceil + 1 - k}$
    \end {enumerate}
\end {thm}

\section{Catalan and Ballot Numbers from Zigzag Snakes} \label{sec:catalan}

\subsection {Catalan and Ballot Numbers}

The Catalan numbers $C_n$ were first systematically studied by Euler, Goldbach, and Segner, and initially were shown to count the number of triangulations of an $n$-gon. This ubiquitous number sequence also counts Dyck paths, linear extensions of the $2 \times n$ rectangle poset, Young tableaux 
that fit into the staircase shape $(n-1,n-2,...,1)$, order ideals of the type $A$ root poset, and $132$-avoiding permutations, to name a few.
For more about the Catalan numbers, see \cite{stanley_catalan} or \cite{pak}.

\begin {defn}
    A permutation $\tau \in S_n$ is called \emph{132-avoiding} if for any $j < k < \ell$, we never have $\tau(j) < \tau(\ell) < \tau(k)$. 
    Let $\mathrm{Cat}_n \subset S_n$ be the set of $132$-avoiding permutations, and $\mathrm{Cat}_{n,k}$ the set of $132$-avoiding permutations
    with $\tau(1) = k$. We define the \emph{ballot numbers} to be
    $C_{n,k} := \left| \mathrm{Cat}_{n,k} \right|$.
\end {defn}

The Catalan numbers are given by the formula $C_n = \frac{1}{n+1}\binom{2n}{n}$, and the ballot numbers are  
$C_{n,k} = \frac{n-(k-1)}{n+k-1}\binom{n+k-1}{n} = \frac{n-k+1}{n}\binom{n+k-2}{k-1}$.
The ballot numbers can be arranged into a triangle similar to the Entringer numbers (see Table \ref{table:ballot}). The following result lists some recurrences satisfied by these numbers. We omit the proofs, as they are well-known and easy to check.

\begin {table}
\centering
\[
    \begin{array}{r|cccccc}
        n \setminus k & 1 & 2 & 3 & 4 & 5 \\ \hline
        1             & 1 &   &   &   &   \\
        2             & 1 & 1 &   &   &   \\
        3             & 1 & 2 & 2 &   &   \\
        4             & 1 & 3 & 5 & 5 &   \\
        5             & 1 & 4 & 9 & 14 & 14
    \end {array}
\]
\caption {The Ballot numbers}
\label {table:ballot}
\end{table}

\begin {lem} \label{lem:ballot}
    The Ballot and Catalan numbers satisfy the following recursions:
    \begin {enumerate}
        \item[(a)] $\displaystyle C_{n,k} = C_{n,k-1} + C_{n-1,k}$.
        \item[(b)] $\displaystyle C_{n,k} = \sum_{i=1}^{k} C_{n-1,i}$.
        \item[(c)] $\displaystyle C_{n+1,n+1} = C_{n}$.
    \end {enumerate}
\end {lem}

\subsection {Enumeration for Zigzag Snakes}

We define the \emph{standard labeling} of the zigzag graph $\mathcal{G}_n^z$ to be the one pictured in Figure \ref{fig:snakes}. The horizontal edges along the bottom boundary
have their endpoints labeled with odd numbers, and the vertical edges along the top boundary have their endpoints labeled with even numbers. 
Note that the top (respectively right) edge of the last square when $n$ is odd (resp. even) is an exception to this rule. 
That is, the right vertical edge of the last tile of $\mathcal{G}^z_{2n}$ is labeled by the odd number $2n+1$, and likewise the top horizontal edge of 
the last tile of $\mathcal{G}^z_{2n-1}$ is labeled by the even number $2n$.

From now on, when considering the standard labeling of $\mathcal{G}_n^z$ we will write $\Omega(\mathcal{G}_n^z)$ rather than $\Omega_\n(\mathcal{G}_n^z)$,
omitting the subscript.
Let the \emph{final edge} be the right vertical (resp. top horizontal) edge of the last tile when $n$ is odd (resp. even). 
Note that with the standard labeling, the endpoint vertices of the final edge are labeled $n$ and $n+1$. Similar to the notation used for straight snakes,
we will write ${}^{\bullet k}\Omega(\mathcal{G}^z_n)$ for the set of mixed dimer covers where the final edge has multiplicity $m_n+1-k$, and we write
$\Omega^{ij}(\mathcal{G}^z_n)$ for the set of mixed dimer covers of $G^z_n$ with vertex multiplicities $(2-i,2,3,\dots,n,n+2-j)$. 

\begin {rmk}
    Similar to Remark 2, we have for zigzags that $\left| {}^{\bullet k}\Omega(\mathcal{G}^z_n) \right| = \left| \Omega^{1,k-1}(\mathcal{G}^z_{n-1}) \right|$,
    and the bijection is obtained by deleting the vertices of the final edge.
\end {rmk}

\begin {thm} \label{thm:dimer_catalan_bijection}
    For the standard labeling on $\mathcal{G}_n^z$, we have: 
    \begin{enumerate}
        \item[(a)] $\left|\Omega(\mathcal{G}^z_{n-1})\right| = C_n$. 
        \item[(b)] $\left|{}^{\bullet k}\Omega(\mathcal{G}^z_{n-1})\right| = \left| \Omega^{1,k-1}(\mathcal{G}^z_{n-2})\right| = C_{n,n+2-k}$  
    \end{enumerate}
\end {thm}

{\color{red}
}

\begin {proof}
    We will exhibit a bijection $\mathrm{Cat}_n \to \Omega(\mathcal{G}^z_{n-1})$, establishing part $(a)$. From the definition of this map,
    part $(b)$ will be clear.

    The map will be defined in much the same way as the bijection discussed in Section \ref{sec:bijection}. 
    We will exhibit bijections realizing the recurrence relation in Lemma \ref{lem:ballot}(b) for both the set of $132$-avoiding permutations
    and the set of mixed dimer covers of zig-zag snake graphs, and the bijection claimed in the statement of the theorem will be the unique one satsifying two conditions:
    (i) the bijection $F$ between permutations and mixed dimer covers commutes with the combinatorial recurrence bijections, and
    (ii) under the bijection, the first entry in the permutation's Lehmer code equals the multiplicity of the final edge in the mixed dimer cover.

    First, we discuss the recurrence for permutations. As in Section \ref{sec:bijection}, we start by removing the first entry of the permutation,
    and then standardizing. That is, $\sigma = (\sigma_1, \sigma_2, \dots, \sigma_n) \mapsto \bar{\sigma} = (\bar{\sigma}_2, \dots, \bar{\sigma}_n)$,
    where $\bar{\sigma}_i = \sigma_i$ if $\sigma_i < \sigma_1$, and $\bar{\sigma}_i = \sigma_i-1$ if $\sigma_i > \sigma_1$. Note that $\bar{\sigma}$
    will still be $132$-avoiding, since the relative inequalities between entries have not changed. Since $\sigma$ was $132$-avoiding, we must have $\sigma_2 \leq \sigma_1+1$, and so $\bar{\sigma}_2 \leq \sigma_1$. This gives the bijection
    realizing the recurrence from Lemma \ref{lem:ballot}(b). Note that on the level of Lehmer codes, we simply have $(L_1,\dots,L_n) \mapsto (L_2,\dots,L_n)$.

    Now we will give a realization of the same recurrence for mixed dimer covers. 
    Let $\mathcal{Z}_{n,k} = {}^{\bullet,n+2-k}\Omega(\mathcal{G}^z_{n-1})$ be the set of mixed dimer covers of $\mathcal{G}^z_{n-1}$ where the multiplicity
    of the final edge is $k-1$. We want to define a bijection $\mathcal{Z}_{n,k} \to \bigcup_{i \leq k} \mathcal{Z}_{n-1,i}$.
    Let $\ell_1,\ell_2,\dots$ be the multiplicities of the boundary edges, as indicated in Figure \ref{fig:ballot_recurrence}.
    More specifically $\ell_1$ is the multiplicity of the final edge, and each $\ell_i$ labels the final edge of one of the subgraphs $\mathcal{G}^z_{n-i}$.

    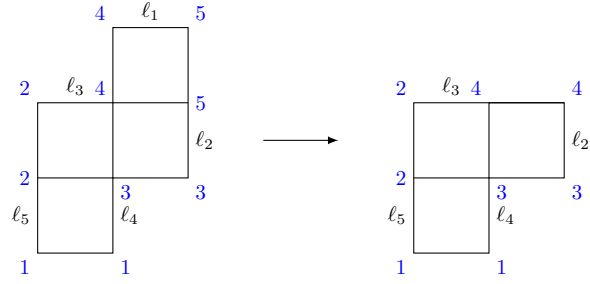
\begin {figure}
    \centering
    \begin {tikzpicture}[scale=1.0, every node/.style={scale=0.8}]
        \draw (0,0) --++ (1,0) --++ (0,1) --++ (1,0) --++ (0,2) --++ (-1,0) --++ (0,-1) --++ (-1,0) -- cycle;
        \draw (0,1) --++ (1,0) --++ (0,1) --++ (1,0);

        \draw[blue] (2,3) node[above right] {$5$};
        \draw[blue] (2,2) node[right]       {$5$};
        \draw[blue] (1,3) node[above left]  {$4$};
        \draw[blue] (1,2) node[above left]  {$4$};
        \draw[blue] (2,1) node[below right] {$3$};
        \draw[blue] (1,1) node[below right] {$3$};
        \draw[blue] (0,2) node[above left]  {$2$};
        \draw[blue] (0,1) node[left]        {$2$};
        \draw[blue] (0,0) node[below left]  {$1$};
        \draw[blue] (1,0) node[below right] {$1$};

        \draw (1.5,3) node[above] {$\ell_1$};
        \draw (2,1.5) node[right] {$\ell_2$};
        \draw (0.5,2) node[above] {$\ell_3$};
        \draw (1,0.5) node[right] {$\ell_4$};
        \draw (0,0.5) node[left]  {$\ell_5$};

        \draw[-latex] (3,1.5) -- (4,1.5);

        \begin {scope}[shift={(5,0)}]
            \draw (0,0) --++ (1,0) --++ (0,1) --++ (1,0) --++ (0,1) --++ (-2,0) -- cycle;
            \draw (0,1) --++ (1,0) --++ (0,1) --++ (1,0);
    
            \draw[blue] (2,2) node[above right] {$4$};
            \draw[blue] (1,2) node[above left]  {$4$};
            \draw[blue] (2,1) node[below right] {$3$};
            \draw[blue] (1,1) node[below right] {$3$};
            \draw[blue] (0,2) node[above left]  {$2$};
            \draw[blue] (0,1) node[left]        {$2$};
            \draw[blue] (0,0) node[below left]  {$1$};
            \draw[blue] (1,0) node[below right] {$1$};
    
            \draw (2,1.5) node[right] {$\ell_2$};
            \draw (0.5,2) node[above] {$\ell_3$};
            \draw (1,0.5) node[right] {$\ell_4$};
            \draw (0,0.5) node[left]  {$\ell_5$};
        \end {scope}
    \end {tikzpicture}
    \caption {The bijection $\mathrm{Cat}_{n,k} \to \bigcup_{i \leq k} \mathrm{Cat}_{n-1,i}$, illustrated for $n=5$.
    Multiplicities of unmarked edges are determined by the $\ell_i$.}
    \label {fig:ballot_recurrence}
    \end {figure}

    The bijection maps a mixed dimer cover of $\mathcal{G}^z_{n-1}$
    to the corresponding mixed dimer cover of $\mathcal{G}^z_{n-2}$ which is identical to the original one except that the edge separating the last and second-to-last
    tiles of $\mathcal{G}^z_{n-1}$ (which is now one of the non-final boundary edges of $\mathcal{G}^z_{n-2}$) has multiplicity $n-1-\ell_2$. Since the $\ell_2$ edge and
    one with multiplicity $n-\ell_1$ in $\mathcal{G}^z_{n-1}$ are adjacent to the same vertex with label $n$, this implies that $\ell_2 \leq \ell_1$. We therefore
    map to some $\mathcal{Z}_{n-1,i}$ with $i \leq k$.

    Comparing the two recurrences, and requiring that our bijection commute with them, implies that $\ell_i = L_i$ for all $i$.
    Knowing the multiplicities of the $\ell_i$ edges uniquely determines the rest of the mixed dimer cover, and so we are done.
\end {proof}

\begin {figure}[h]
\centering
\begin {tikzpicture}[scale=1.2]
\draw (0,0) --++ (1,0) --++ (0,1) --++ (1,0) --++ (0,1) --++ (1,0) --++ (0,1) --++ (1,0) --++ (0,1);
\draw (0,0) --++ (0,2) --++ (1,0) --++ (0,1) --++ (1,0) --++ (0,1) --++ (2,0);
\draw (0,1) --++ (1,0) --++ (0,1) --++ (1,0) --++ (0,1) --++ (1,0) --++ (0,1);

\draw (4,3.5) node[right] {$L_1$};
\draw (3,2.5) node[right] {$L_3$};
\draw (2,1.5) node[right] {$L_5$};
\draw (1,0.5) node[right] {$L_7$};

\draw (2.5,4) node[above] {$L_2$};
\draw (1.5,3) node[above] {$L_4$};
\draw (0.5,2) node[above] {$L_6$};

\draw (3,3.5) node[right] {$\Delta_1$};
\draw (2,2.5) node[right] {$\Delta_3$};
\draw (1,1.5) node[right] {$\Delta_5$};

\draw (2.5,3) node[above] {$\Delta_2$};
\draw (1.5,2) node[above] {$\Delta_4$};
\draw (0.5,1) node[above] {$\Delta_6$};
\end {tikzpicture}
\caption {The bijection $\mathrm{Cat}_8 \to \Omega(\mathcal{G}^z_7)$, with edge multiplicities given in terms of the Lehmer code.
The multiplicities of the internal edges are given by $\Delta_i := L_i - L_{i+1}$.}
\label {fig:catalan_bijection}
\end {figure}
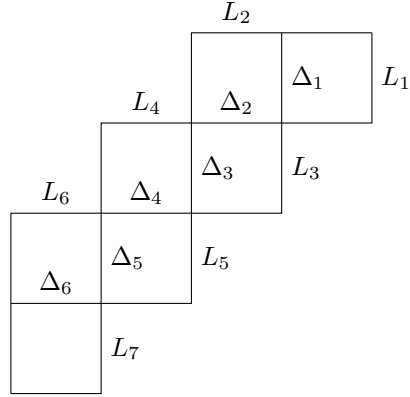

An example of the bijection constructed in the proof is illustrated in Figure \ref{fig:catalan_bijection}.

A matrix product formula similar to Theorem \ref{thm:matrix_formula} can be given in this case as well. Recall the matrices $R_{a,b}$ from before.
Let $W$ be the $(a+1) \times (a+1)$ matrix with $W_{ij} = \delta_{i,a+2-i}$. It is the anti-diagonal matrix with $1$'s along the anti-diagonal.
We define $L_{a,b} := W R_{a,b}$, obtained by reversing the order of the rows of $R_{a,b}$.

\begin {thm} \label{thm:matrix_formula_zigzag}
    For an integer sequence $m = (m_0,m_1,\dots,m_n)$, define 
    the matrix $X_m$ by
    \[ X_m := R_{m_0,m_0}R_{m_0,m_1}L_{m_1,m_2}L_{m_2,m_3} \cdots L_{m_{n-1},m_n} \]
    Then the $(i,j)$-entry of $X_m$ is equal to $\left| \Omega^{ij}_{\n_{m}}(\mathcal{G}^z_n) \right|$. 
    In particular, the $(1,1)$-entry is $\left| \Omega_{\n_m}(\mathcal{G}^z_n) \right|$.
\end {thm}
\begin {proof}
    Note that for a single square, there is no difference between a straight and a zig-zag snake graph. So the base case is the same
    as the matrix formula for Theorem \ref{thm:matrix_formula}. We must simply check that after the first square, the recurrence is
    given by right-multiplication by an $L$-matrix (rather than an $R$-matrix).

    Note that $(L_{m_n,m_{n+1}})_{ij} = 1$ if $j-i \leq m_{n+1}-m_n$. Therefore the entries of the matrix product are given by
    \[ (X L_{m_n,m_{n+1}})_{ij} = \sum_k X_{ik} (L_{m_n,m_{n+1}})_{kj} = \sum_{k \geq j-(m_n-m_{n+1})} X_{ik} \]
    Therefore the thoerem is equivalent to showing that $\left| \Omega_{\n_m}(\mathcal{G}^z_n) \right|$ satisfy the following recurrence:
    \[ \left| \Omega_{\n_m}^{ij}(\mathcal{G}^z_{n+1}) \right| = \sum_{k \geq j-(m_{n+1}-m_n)} \left| \Omega_{\n_m}^{ik}(\mathcal{G}^z_n) \right| \]
    This recurrence follows from a similar bijection to the one pictured in Figure \ref{fig:ballot_recurrence}. 
    Suppose we start with an element of $\Omega_{\n_m}^{ij}(\mathcal{G}^z_{n+1})$, with multiplicity $\ell_1$ on the final edge (as in the figure).
    By removing the three edges of the last tile which do not border the previous tile, we obtain an element of $\Omega_{\n_m}^{i,m_n+1-\ell_1}$.
    Since $\ell_1 \leq m_{n+1}+1-j$, this means $m_n+1-\ell_1 \geq j-(m_{n+1}-m_n)$, as desired.
\end {proof}

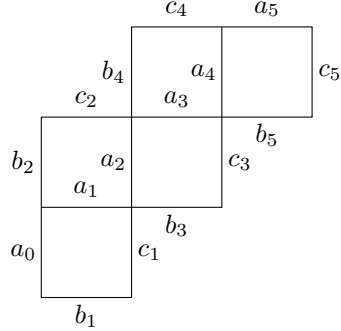
\begin {figure}
\centering
\begin {tikzpicture}[scale=1.2]
    \draw (0,0) -- (1,0) -- (1,1) -- (2,1) -- (2,2) -- (3,2) -- (3,3) -- (1,3) -- (1,2) -- (0,2) -- cycle;
    \draw (0,1) -- (1,1) -- (1,2) -- (2,2) -- (2,3);

    \draw (-0.2,0.5) node{$a_0$};
    \draw (0.5,1.2) node{$a_1$};
    \draw (0.8,1.5) node{$a_2$};
    \draw (1.5,2.2) node{$a_3$};
    \draw (1.8,2.5) node{$a_4$};
    \draw (2.5,3.2) node{$a_5$};

    \draw (0.5,-0.2) node{$b_1$};
    \draw (1.2,0.5) node{$c_1$};
    \draw (1.5,0.8) node{$b_3$};
    \draw (2.2,1.5) node{$c_3$};
    \draw (2.5,1.8) node{$b_5$};
    \draw (3.2,2.5) node{$c_5$};

    \draw (-0.2,1.5) node{$b_2$};
    \draw (0.5,2.2) node{$c_2$};
    \draw (0.8,2.5) node{$b_4$};
    \draw (1.5,3.2) node{$c_4$};
\end {tikzpicture}
\caption {Edge weights used in the weighted version of the matrix product.}
\label {fig:zigzag_edge_weights}
\end {figure}

We can also upgrade this to a weighted version. That is, if the edges of $\mathcal{G}^z_{n-1}$ are given weights, we can modify the $L_{a,b}$ matrices
so that each entry of the matrix product is the weighted sum over the associated set of mixed dimer covers. Define $W_{k}(t)$ as the $(k+1) \times (k+1)$
anti-diagonal matrix with $(i,k+2-i)$-entry $t^k$, and all others zero. Also recall the $U$ and $T$ matrices from section \ref{sec:straight_snakes}.
We assign edge weights $a_i,b_i,c_i$ as in Figure \ref{fig:zigzag_edge_weights}. Then we have the following. Again, we omit the proof, as it is
a straightforward generalization of the calculation from Theorem \ref{thm:matrix_formula_zigzag}.

\begin {thm}
    Let $m = (m_0,\dots,m_n)$ be a sequence of vertex multiplicities for $\mathcal{G}^z_n$, and let $a_i,b_i,c_i$ be edge weights
    as in Figure \ref{fig:zigzag_edge_weights}. Define the matrix
    \[ X =  U_{m_0}(b_1) T_{m_0,m_1}(a_0c_1) U_{m_1}(a_1) \prod_{k=2}^n W_{m_{k-1}}(b_k) T_{m_{k-1},m_k}(b_k^{-1}c_k) U_{m_k}(a_k) \]
    where the product is taken in left-to-right order. Then $X_{ij}$ is the weighted sum of mixed dimer covers in $\Omega^{ij}(\mathcal{G}^z_n)$.
\end {thm}
    
In particular, when all edge weights are $1$, then $T_{m_{k-1},m_k}U_{m_k} = R_{m_{k-1},m_k}$
and $W_{m_{k-1}}T_{m_{k-1},m_k}S_{m_k} = L_{m_{k-1},m_k}$, recovering the formula from Theorem \ref{thm:matrix_formula_zigzag}.

\begin {ex}
    For the standard labeling $m=(1,2,3,4)$ on $\mathcal{G}^z_{3}$, we get
    \[ X = U_1(b_1)T_{1,2}(a_0c_1) U_2(a_1) W_2(b_2)T_{2,3}(b_2^{-1}c_2)U_3(a_2) W_3(b_3)T_{3,4}(b_3^{-1}c_3)U_4(a_3) \]
    The $(1,1)$-entry of this matrix is
    \begin {align*}
        X_{11} &= a_3^4b_1b_2^2b_3^3 + a_2a_3^3b_1b_2^2b_3^2c_3 + a_1a_3^3b_1b_2b_3^2c_2c_3 + a_2^2a_3^2b_1b_2^2b_3c_3^2 + a_0a_3^3b_2b_3^2c_1c_2c_3 \\
        &\phantom{=} + a_1a_2a_3^2b_1b_2b_3c_2c_3^2 + a_2^3a_3b_1b_2^2c_3^3 + a_0a_2a_3^2b_2b_3c_1c_2c_3^2 + a_1^2a_3^2b_1b_3c_2^2c_3^2 + a_1a_2^2a_3b_1b_2c_2c_3^2 \\
        &\phantom{=} + a_0a_1a_3^2b_3c_1c_2^2c_3^2 + a_0a_2^2a_3b_2c_1c_2c_3^3 + a_1^2a_2a_3b_1c_2^2c_3^3 + a_0a_1a_2a_3c_1c_2^2c_3^3
    \end {align*}
    These are precisely the weights of the 14 mixed dimer configurations pictured in Figure \ref{fig:catalan_posets}.
\end {ex}

\section {Distributive Lattice Structures and $q$-Analogs} \label{sec:lattices}

It is known that the set of perfect matchings (and more generally the set of $\mathbf{d}$-factors) of a planar bipartite graph has a partial order which
makes it a distributive lattice \cite{propp_02}. This partial order generalizes to the current setting of mixed dimer covers,
and we will now describe this partial order. 

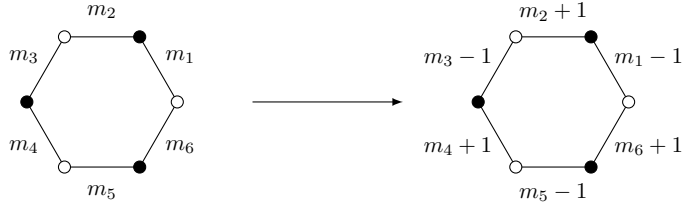
\begin {figure}
\centering
\begin {tikzpicture}[every node/.style={scale=0.85}]
\draw (0:1) -- (60:1) -- (120:1) -- (180:1) -- (240:1) -- (300:1) -- cycle;
\foreach \t in {0,120,240} {
    \draw[fill=white] (\t:1) circle (0.08);
    \draw[fill=black] (\t+60:1) circle (0.08);
}

\draw (30:1.2) node{$m_1$};
\draw (90:1.2) node{$m_2$};
\draw (150:1.2) node{$m_3$};
\draw (210:1.2) node{$m_4$};
\draw (270:1.2) node{$m_5$};
\draw (330:1.2) node{$m_6$};

\draw[-latex] (2,0) -- (4,0);

\begin {scope}[shift={(6,0)}]
\draw (0:1) -- (60:1) -- (120:1) -- (180:1) -- (240:1) -- (300:1) -- cycle;
\foreach \t in {0,120,240} {
    \draw[fill=white] (\t:1) circle (0.08);
    \draw[fill=black] (\t+60:1) circle (0.08);
}

\draw (25:1.4) node{$m_1-1$};
\draw (90:1.2) node{$m_2+1$};
\draw (155:1.4) node{$m_3-1$};
\draw (205:1.4) node{$m_4+1$};
\draw (270:1.2) node{$m_5-1$};
\draw (335:1.4) node{$m_6+1$};
\end {scope}
\end {tikzpicture}
\caption {A positive face twist of a mixed dimer cover, with edge multiplicities $m_i$.}
\label{fig:face_twist}
\end {figure}

Let $\mathcal{G}$ be a planar bipartite graph, with vertices colored black and white.
For a given face of the graph, if we traverse the edges around the face in the counter-clockwise direction, we will alternate going black-to-white
and white-to-black. Call the former edges \emph{even} and the latter \emph{odd}. For a given mixed dimer configuration $\omega \in \Omega_\n(\mathcal{G})$, suppose that the multiplicities of all odd edges around a face are positive. We can decrease all these multiplicities by 1, and simultaneously increase the multiplicities
of the even edges around this face, producing a new mixed dimer cover. We call such move a \emph{(positive) face twist} (and its inverse a \emph{negative face twist}).
See Figure \ref{fig:face_twist} for an illustration.
The partial order on $\Omega_\n(\mathcal{G})$ has the positive face twists as the covering relations.

\begin {rmk}
    Note that the definition depends on the black/white coloring. Interchanging the black/white colors turns the poset into its dual
    (i.e. the Hasse diagram is turned upside down). To remove this ambiguity, we will adopt the following conventions.
    For the straight snake $\mathcal{G}^s_n$, we always assume the vertices are
    colored so that the bottom-right vertex is black. For the zig-zag $\mathcal{G}^z_n$,
    if $n$ is odd (so the final edge is the right edge of the last square), we color the bottom-right vertex of the last square black,
    and if $n$ is even (so the final edge is the top edge of the last square), we color the top-right vertex of the last square black.
\end {rmk}

\begin {thm} \label{thm:lattice}
    The set $\Omega_\n(\mathcal{G})$ is a distributive lattice under the face-twist order.
\end {thm}
%
%

We will delay the proof of this until later (the end of Section \ref{sec:dimers_paths}), after some other helpful results and observations.

The rank function of this poset is defined as follows. For an element $\omega \in \Omega_\n(\mathcal{G})$, its rank $\mathrm{rk}(\omega)$ is the
minimal number of face twists to reach $\omega$ from the minimal element. The rank generating function is denoted 
\[ \Omega_\n(\mathcal{G},q) := \sum_{\omega \in \Omega_\n(\mathcal{G})} q^{\mathrm{rk}(\omega)} \]

There are natural $q$-analogs of the Euler and Entringer numbers, which count the inversions in
alternating permutations. That is,
\[ E_{n,k}(q) := \sum_{\sigma \in \mathrm{Alt}_{n,k}} q^{\mathrm{inv}(\sigma)} \]


We will show that the rank generating function for the poset of mixed dimer covers on the straight snake graph $\mathcal{G}^s_{n-2}$
is (a multiple of) the $q$-Euler number $E_n(q)$. This will follow from establishing a poset isomorphism with a certain partial order on alternating permutations.
We will describe this partial order in the next section.

\subsection {The Middle Order on Permutations}

Consider the following inversion code for permutations:
\[ I(\sigma) = (x_1,\dots,x_n), \quad \text{ where } x_i := \#\{j<i ~|~ \sigma^{-1}(j) > \sigma^{-1}(i)\}. \]
This encoding gives a bijection $S_n \to [0,0] \times [0,1] \times \cdots \times [0,n-1]$. 
In \cite{bft}, the \emph{middle order} on $S_n$ was defined as the
partial order (induced from this bijection) coming from the standard order on the product of chains.
That is, if $I(\sigma) = (x_1,\dots,x_n)$ and $I(\tau) = (y_1,\dots,y_n)$,
we say that $\sigma \leq \tau$ if $x_i \leq y_i$ for all $i$. Moreover, $\sigma < \tau$ is a covering relation if $I(y)$ and $I(x)$ 
differ in only one entry, with $y_i = x_i+1$.

Several facts about the middle order are established in \cite{bft}, including the fact that it makes $S_n$ into a distributive lattice, and that it
is a refinement of the right weak order and a coarsening of the Bruhat order (hence the name \emph{middle order}). 
We will henceforth refer to it as the \emph{right middle order}, to emphasize that it is a refinement of the right weak order (rather than the left).

Now, let us make a very similar definition, but using the Lehmer code rather than $I(\sigma)$.
As before, we define a partial order on $S_n$ whose cover relations are given by $\sigma < \tau$ when $L(\tau)$ is obtained from $L(\sigma)$
by incrementing a single entry. We will call this the \emph{left middle order} on $S_n$. Shortly, we will give a justification of this name.

Note that since both left and right middle orders come from the corresponding partial orders on the product of chains, it is clear
that they are isomorphic as posets, and the isomorphism $S_n \to S_n$ is given by $I^{-1} \circ \mathrm{rev} \circ L$, where $\mathrm{rev}$
is the reversal of sequences (i.e. $\mathrm{rev}(a_1,a_2,\dots,a_n) = (a_n,a_{n-1},\dots,a_1)$). We will now describe more explicitly 
what this mapping looks like, and give a justification of the name ``left middle order''.

\begin {prop}
    Let $L$ and $I$ be the Lehmer code and inversion code maps defined above. The
    composition $I^{-1} \circ \mathrm{rev} \circ L \colon S_n \to S_n$ is given by $\sigma \mapsto w_0 \sigma^{-1} w_0$, where
    $w_0$ is the top element in the Bruhat order ($w_0(k) = n+1-k$).
\end {prop}
\begin {proof}
    The claim is equivalent to saying that $L(\sigma)$ and $I(w_0\sigma^{-1}w_0)$ are reversals of each other.
    In other words, if $L(\sigma) = (L_1,\dots,L_n)$ and $I(w_0\sigma^{-1}w_0) = (x_1,\dots,x_n)$, then we need to show
    that $L_i = x_{n+1-i}$. 

    We will do this by observing how inversions are transformed under first $\sigma \mapsto \sigma^{-1}$, and then $\tau \mapsto w_0 \tau w_0$.
    Suppose $\sigma$ has an inversion in positions $i<j$ with values $\sigma_i > \sigma_j$. Then $\sigma^{-1}$ has a corresponding
    inversion in positions $\sigma_j < \sigma_i$ with values $\sigma^{-1}(\sigma_j) = j > i = \sigma^{-1}(\sigma_i)$.
    
    Next, we examine conjugation by $w_0$. Note that multiplying by $w_0$ on the right reverses the one-line notation of a permutation,
    while multiplying by $w_0$ on the left takes complementary values (that is, occurrences of $k$ in the one-line notation become $n+1-k$).
    Conjugation by $w_0$ therefore does both. One sees that an inversion in a permutation $\tau$ with positions $i<j$ and values $\tau_i > \tau_j$
    corresponds to an inversion in $w_0 \tau w_0$ with positions $n+1-j < n+1-i$ and values $n+1-\tau_j > n+1-\tau_i$.

    Combining the observations in the previous two paragraphs, we see that an inversion of $\sigma$ in positions $i<j$ 
    and values $\sigma_i > \sigma_j$ corresponds to the inversion
    of $w_0 \sigma^{-1} w_0$ in positions $n+1-\sigma_i < n+1-\sigma_j$ and values $n+1-i > n+1-j$. To summarize, the number of inversions of $\sigma$
    whose left \emph{position} is $i$ is equal to the number of inversions in $w_0 \sigma^{-1} w_0$ whose smaller \emph{value} is $n+1-i$.
    This is equivalent to the desired statement $L_i = x_{n+1-i}$.
\end {proof}

The map $\sigma \mapsto w_0 \sigma^{-1} w_0$ is a poset isomorphism between the left and right weak orders, and an automorphism
of the Bruhat order (see e.g. \cite{bb_05} chapters 2 and 3). 
From this, we immediately get the following, which justifies our name of ``left middle order''.

\begin {cor}
    The left middle order is a refinement of the left weak order, and a coarsening of the Bruhat order.
\end {cor}

\begin {rmk} \label{rmk:left-right-transfer}
    If $X \subseteq S_n$ is any subset of permutations, the left middle order on $X$ is isomorphic to the 
    right middle order on the set $w_0 X^{-1} w_0$. As an example, the map $\sigma \mapsto w_0 \sigma^{-1} w_0$ interchanges
    the sets of $132$-avoiding and $213$-avoiding permutations, and hence the statements in remarks 1.10 and 1.11 in \cite{bft} are interchanged if we
    consider the left (rather than right) middle order. Specifically, the remarks from \cite{bft} imply that the left middle order on $132$-avoiding permutations
    is the same as the Bruhat order, and the left middle order on $213$-avoiding permutations is the same as the left weak order.
\end {rmk}

\subsection {Rank Functions for Mixed Dimer Lattices}

We will now show that the face twist partial order on mixed dimer covers is isomorphic to the restriction of the middle order on the appropriate
subset of permutations.

\begin {thm} \label{thm:rank_equals_q_euler}
    Let $m = (1,2,\dots,n-1)$ be the standard vertex labeling of $\mathcal{G}^s_{n-2}$. 
    \begin {enumerate}
        \item[(a)] The bijection $\varphi \colon \Omega(\mathcal{G}^s_{n-2}) \to \mathrm{Alt}_n$ from Theorem \ref{thm:euler_bijection}
                   is a poset isomorphism, where the partial order on $\mathrm{Alt}_n$ is the restriction of the left middle order. 
        \item[(b)] The rank generating function for the face twist order on $\Omega(\mathcal{G}^s_{n-2})$ is equal to the corresponding $q$-Euler number
                  (up to a monomial factor). More precisely,
                  \[ E_n(q) = \sum_{\sigma \in \mathrm{Alt}_n} q^{\mathrm{inv}(\sigma)} = q^{\lfloor n/2 \rfloor} \, \Omega(\mathcal{G}^s_{n-2},q) \]
    \end {enumerate}
\end {thm}
\begin {proof}
    Order the square faces of $\mathcal{G}^s_{n-2}$ from right-to-left. 
    From the bijection in Theorem \ref{thm:euler_bijection}, and the labels illustrated in Figure \ref{fig:lehmer}
    giving the edge multipicities in terms of the Lehmer code $L(\sigma)$, it is easy to see that a positive face twist at 
    the $k^\mathrm{th}$ face increases $L_{k}$ by one (and all other $L_i$ are unaffected).  
    This is a covering relation in the middle order on $S_n$.

    Since the bijection in part $(a)$ sends covering relations to covering relations, it is clear that $\Omega(\mathcal{G}^s_{n-2},q)$ and $E_n(q)$ 
    are the same up to some factor of $q^N$.
    To see that $N = \lfloor \frac{n}{2} \rfloor$, we just note that the minimal alternating permutation
    is $\sigma = 214365\cdots = \prod_{i \geq 1} (2i-1,2i)$, which has $\lfloor \frac{n}{2} \rfloor$ inversions.
\end {proof}

\begin {rmk}
    By restricting to $\sigma \in \mathrm{Alt}_{n,k}$, the theorem tells us that the $q$-analog of the Entringer number $E_{n,k}(q)$ 
    is the rank function of the induced sub-poset consisting of mixed dimer covers with the right-most vertical edge covered by $k-1$ dimers.
\end {rmk}

\begin {rmk}
    Let $\mathrm{Alt}_n^{-1} = \{\sigma \in S_n ~|~ \sigma^{-1} \in \mathrm{Alt}_n\}$. These are permutations such that in one-line notation,
    $2$ appears to the left of $1$, $3$ to the right of $2$, $4$ to the left of $3$, $5$ to the right of $4$, etc. 
    Since $w_0 \mathrm{Alt}_n^{-1}w_0 = \mathrm{Alt}_n^{-1}$ when $n$ is even, we have by Remark \ref{rmk:left-right-transfer} 
    that $\Omega(\mathcal{G}^s_{n-2})$ is also isomorphic to the right middle order on $\mathrm{Alt}_n^{-1}$.
    On the other hand, when $n$ is odd, $w_0 \mathrm{Alt}_n^{-1}w_0 = \mathrm{RAlt}_n^{-1}$, where $\mathrm{RAlt}_n$ is the
    set of reverse alternating permutations (see Remark \ref{rmk:ralt}),
    and so $\Omega(\mathcal{G}^s_{n-2})$ is isomorphic to
    the right middle order on $\mathrm{RAlt}_n^{-1}$.
\end {rmk}

Since the mixed dimer posets are always distributive lattices, we obtain the following fact about the
middle order on alternating permutations.

\begin {cor}
    The restriction of the left middle order on $\mathrm{Alt}_n$ is a distributive lattice.
\end {cor}

\begin {figure}[h!]
\centering
\begin {tikzpicture}[scale=0.5]
    \foreach \x/\y in {0/0, 0/1,-1/1,-2/2,-1/2,0/2,-2/3,-1/3,0/3,1/3,-1/4,0/4,1/4,0/5,-1/5,0/6} {
        \draw (0+4*\x,0+3*\y) -- (3+4*\x,0+3*\y) -- (3+4*\x,1+3*\y) -- (0+4*\x,1+3*\y) -- cycle;
        \draw (1+4*\x,0+3*\y) -- (1+4*\x,1+3*\y);
        \draw (2+4*\x,0+3*\y) -- (2+4*\x,1+3*\y);
    }

    \foreach \x/\y in {0/0,-1/1,-2/2,0/2,-1/3,1/3,0/4,-1/4,0/5} {
        \draw (4*\x+1.5,3*\y+1.5) -- (4*\x+1.5,3*\y+2.5);
    }
    \foreach \x/\y in {0/0,0/1,-1/1,-1/2,0/2,0/3,1/3,1/4} {
        \draw (4*\x+0.5,3*\y+1.5) -- (4*\x-1.5,3*\y+2.5);
    }
    \foreach \x/\y in {-1/1,-2/2,-1/2,0/2,0/3,-1/3,-2/3,-1/4,-1/5} {
        \draw (4*\x+2.5,3*\y+1.5) -- (4*\x+4.5,3*\y+2.5);
    }

    \begin {scope}[shift={(0,0)}]
        \draw[blue, line width=1.2] (0,0) -- (1,0) -- (1,1) -- (0,1);
        \draw[blue, line width=1.2] (2,0) -- (3,0) -- (3,1) -- (2,1);
        \draw[blue, line width=1.2] (2,-0.12) -- (3,-0.12);
        \draw[blue, line width=1.2] (2,-0.24) -- (3,-0.24);
        \draw[blue, line width=1.2] (2,1.12) -- (3,1.12);
        \draw[blue, line width=1.2] (2,1.24) --  (3,1.24);
    \end {scope}

    \begin {scope}[shift={(0,3)}]
        \draw[blue, line width=1.2] (0,0) -- (0,1);
        \draw[blue, line width=1.2] (1,0) -- (1,1);
        \draw[blue, line width=1.2] (1.12,0) -- (1.12,1);
        \draw[blue, line width=1.2] (2,0) -- (3,0) -- (3,1) -- (2,1);
        \draw[blue, line width=1.2] (2,-0.12) -- (3,-0.12);
        \draw[blue, line width=1.2] (2,-0.24) -- (3,-0.24);
        \draw[blue, line width=1.2] (2,1.12) -- (3,1.12);
        \draw[blue, line width=1.2] (2,1.24) --  (3,1.24);
    \end {scope}

    \begin {scope}[shift={(-4,3)}]
        \draw[blue, line width=1.2] (0,0) -- (1,0) -- (1,1) -- (0,1);
        \draw[blue, line width=1.2] (2,0) -- (3,0) -- (3,1) -- (2,1) -- cycle;
        \draw[blue, line width=1.2] (2,-0.12) -- (3,-0.12);
        \draw[blue, line width=1.2] (2,1.12) -- (3,1.12);
        \draw[blue, line width=1.2] (3.12,0) -- (3.12,1);
    \end {scope}

    \begin {scope}[shift={(-8,6)}]
        \draw[blue, line width=1.2] (0,0) -- (3,0) -- (3,1) -- (0,1);
        \draw[blue, line width=1.2] (2,-0.12) -- (3,-0.12);
        \draw[blue, line width=1.2] (2,1.12) -- (3,1.12);
        \draw[blue, line width=1.2] (3.12,0) -- (3.12,1);
    \end {scope}

    \begin {scope}[shift={(-4,6)}]
        \draw[blue, line width=1.2] (0,0) -- (0,1);
        \draw[blue, line width=1.2] (1,0) -- (1,1);
        \draw[blue, line width=1.2] (1.12,0) -- (1.12,1);
        \draw[blue, line width=1.2] (2,0) -- (3,0) -- (3,1) -- (2,1) -- cycle;
        \draw[blue, line width=1.2] (2,-0.12) -- (3,-0.12);
        \draw[blue, line width=1.2] (2,1.12) -- (3,1.12);
        \draw[blue, line width=1.2] (3.12,0) -- (3.12,1);
    \end {scope}

    \begin {scope}[shift={(0,6)}]
        \draw[blue, line width=1.2] (0,0) -- (1,0) -- (1,1) -- (0,1);
        \draw[blue, line width=1.2] (2,0) -- (3,0) -- (3,1) -- (2,1) -- cycle;
        \draw[blue, line width=1.2] (1.88,0) -- (1.88,1);
        \draw[blue, line width=1.2] (3.12,0) -- (3.12,1);
        \draw[blue, line width=1.2] (3.24,0) -- (3.24,1);
    \end {scope}

    \begin {scope}[shift={(-8,9)}]
        \draw[blue, line width=1.2] (0,0) -- (0,1);
        \draw[blue, line width=1.2] (1,0) -- (3,0) -- (3,1) -- (1,1) -- cycle;
        \draw[blue, line width=1.2] (2,-0.12) -- (3,-0.12);
        \draw[blue, line width=1.2] (2,1.12) -- (3,1.12);
        \draw[blue, line width=1.2] (3.12,0) -- (3.12,1);
    \end {scope}

    \begin {scope}[shift={(-4,9)}]
        \draw[blue, line width=1.2] (0,0) -- (3,0) -- (3,1) -- (0,1);
        \draw[blue, line width=1.2] (2,0) -- (2,1);
        \draw[blue, line width=1.2] (3.12,0) -- (3.12,1);
        \draw[blue, line width=1.2] (3.24,0) -- (3.24,1);
    \end {scope}

    \begin {scope}[shift={(0,9)}]
        \draw[blue, line width=1.2] (0,0) -- (0,1);
        \draw[blue, line width=1.2] (1,0) -- (1,1);
        \draw[blue, line width=1.2] (1.12,0) -- (1.12,1);
        \draw[blue, line width=1.2] (2,0) -- (3,0) -- (3,1) -- (2,1) -- cycle;
        \draw[blue, line width=1.2] (1.88,0) -- (1.88,1);
        \draw[blue, line width=1.2] (3.12,0) -- (3.12,1);
        \draw[blue, line width=1.2] (3.24,0) -- (3.24,1);
    \end {scope}

    \begin {scope}[shift={(4,9)}]
        \draw[blue, line width=1.2] (0,0) -- (1,0) -- (1,1) -- (0,1);
        \draw[blue, line width=1.2] (2,0) -- (2,1);
        \draw[blue, line width=1.2] (2.12,0) -- (2.12,1);
        \draw[blue, line width=1.2] (2.24,0) -- (2.24,1);
        \draw[blue, line width=1.2] (3,0) -- (3,1);
        \draw[blue, line width=1.2] (3.12,0) -- (3.12,1);
        \draw[blue, line width=1.2] (3.24,0) -- (3.24,1);
        \draw[blue, line width=1.2] (3.36,0) -- (3.36,1);
    \end {scope}

    \begin {scope}[shift={(-4,12)}]
        \draw[blue, line width=1.2] (0,0) -- (0,1);
        \draw[blue, line width=1.2] (1,0) -- (3,0) -- (3,1) -- (1,1) -- cycle;
        \draw[blue, line width=1.2] (2,0) -- (2,1);
        \draw[blue, line width=1.2] (3.12,0) -- (3.12,1);
        \draw[blue, line width=1.2] (3.24,0) -- (3.24,1);
    \end {scope}

    \begin {scope}[shift={(0,12)}]
        \draw[blue, line width=1.2] (0,0) -- (2,0) -- (2,1) -- (0,1);
        \draw[blue, line width=1.2] (2.12,0) -- (2.12,1);
        \draw[blue, line width=1.2] (3,0) -- (3,1);
        \draw[blue, line width=1.2] (3.12,0) -- (3.12,1);
        \draw[blue, line width=1.2] (3.24,0) -- (3.24,1);
        \draw[blue, line width=1.2] (3.36,0) -- (3.36,1);
    \end {scope}

    \begin {scope}[shift={(4,12)}]
        \draw[blue, line width=1.2] (0,0) -- (0,1);
        \draw[blue, line width=1.2] (1,0) -- (1,1);
        \draw[blue, line width=1.2] (1.12,0) -- (1.12,1);
        \draw[blue, line width=1.2] (2,0)    -- (2,1);
        \draw[blue, line width=1.2] (2.12,0) -- (2.12,1);
        \draw[blue, line width=1.2] (2.24,0) -- (2.24,1);
        \draw[blue, line width=1.2] (3,0) -- (3,1);
        \draw[blue, line width=1.2] (3.12,0) -- (3.12,1);
        \draw[blue, line width=1.2] (3.24,0) -- (3.24,1);
        \draw[blue, line width=1.2] (3.36,0) -- (3.36,1);
    \end {scope}

    \begin {scope}[shift={(-4,15)}]
        \draw[blue, line width=1.2] (0,0) -- (0,1);
        \draw[blue, line width=1.2] (1,0) -- (3,0) -- (3,1) -- (1,1);
        \draw[blue, line width=1.2] (1,-0.12) -- (2,-0.12);
        \draw[blue, line width=1.2] (1,1.12) -- (2,1.12);
        \draw[blue, line width=1.2] (3.12,0) -- (3.12,1);
        \draw[blue, line width=1.2] (3.24,0) -- (3.24,1);
    \end {scope}

    \begin {scope}[shift={(0,15)}]
        \draw[blue, line width=1.2] (0,0) -- (0,1);
        \draw[blue, line width=1.2] (1,0) -- (2,0) -- (2,1) -- (1,1) -- cycle;
        \draw[blue, line width=1.2] (2.12,0) -- (2.12,1);
        \draw[blue, line width=1.2] (3,0) -- (3,1);
        \draw[blue, line width=1.2] (3.12,0) -- (3.12,1);
        \draw[blue, line width=1.2] (3.24,0) -- (3.24,1);
        \draw[blue, line width=1.2] (3.36,0) -- (3.36,1);
    \end {scope}

    \begin {scope}[shift={(0,18)}]
        \draw[blue, line width=1.2] (0,0) -- (0,1);
        \draw[blue, line width=1.2] (1,0) -- (2,0) -- (2,1) -- (1,1);
        \draw[blue, line width=1.2] (1,-0.12) -- (2,-0.12);
        \draw[blue, line width=1.2] (1,1.12) -- (2,1.12);
        \draw[blue, line width=1.2] (3,0) -- (3,1);
        \draw[blue, line width=1.2] (3.12,0) -- (3.12,1);
        \draw[blue, line width=1.2] (3.24,0) -- (3.24,1);
        \draw[blue, line width=1.2] (3.36,0) -- (3.36,1);
    \end {scope}


    \begin {scope}[shift={(18,0)}]
        \foreach \x/\y in {0/0,-1/1,-2/2,0/2,-1/3,1/3,0/4,-1/4,0/5} {
            \draw (4*\x+1.5,3*\y+1.5) -- (4*\x+1.5,3*\y+2.5);
        }
        \foreach \x/\y in {0/0,0/1,-1/1,-1/2,0/2,0/3,1/3,1/4} {
            \draw (4*\x+0.5,3*\y+1.5) -- (4*\x-1.5,3*\y+2.5);
        }
        \foreach \x/\y in {-1/1,-2/2,-1/2,0/2,0/3,-1/3,-2/3,-1/4,-1/5} {
            \draw (4*\x+2.5,3*\y+1.5) -- (4*\x+4.5,3*\y+2.5);
        }

        \draw (0+1.5,    0.5)     node {$21435$};
        \draw (0+1.5,    1*3+0.5) node {$21534$};
        \draw (-1*4+1.5, 1*3+0.5) node {$31425$};
        \draw (-2*4+1.5, 2*3+0.5) node {$32415$};
        \draw (-1*4+1.5, 2*3+0.5) node {$31524$};
        \draw (0*4+1.5,  2*3+0.5) node {$41325$};
        \draw (-2*4+1.5, 3*3+0.5) node {$32514$};
        \draw (-1*4+1.5, 3*3+0.5) node {$42315$};
        \draw (0*4+1.5,  3*3+0.5) node {$41523$};
        \draw (1*4+1.5,  3*3+0.5) node {$51324$};
        \draw (-1*4+1.5, 4*3+0.5) node {$42513$};
        \draw (0*4+1.5,  4*3+0.5) node {$52314$};
        \draw (1*4+1.5,  4*3+0.5) node {$51423$};
        \draw (-1*4+1.5, 5*3+0.5) node {$43512$};
        \draw (0*4+1.5,  5*3+0.5) node {$52413$};
        \draw (0*4+1.5,  6*3+0.5) node {$53412$};
    \end {scope}
\end {tikzpicture}
\caption{(Left) The face twist partial order on $\Omega(\mathcal{G}^s_3)$; (Right) The left middle order on $\mathrm{Alt}_5$.}
\label{fig:alt5_poset}
\end {figure}
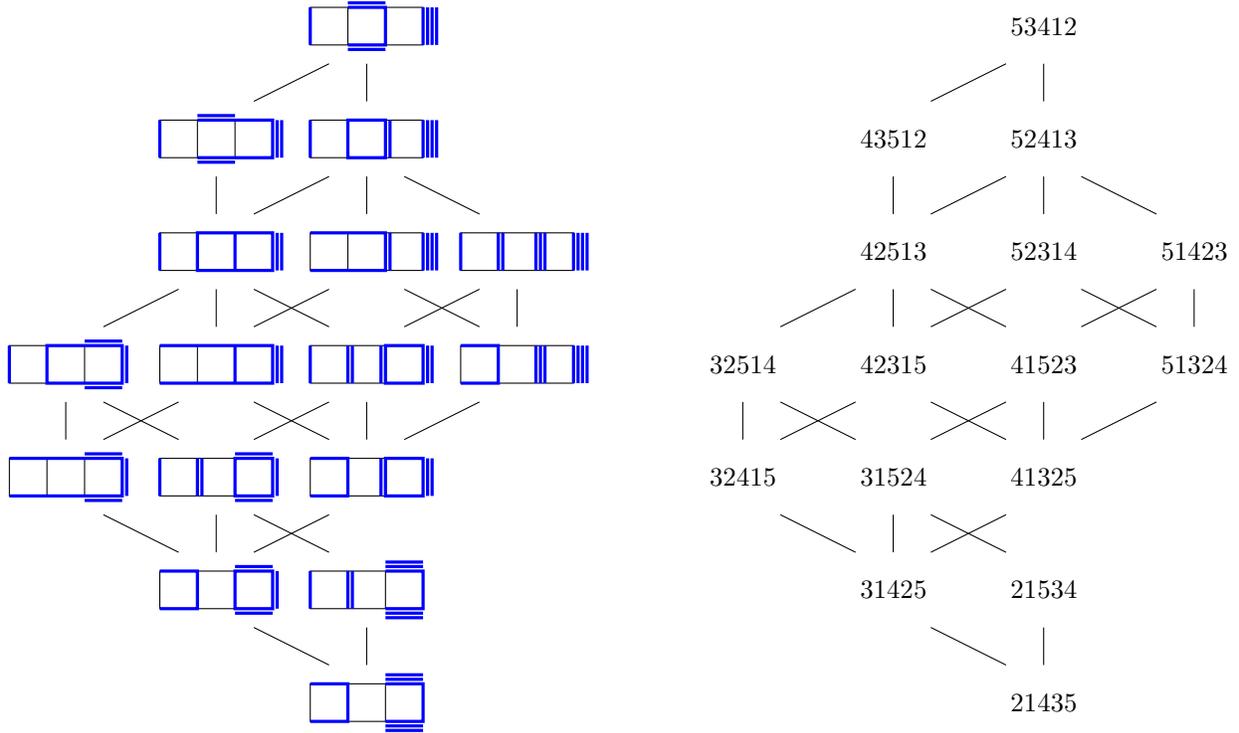

Recall from Theorem \ref{thm:weighted_matrix_product} that given arbitrary edge weights on $\mathcal{G}^s_n$, we can compute
the weighted sum of mixed dimer covers using a matrix product. There is a particular choice of edge weights which will allow us
to easily compute the $q$-Euler and $q$-Entringer numbers. 

Suppose that a face of a planar bipartite graph has $2\ell$ edges $e_1,e_2,\dots,e_{2\ell}$, with edge weights $x_1,\dots,x_{2\ell}$.
Since a face twist decreases the multiplicity of the odd edges by one, and increases the even edges, the weight of the mixed dimer covers
before and after the face twist differ by a factor of $f := \prod_{i=1}^\ell \frac{x_{2i}}{x_{2i-1}}$, which we call the \emph{face weight}.
If one chooses the edge weights so that every face has weight $q$ (which is always possible), then all positive face twists will multiply
the weight of a mixed dimer cover by $q$. The weighted sum of mixed dimer covers will thus be $q^N \Omega_{\n}(\mathcal{G},q)$, where
$q^N$ is the weight of the minimal element using this choice of edge weights.

For the graph $\mathcal{G}^s_n$, one particular choice of edge weights with this property is given as follows. 
Label the horizontal edges along the bottom of $\mathcal{G}^s_n$ alternately $q$, $q^{-1}$, $q$, $q^{-1}$, $\dots$,
and all other edges have weight $1$. Comparing the weight of the minimal mixed dimer cover with the inversion number of the minimal
alternating permutation, one can see that if $Z$ is the weighted sum of mixed dimer covers with this choice of edge weights, then we will have
\[ E_n(q) = q^{\lfloor n^2/4 \rfloor} Z \]

\begin {ex}
    To compute $E_5(q)$ (the rank function of the poset in Figure \ref{fig:alt5_poset}), we use the matrix product
    \[
        \begin{pmatrix} 1&1\\1&0 \end{pmatrix}
        \begin{pmatrix} 1&1&1\\q^{-1}&q^{-1}&0\end{pmatrix}
        \begin{pmatrix} 1&1&1&1\\q&q&q&0\\q^2&q^2&0&0 \end{pmatrix}
        \begin{pmatrix} 1&1&1&1&1\\q^{-1}&q^{-1}&q^{-1}&q^{-1}&0\\q^{-2}&q^{-2}&q^{-2}&0&0\\q^{-3}&q^{-3}&0&0&0 \end{pmatrix}
    \]
    The $(1,1)$-entry of this matrix is $Z = q^{-4} + 2q^{-3} + 3q^{-2} + 4q^{-1} + 3 + 2q + q^2 = q^{-6}E_5(q)$.
\end {ex}

\begin {rmk}
    This method of computing the rank generating function using $q$-deformed versions of the matrices is a generalized
    version (using non-square matrices) of the method used in \cite{bosz_24}, which itself is a generalization
    of the matrix product formulas from \cite{mgo} used to compute $q$-analogs of rational numbers.
\end {rmk}

We now discuss a similar correspondence in the Catalan case. There are several different known $q$-analogs of the Catalan numbers,
but here we focus on the one due originally to Carlitz and Riordan \cite{carlitz_riordan}, which we write as $C_n(q)$. 
It was noted in \cite[Thm 3.1]{sagan_patterns} that $C_n(q)$ counts inversions in $132$-avoiding permutations. That is,
\[ C_n(q) = \sum_{\sigma \in \mathrm{Cat}_n} q^{\mathrm{inv}(\sigma)} \]
Similar to the $q$-Euler numbers, the $q$-Catalan numbers are the rank functions for the partial order on mixed dimer covers.

\begin {thm} \label{thm:rank_equals_q_catalan}
    For the standard vertex labeling of the zigzag snake $\mathcal{G}^z_{n-1}$, we have: 
    \begin {enumerate}
        \item[(a)] The bijection $\varphi \colon \Omega(\mathcal{G}^z_{n-1}) \to \mathrm{Cat}_n$ from Theorem \ref{thm:dimer_catalan_bijection}
                   is a poset isomorphism with the Bruhat order on $\mathrm{Cat}_n$. 
        \item[(b)] The rank generating function for the face twist order on $\Omega(\mathcal{G}^z_{n-1})$ is equal to the $q$-Catalan number $C_n(q)$.
    \end {enumerate}
\end {thm}
\begin {proof}
    The proof is the same as that of Theorem \ref{thm:rank_equals_q_euler}. 
    One easily sees that under the bijection from Theorem \ref{thm:dimer_catalan_bijection}, a face twist corresponds to
    incrementing one entry in the Lehmer code, which is a covering relation in the left middle order. 
    By Remark \ref{rmk:left-right-transfer}, the left middle order on $132$-avoiding permutations
    agrees with the Bruhat order.

    For part (b), we simply note that the identity permutation is 132-avoiding, so $C_n(q)$ and the rank generating function both have constant term 1,
    and so they must agree.
\end {proof}

\begin {figure}
\centering
\begin {tikzpicture}[scale=0.4]
    \foreach \x/\y in {0/0, 0/1, -1/2, 1/2, -2/3, 0/3, 2/3, -1/4, 0/4, 1/4, -1/5, 0/5, 1/5, 0/6} {
        \begin {scope}[shift={(4*\x,4*\y)}]
            \draw (0,0) -- (1,0) -- (1,1) -- (2,1) -- (2,2) -- (0,2) -- cycle;
            \draw (0,1) -- (1,1) -- (1,2);
        \end {scope}
    }

    \foreach \x/\y in {0/0, 0/3, -1/4, 1/4, 0/5} {
        \begin {scope}[shift={(4*\x,4*\y)}]
            \draw (1,2.5) -- (1,3.5);
        \end {scope}
    }

    \foreach \x/\y in {0/1, -1/2, 1/2, -2/3, 0/3, -1/4, 0/4, -1/5} {
        \begin {scope}[shift={(4*\x,4*\y)}]
            \draw (1.5,2.5) -- (3.5,4.5);
        \end {scope}
    }

    \foreach \x/\y in {0/1, -1/2, 1/2, 0/3, 2/3, 0/4, 1/4, 1/5} {
        \begin {scope}[shift={(4*\x,4*\y)}]
            \draw (0.5,2.5) -- (-1.5,4.5);
        \end {scope}
    }

    \begin {scope}[shift={(4*0, 4*0)}]
        \draw[blue, line width=1] (0,0) -- (1,0);
        \draw[blue, line width=1] (0,1) -- (0,2);
        \draw[blue, line width=1] (-0.12,1) -- (-0.12,2);
        \draw[blue, line width=1] (1,1) -- (2,1);
        \draw[blue, line width=1] (1,0.88) -- (2,0.88);
        \draw[blue, line width=1] (1,0.76) -- (2,0.76);
        \draw[blue, line width=1] (1,2) -- (2,2);
        \draw[blue, line width=1] (1,2.12) -- (2,2.12);
        \draw[blue, line width=1] (1,2.24) -- (2,2.24);
        \draw[blue, line width=1] (1,2.36) -- (2,2.36);
    \end {scope}

    \begin {scope}[shift={(4*0, 4*1)}]
        \draw[blue, line width=1] (0,0) -- (1,0);
        \draw[blue, line width=1] (0,1) -- (0,2);
        \draw[blue, line width=1] (-0.12,1) -- (-0.12,2);
        \draw[blue, line width=1] (1,1) -- (2,1) -- (2,2) -- (1,2) -- cycle;
        \draw[blue, line width=1] (1,0.88) -- (2,0.88);
        \draw[blue, line width=1] (1,2.12) -- (2,2.12);
        \draw[blue, line width=1] (1,2.24) -- (2,2.24);
    \end {scope}

    \begin {scope}[shift={(4*-1, 4*2)}]
        \draw[blue, line width=1] (0,0) -- (1,0);
        \draw[blue, line width=1] (0,1) -- (2,1) -- (2,2) -- (0,2) -- cycle;
        \draw[blue, line width=1] (1,0.88) -- (2,0.88);
        \draw[blue, line width=1] (1,2.12) -- (2,2.12);
        \draw[blue, line width=1] (1,2.24) -- (2,2.24);
    \end {scope}

    \begin {scope}[shift={(4*1, 4*2)}]
        \draw[blue, line width=1] (0,0) -- (1,0);
        \draw[blue, line width=1] (0,1) -- (0,2);
        \draw[blue, line width=1] (-0.12,1) -- (-0.12,2);
        \draw[blue, line width=1] (1,1) -- (2,1) -- (2,2) -- (1,2) -- cycle;
        \draw[blue, line width=1] (1,2.12) -- (2,2.12);
        \draw[blue, line width=1] (0.88,1) -- (0.88,2);
        \draw[blue, line width=1] (2.12,1) -- (2.12,2);
    \end {scope}

    \begin {scope}[shift={(4*-2, 4*3)}]
        \draw[blue, line width=1] (0,0) -- (0,2) -- (2,2) -- (2,1) -- (1,1) -- (1,0);
        \draw[blue, line width=1] (1,0.88) -- (2,0.88);
        \draw[blue, line width=1] (1,2.12) -- (2,2.12);
        \draw[blue, line width=1] (1,2.24) -- (2,2.24);
    \end {scope}

    \begin {scope}[shift={(4*0, 4*3)}]
        \draw[blue, line width=1] (0,0) -- (1,0);
        \draw[blue, line width=1] (0,1) -- (2,1) -- (2,2) -- (0,2) -- cycle;
        \draw[blue, line width=1] (1,1) -- (1,2);
        \draw[blue, line width=1] (2.12,1) -- (2.12,2);
        \draw[blue, line width=1] (1,2.12) -- (2,2.12);
    \end {scope}

    \begin {scope}[shift={(4*2, 4*3)}]
        \draw[blue, line width=1] (0,0) -- (1,0);
        \draw[blue, line width=1] (0,1) -- (0,2);
        \draw[blue, line width=1] (-0.12,1) -- (-0.12,2);
        \draw[blue, line width=1] (1,1) -- (1,2) -- (2,2) -- (2,1);
        \draw[blue, line width=1] (0.88,1) -- (0.88,2);
        \draw[blue, line width=1] (0.76,1) -- (0.76,2);
        \draw[blue, line width=1] (2.12,1) -- (2.12,2);
        \draw[blue, line width=1] (2.24,1) -- (2.24,2);
    \end {scope}

    \begin {scope}[shift={(4*-1, 4*4)}]
        \draw[blue, line width=1] (0,0) -- (0,2) -- (2,2) -- (2,1) -- (1,1) -- (1,0);
        \draw[blue, line width=1] (1,1) -- (1,2);
        \draw[blue, line width=1] (1,2.12) -- (2,2.12);
        \draw[blue, line width=1] (2.12,1) -- (2.12,2);
    \end {scope}

    \begin {scope}[shift={(4*0, 4*4)}]
        \draw[blue, line width=1] (0,0) -- (1,0);
        \draw[blue, line width=1] (0,1) -- (2,1) -- (2,2) -- (0,2);
        \draw[blue, line width=1] (0,0.88) -- (1,0.88);
        \draw[blue, line width=1] (0,2.12) -- (2,2.12);
        \draw[blue, line width=1] (2.12,1) -- (2.12,2);
    \end {scope}

    \begin {scope}[shift={(4*1, 4*4)}]
        \draw[blue, line width=1] (0,0) -- (1,0);
        \draw[blue, line width=1] (0,1) -- (1,1) -- (1,2) -- (0,2) -- cycle;
        \draw[blue, line width=1] (1,2) -- (2,2) -- (2,1);
        \draw[blue, line width=1] (1.12,1) -- (1.12,2);
        \draw[blue, line width=1] (2.12,1) -- (2.12,2);
        \draw[blue, line width=1] (2.24,1) -- (2.24,2);
    \end {scope}

    \begin {scope}[shift={(4*-1, 4*5)}]
        \draw[blue, line width=1] (0,0) -- (0,1) -- (2,1) -- (2,2) -- (0,2);
        \draw[blue, line width=1] (1,0) -- (1,1);
        \draw[blue, line width=1] (2.12,1) -- (2.12,2);
        \draw[blue, line width=1] (0,2.12) -- (2,2.12);
    \end {scope}

    \begin {scope}[shift={(4*0, 4*5)}]
        \draw[blue, line width=1] (0,0) -- (0,2) -- (2,2) -- (2,1);
        \draw[blue, line width=1] (1,0) -- (1,2);
        \draw[blue, line width=1] (2.12,1) -- (2.12,2);
        \draw[blue, line width=1] (2.24,1) -- (2.24,2);
        \draw[blue, line width=1] (0.88,1) -- (0.88,2);
    \end {scope}
    
    \begin {scope}[shift={(4*1, 4*5)}]
        \draw[blue, line width=1] (0,0) -- (1,0);
        \draw[blue, line width=1] (0,1) -- (1,1) -- (1,2) -- (2,2) -- (2,1);
        \draw[blue, line width=1] (0,2) -- (1,2);
        \draw[blue, line width=1] (0,0.88) -- (1,0.88);
        \draw[blue, line width=1] (0,2.12) -- (1,2.12);
        \draw[blue, line width=1] (2.12,1) -- (2.12,2);
        \draw[blue, line width=1] (2.24,1) -- (2.24,2);
    \end {scope}

    \begin {scope}[shift={(4*0, 4*6)}]
        \draw[blue, line width=1] (0,0) -- (0,1) -- (1,1) -- (1,2) -- (2,2) -- (2,1);
        \draw[blue, line width=1] (1,0) -- (1,1);
        \draw[blue, line width=1] (0,2) -- (1,2);
        \draw[blue, line width=1] (0,2.12) -- (1,2.12);
        \draw[blue, line width=1] (2.12,1) -- (2.12,2);
        \draw[blue, line width=1] (2.24,1) -- (2.24,2);
    \end {scope}


    \begin {scope}[shift={(20,0)}]
        \foreach \x/\y in {0/0, 0/3, -1/4, 1/4, 0/5} {
            \begin {scope}[shift={(4*\x,4*\y)}]
                \draw (1,2.5) -- (1,3.5);
            \end {scope}
        }

        \foreach \x/\y in {0/1, -1/2, 1/2, -2/3, 0/3, -1/4, 0/4, -1/5} {
            \begin {scope}[shift={(4*\x,4*\y)}]
                \draw (1.5,2.5) -- (3.5,4.5);
            \end {scope}
        }

        \foreach \x/\y in {0/1, -1/2, 1/2, 0/3, 2/3, 0/4, 1/4, 1/5} {
            \begin {scope}[shift={(4*\x,4*\y)}]
                \draw (0.5,2.5) -- (-1.5,4.5);
            \end {scope}
        }

        \begin {scope}[shift={(4*0, 4*0)}]
            \draw (1,1) node {$1234$};
        \end {scope}

        \begin {scope}[shift={(4*0, 4*1)}]
            \draw (1,1) node {$2134$};
        \end {scope}

        \begin {scope}[shift={(4*-1, 4*2)}]
            \draw (1,1) node {$2314$};
        \end {scope}

        \begin {scope}[shift={(4*1, 4*2)}]
            \draw (1,1) node {$3124$};
        \end {scope}

        \begin {scope}[shift={(4*-2, 4*3)}]
            \draw (1,1) node {$2341$};
        \end {scope}

        \begin {scope}[shift={(4*0, 4*3)}]
            \draw (1,1) node {$3214$};
        \end {scope}

        \begin {scope}[shift={(4*2, 4*3)}]
            \draw (1,1) node {$4123$};
        \end {scope}

        \begin {scope}[shift={(4*-1, 4*4)}]
            \draw (1,1) node {$3241$};
        \end {scope}

        \begin {scope}[shift={(4*0, 4*4)}]
            \draw (1,1) node {$3412$};
        \end {scope}

        \begin {scope}[shift={(4*1, 4*4)}]
            \draw (1,1) node {$4213$};
        \end {scope}

        \begin {scope}[shift={(4*-1, 4*5)}]
            \draw (1,1) node {$3421$};
        \end {scope}

        \begin {scope}[shift={(4*0, 4*5)}]
            \draw (1,1) node {$4231$};
        \end {scope}

        \begin {scope}[shift={(4*1, 4*5)}]
            \draw (1,1) node {$4312$};
        \end {scope}

        \begin {scope}[shift={(4*0, 4*6)}]
            \draw (1,1) node {$4321$};
        \end {scope}
    \end {scope}
\end {tikzpicture}
\caption {The partial orders on $\Omega(\mathcal{G}^z_3)$ (left) and $\mathrm{Cat}_4$ (right).}
\label {fig:catalan_posets}
\end {figure}
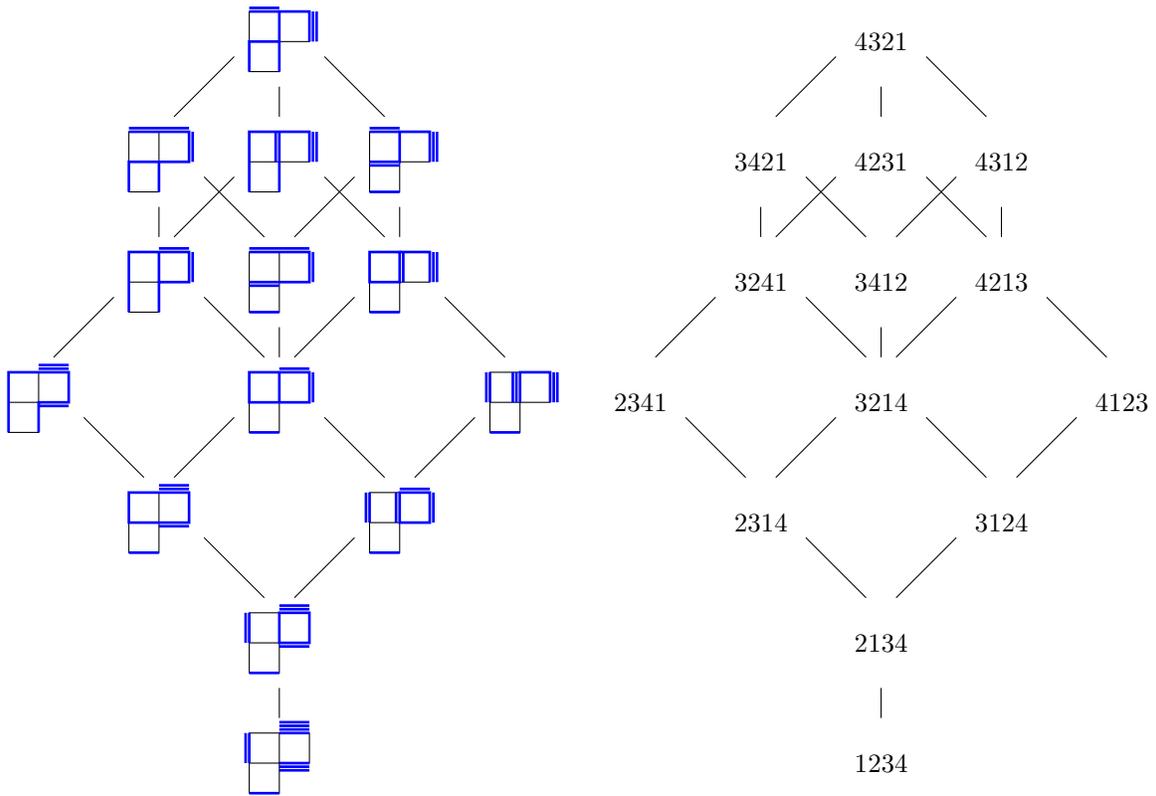

\begin {ex}
    The Hasse diagrams of the posets $\mathcal{G}^z_3$ and $\mathrm{Cat}_4$ are pictured in Figure \ref{fig:catalan_posets}.
\end {ex}

Just as with the Euler numbers, we can use the weighted version of Theorem \ref{thm:matrix_formula_zigzag} to compute the Carlitz-Riordan $q$-Catalan numbers.
In this case, if we set $a_i = b_i = 1$ and $c_i = q$ for all $i$ (using the edge labels from Figure \ref{fig:zigzag_edge_weights}), then every face weight will be equal to $q$, and the minimal mixed dimer cover
will have weight $1$ (since it does not use any $c$-edges). 

\begin {ex}
    To compute $C_4(q)$, corresponding to Figure \ref{fig:catalan_posets}, we use the matrix product
    \[ 
        X = \begin{pmatrix} 1&1 \\ 1&0 \end{pmatrix}
            \begin{pmatrix} 1&1&1 \\ q&q&0 \end{pmatrix}
            \begin{pmatrix} q^2&q^2&0&0 \\ q&q&q&0 \\ 1&1&1&1 \end{pmatrix}
            \begin{pmatrix} q^3&q^3&0&0&0 \\ q^2&q^2&q^2&0&0 \\ q&q&q&q&0 \\ 1&1&1&1&1 \end{pmatrix}
    \]
    The $(1,1)$-entry of this matrix is the $q$-Catalan number $C_4(q) = 1+q+2q^2+3q^3+3q^4+3q^5+q^6$.
\end {ex}

\bigskip

\section {Other Combinatorial Models}

\subsection {Networks} \label{sec:networks}

Let $N$ be an acyclic, edge-weighted planar directed graph (which we will refer to as a \emph{network}).
Suppose there are $k$ source vertices $u_1,\dots,u_k$ and $n$ sink vertices $v_1,\dots,v_n$.
We may form a $k \times n$ \emph{path-weight matrix} $M$ whose $(i,j)$ entry is the weighted sum of paths from $u_i$ to $v_j$.
It is easy to see that multiplying such matrices corresponds to glueing/concatenation of networks. 

Recall that $W_a$ is the $(a+1) \times (a+1)$ anti-diagonal matrix, and that $L_{a,b} := W_a R_{a,b}$. 
Similarly, we define $U_{a,b} := R_{a,b} W_b$. Note that since $W_a^2 = \mathrm{Id}_a$, then $R_{a,b}R_{b,c} = U_{a,b}L_{b,c}$.
Therefore if we can identify networks which represent the $L$ and $U$ matrices, then by concatenating them in the proper sequences, we
can draw networks which represent the matrices from Theorem \ref{thm:matrix_formula} and Theorem \ref{thm:matrix_formula_zigzag}.

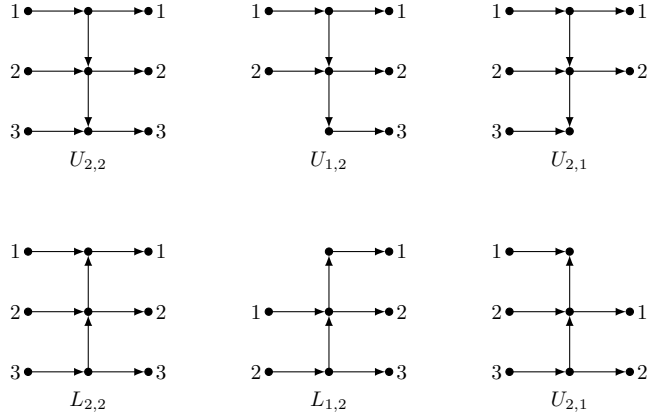
\begin {figure}
\centering
\begin {tikzpicture}[scale=0.8, every node/.style={scale=0.8}]
    \foreach \x in {0,1,2} {
        \foreach \y in {0,1,2} {
            \draw[fill=black] (\x,\y) circle (0.06);
        }
    }
    \foreach \x in {0,1} {
        \foreach \y in {0,1,2} {
            \draw[-latex] (\x,\y) -- (\x+0.94,\y);
        }
    }
    \foreach \y in {1,2} {
        \draw[-latex] (1,\y) -- (1,\y-0.94);
    }
    \draw (1,-0.5) node {$U_{2,2}$};
    \draw (0,2) node[left] {$1$};
    \draw (0,1) node[left] {$2$};
    \draw (0,0) node[left] {$3$};
    \draw (2,2) node[right] {$1$};
    \draw (2,1) node[right] {$2$};
    \draw (2,0) node[right] {$3$};

    \begin {scope}[shift={(4,0)}]
        \draw[fill=black] (0,1) circle (0.06);
        \draw[fill=black] (0,2) circle (0.06);
        \draw[fill=black] (1,0) circle (0.06);
        \draw[fill=black] (1,1) circle (0.06);
        \draw[fill=black] (1,2) circle (0.06);
        \draw[fill=black] (2,0) circle (0.06);
        \draw[fill=black] (2,1) circle (0.06);
        \draw[fill=black] (2,2) circle (0.06);
        
        \draw[-latex] (0,1) -- (0.94,1);
        \draw[-latex] (0,2) -- (0.94,2);
        \draw[-latex] (1,0) -- (1.94,0);
        \draw[-latex] (1,1) -- (1.94,1);
        \draw[-latex] (1,2) -- (1.94,2);
        
        \foreach \y in {1,2} {
            \draw[-latex] (1,\y) -- (1,\y-0.94);
        }
        \draw (1,-0.5) node {$U_{1,2}$};
        \draw (0,2) node[left] {$1$};
        \draw (0,1) node[left] {$2$};
        \draw (2,2) node[right] {$1$};
        \draw (2,1) node[right] {$2$};
        \draw (2,0) node[right] {$3$};
    \end {scope}
    
    \begin {scope}[shift={(8,0)}]
        \draw[fill=black] (0,0) circle (0.06);
        \draw[fill=black] (0,1) circle (0.06);
        \draw[fill=black] (0,2) circle (0.06);
        \draw[fill=black] (1,0) circle (0.06);
        \draw[fill=black] (1,1) circle (0.06);
        \draw[fill=black] (1,2) circle (0.06);
        \draw[fill=black] (2,1) circle (0.06);
        \draw[fill=black] (2,2) circle (0.06);
        
        \draw[-latex] (0,1) -- (0.94,1);
        \draw[-latex] (0,2) -- (0.94,2);
        \draw[-latex] (0,0) -- (0.94,0);
        \draw[-latex] (1,1) -- (1.94,1);
        \draw[-latex] (1,2) -- (1.94,2);
        
        \foreach \y in {1,2} {
            \draw[-latex] (1,\y) -- (1,\y-0.94);
        }
        \draw (1,-0.5) node {$U_{2,1}$};
        \draw (0,2) node[left] {$1$};
        \draw (0,1) node[left] {$2$};
        \draw (0,0) node[left] {$3$};
        \draw (2,2) node[right] {$1$};
        \draw (2,1) node[right] {$2$};
    \end {scope}

    \begin {scope}[shift={(0,-4)}]
        \foreach \x in {0,1,2} {
            \foreach \y in {0,1,2} {
                \draw[fill=black] (\x,\y) circle (0.06);
            }
        }
        \foreach \x in {0,1} {
            \foreach \y in {0,1,2} {
                \draw[-latex] (\x,\y) -- (\x+0.94,\y);
            }
        }
        \foreach \y in {0,1} {
            \draw[-latex] (1,\y) -- (1,\y+0.94);
        }
        \draw (1,-0.5) node {$L_{2,2}$};
        \draw (0,2) node[left] {$1$};
        \draw (0,1) node[left] {$2$};
        \draw (0,0) node[left] {$3$};
        \draw (2,2) node[right] {$1$};
        \draw (2,1) node[right] {$2$};
        \draw (2,0) node[right] {$3$};
    \end {scope}
    
    \begin {scope}[shift={(4,-4)}]
        \draw[fill=black] (0,0) circle (0.06);
        \draw[fill=black] (0,1) circle (0.06);
        \draw[fill=black] (1,0) circle (0.06);
        \draw[fill=black] (1,1) circle (0.06);
        \draw[fill=black] (1,2) circle (0.06);
        \draw[fill=black] (2,0) circle (0.06);
        \draw[fill=black] (2,1) circle (0.06);
        \draw[fill=black] (2,2) circle (0.06);
        
        \draw[-latex] (0,0) -- (0.94,0);
        \draw[-latex] (0,1) -- (0.94,1);
        \draw[-latex] (1,0) -- (1.94,0);
        \draw[-latex] (1,1) -- (1.94,1);
        \draw[-latex] (1,2) -- (1.94,2);
        
        \foreach \y in {0,1} {
            \draw[-latex] (1,\y) -- (1,\y+0.94);
        }
        \draw (1,-0.5) node {$L_{1,2}$};
        \draw (0,0) node[left] {$2$};
        \draw (0,1) node[left] {$1$};
        \draw (2,2) node[right] {$1$};
        \draw (2,1) node[right] {$2$};
        \draw (2,0) node[right] {$3$};
    \end {scope}
    
    \begin {scope}[shift={(8,-4)}]
        \draw[fill=black] (0,0) circle (0.06);
        \draw[fill=black] (0,1) circle (0.06);
        \draw[fill=black] (0,2) circle (0.06);
        \draw[fill=black] (1,0) circle (0.06);
        \draw[fill=black] (1,1) circle (0.06);
        \draw[fill=black] (1,2) circle (0.06);
        \draw[fill=black] (2,1) circle (0.06);
        \draw[fill=black] (2,0) circle (0.06);
        
        \draw[-latex] (0,1) -- (0.94,1);
        \draw[-latex] (0,2) -- (0.94,2);
        \draw[-latex] (0,0) -- (0.94,0);
        \draw[-latex] (1,1) -- (1.94,1);
        \draw[-latex] (1,0) -- (1.94,0);
        
        \foreach \y in {0,1} {
            \draw[-latex] (1,\y) -- (1,\y+0.94);
        }
        \draw (1,-0.5) node {$U_{2,1}$};
        \draw (0,2) node[left] {$1$};
        \draw (0,1) node[left] {$2$};
        \draw (0,0) node[left] {$3$};
        \draw (2,0) node[right] {$2$};
        \draw (2,1) node[right] {$1$};
    \end {scope}
\end {tikzpicture}
\caption {Examples of the networks for $U_{a,b}$ and $L_{a,b}$.}
\label {fig:networks}
\end {figure}

It is quite simple to find such networks. Let us first describe the case for square matrices. The matrix $U_{n,n}$ is represented by the network
with $n$ horizontal strands directed left-to-right, and downward-pointing vertical edges connecting adjacent strands. The network for $L_{n,n}$ is
similar, but with the vertical edges oriented upwards. For non-square matrices, the associated networks can be obtained from these by deleting 
some of the source or sink vertices. The different possibilities are pictured in Figure \ref{fig:networks}.

\begin {rmk}
    One of the standard combinatorial models of the Catalan numbers is counting north-east lattice paths from $(0,0)$ to $(n,n)$ which do not go above
    the line $y=x$. This corresponds to paths through the network $L_{11}L_{12}L_{23} \cdots L_{n-1,n}$ starting at vertex 2 on the left, and ending
    at vertex 1 on the right. Note that $U_{11} = XL_{11}$, where $X = \left( \begin{smallmatrix} 0&1\\-1&1\end{smallmatrix} \right)$, which means that
    the first row of $U_{11}L_{12}L_{23} \cdots L_{n-1,n}$ is the same as the second row of $L_{11}L_{12} \cdots L_{n-1,n}$, and the former is precisely
    the matrix given in Theorem \ref{thm:matrix_formula_zigzag}.
\end {rmk}

\begin {rmk}
    The \emph{Fuss-Catalan numbers} $\frac{1}{kn+1}\binom{(k+1)n}{n}$ (which include the ordinary Catalan numbers in the case $k=1$) 
    are known to count north-east lattice paths from $(0,0)$ to $(n,kn)$ which do not go above the line $y=kx$ 
    (for example, see \cite[Prop 6.2.1($v$)]{stanley_EC2}). In the notations introduced in this section, this is the number
    of paths in the network $L_{k,k} L_{k,2k} L_{2k,3k} \cdots L_{(n-1)k,nk}$ from vertex $k$ on the left to vertex $1$ on the right.
    By the same argument as in the previous remark, this is the number of mixed dimer covers on the zigzag snake $\mathcal{G}^z_{n-1}$
    with vertex labels $m = (k, 2k, 3k, \dots, nk)$.
\end {rmk}

\begin {rmk}
    By Theorem \ref{thm:matrix_formula}, the network for the Euler numbers is given by $U_{11}L_{12}U_{23}L_{34} \cdots$, with alternating $U$'s and $L$'s.
    This interpretation of the Euler numbers in terms of paths in this network appeared in \cite{boustrophedon}. An example of such a network
    appears in Figure \ref{fig:plabic_transform}.
\end {rmk}

\subsection {Single Dimer Covers}

We will give families of planar graphs whose number of perfect matchings are the Euler and Catalan numbers. To do so,
we will use the networks of the previous section, together with results from \cite{psw} which give bijections between paths 
in a network to perfect matchings. 

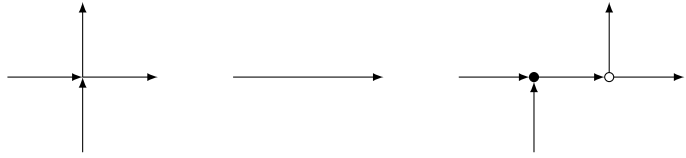
\begin {figure}
\centering
\begin {tikzpicture}
    \draw[-latex] (0,-1) -- (0,0);
    \draw[-latex] (-1,0) -- (0,0);
    \draw[-latex] (0,0) -- (1,0);
    \draw[-latex] (0,0) -- (0,1);

    \draw[-latex] (2,0) -- (4,0);

    \begin {scope}[shift={(6,0)}]
        \draw[-latex] (-1,0) -- (-0.06,0);
        \draw[-latex] (0,-1) -- (0,-0.06);
        \draw[-latex] (0,0) -- (0.94,0);
        \draw[-latex] (1,0) -- (2,0);
        \draw[-latex] (1,0) -- (1,1);
        \draw[fill=black] (0,0) circle (0.06);
        \draw[fill=white] (1,0) circle (0.06);
    \end {scope}
\end {tikzpicture}
\caption {Splitting a 4-valent vertex into two 3-valent vertices.}
\label {fig:vertex_split}
\end {figure}

We first start with a network obtained by concatenation of the $U_{a,b}$ and $L_{a,b}$ from Figure \ref{fig:networks}.
In order to conform to the setup in \cite{psw}, we must have a \emph{perfectly oriented network}. This means the graph must be
bipartite, and each white vertex must have a unique incoming arrow, and each black vertex must have a unique outgoing arrow. 
We can obtain such a network by splitting each 4-valent vertex into two 3-valent vertices (while preserving the orientations).
See Figure \ref{fig:vertex_split} for an illustration. Since the original network had all square faces, the resulting graph
can be drawn (after possibly inserting some more 2-valent vertices) so that all faces are hexagons.

Since we only care about paths from vertex 1 on the left to vertex 1 on the right, 
we may delete the remaining boundary vertices. Corollary 4.6 from \cite{psw} then says that there is a bijection between paths in this
network (from the source to the sink) and the set of perfect matchings. Explicitly, given a path, the corresponding matching consists of those edges
off the path directed from black-to-white, together with the edges along the path directed from white-to-black.

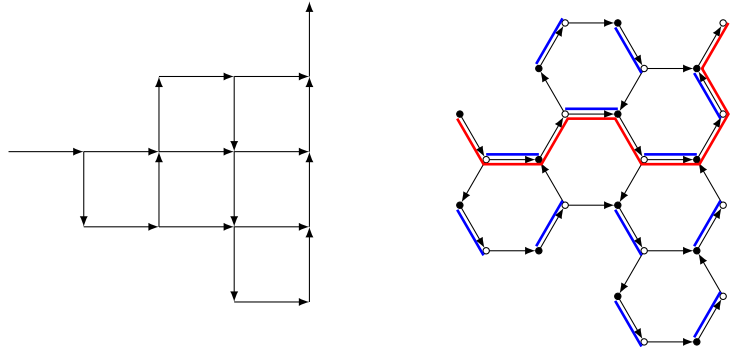
\begin {figure}
\centering
\begin {tikzpicture}
    \foreach \x/\y in {1/0, 2/0, 3/0, 3/-1, 0/1, 1/1, 2/1, 3/1, 2/2, 3/2} {
        \draw[-latex] (\x,\y) -- (\x+1,\y);
    }
    \foreach \x/\y in {1/1, 3/0, 3/1, 3/2} {
        \draw[-latex] (\x,\y) -- (\x,\y-1);
    }
    \foreach \x/\y in {2/0, 2/1, 4/-1, 4/0, 4/1, 4/2} {
        \draw[-latex] (\x,\y) -- (\x,\y+1);
    }

    \begin {scope}[shift={(6,1.5)}, scale=0.7]
        \foreach \x/\y in {-1/2, 0/0, 0/1, 0/2, 1/-1, 1/0, 1/1, 2/0, 3/-1} {
            \draw[-latex] ($(0.5, {-sqrt(3)*0.5}) + ({1.5*(\y+\x)}, {sqrt(3)*0.5*(\y-\x)})$) -- ($(0.5, {-sqrt(3)*0.5}) + ({1.5*(\y+\x)}, {sqrt(3)*0.5*(\y-\x)}) + (0.94,0)$);
        }
        \foreach \x/\y in {0/0, 1/-1, 0/2, 1/1, 2/0, 3/-1} {
            \draw[-latex] ($({1.5*(\y+\x)}, {sqrt(3)*0.5*(\y-\x)})$) -- ($({1.5*(\y+\x)}, {sqrt(3)*0.5*(\y-\x)}) + 0.94*(-60:1)$);
        }
        \foreach \x/\y in {0/0, 0/2, 1/1, 2/0} {
            \draw[-latex] ($(0.5, {-sqrt(3)*0.5}) + ({1.5*(\y+\x)}, {sqrt(3)*0.5*(\y-\x)})$) -- ($(0.5, {-sqrt(3)*0.5}) + ({1.5*(\y+\x)}, {sqrt(3)*0.5*(\y-\x)}) + 0.94*(-120:1)$);
        }
        \foreach \x/\y in {2/-1, 1/0, 0/1, 4/-1, 3/0, 2/1, 1/2} {
            \draw[-latex] ($({1.5*(\y+\x)}, {sqrt(3)*0.5*(\y-\x)})$) -- ($({1.5*(\y+\x)}, {sqrt(3)*0.5*(\y-\x)}) + 0.94*(60:1)$);
        }
        \foreach \x/\y in {1/0, 0/1, 3/0, 2/1, 1/2} {
            \draw[-latex] ($(0.5, {-sqrt(3)*0.5}) + ({1.5*(\y+\x)}, {sqrt(3)*0.5*(\y-\x)})$) -- ($(0.5, {-sqrt(3)*0.5}) + ({1.5*(\y+\x)}, {sqrt(3)*0.5*(\y-\x)}) + 0.94*(120:1)$);
        }
        
        \foreach \x/\y in {0/0, 0/1, 0/2, 1/-1, 1/0, 1/1, 1/2, 2/-1, 2/0, 2/1, 3/-1, 3/0, 4/-1} {
            \draw[fill=black] ($({1.5*(\y+\x)}, {sqrt(3)*0.5*(\y-\x)})$) circle (0.06);
        }

        \foreach \x/\y in {-1/2, 0/0, 0/1, 0/2, 0/3, 1/-1, 1/0, 1/1, 1/2, 2/0, 2/1, 3/-1, 3/0} {
            \draw[fill=white] ($(0.5, {-sqrt(3)*0.5}) + ({1.5*(\y+\x)}, {sqrt(3)*0.5*(\y-\x)})$) circle (0.06);
        }

        \draw[red, line width=1] ($(0,0) + 0.1*(-120:1)$) --++ (-60:1) --++ (1.1,0) --++ (60:1) --++ (0.9,0) --++ (-60:1) --++ (1.1,0) --++ (60:1.1) --++ (120:1) --++ (60:1);

        \foreach \x/\y in {1/-1, 0/2, 2/0, 3/-1, 1/2} {
            \draw[blue, line width=1] ($({1.5*(\y+\x)}, {sqrt(3)*0.5*(\y-\x)}) + 0.1*(-120:1)$) --++ (-60:1);
        }
        \foreach \x/\y in {0/1, 2/-1, 3/0, 4/-1} {
            \draw[blue, line width=1] ($({1.5*(\y+\x)}, {sqrt(3)*0.5*(\y-\x)}) + 0.1*(120:1)$) --++ (60:1);
        }
        \foreach \x/\y in {1/0, 1/1, 2/1} {
            \draw[blue, line width=1] ($({1.5*(\y+\x)}, {sqrt(3)*0.5*(\y-\x)}) + (0,0.1)$) --++ (180:1);
        }
    \end {scope}
\end {tikzpicture}
\caption {Transforming a network into a perfectly oriented bipartite network. The right figure also illustrates the bijection between paths and matchings.}
\label {fig:plabic_transform}
\end {figure}

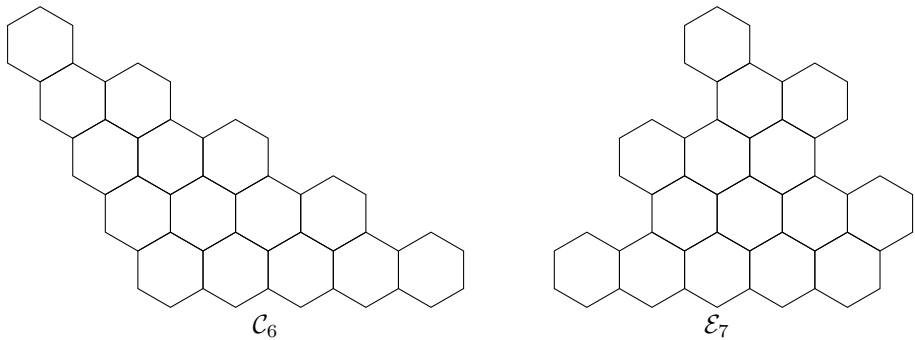
\begin {figure}
\centering
\begin {tikzpicture}[scale=0.5]
    \foreach \x/\y in {5/0, 6/0, 7/0, 3/0, 5/1, 3/1, 2/1, 2/2, 1/2, 0/3, 4/0, 3/2, -1/4, 4/1, 1/3} {
        \begin {scope}[shift={({sqrt(3)*\x + 0.5*sqrt(3)*\y},{1.5*\y})}]
            \draw (30:1) -- (90:1) -- (150:1) -- (210:1) -- (270:1) -- (330:1) -- cycle;
        \end {scope}
    }
    \draw ({6+sqrt(3)},-1.5) node {$\mathcal{C}_6$};

    \begin {scope}[shift={(18,0)}]
    \foreach \x/\y in {0/0, 1/0, 2/0, 3/0, 0/1, 1/1, 2/1, 0/2, 1/2, 0/3, -1/0, -1/2, -1/4, 3/1, 1/3} {
        \begin {scope}[shift={({sqrt(3)*\x + 0.5*sqrt(3)*\y},{1.5*\y})}]
            \draw (30:1) -- (90:1) -- (150:1) -- (210:1) -- (270:1) -- (330:1) -- cycle;
        \end {scope}
    }
    \draw ({sqrt(3)},-1.5) node {$\mathcal{E}_7$};
    \end {scope}
\end {tikzpicture}
\caption {Examples of the graphs $\mathcal{C}_n$ and $\mathcal{E}_n$}
\label {fig:dimer_graphs}
\end {figure}

Doing this procedure for the networks obtained from the matrix product in Theorem \ref{thm:matrix_formula} in the case of the standard labeling
gives a family of planar bipartite graphs $\mathcal{E}_n$ whose number of perfect matchings are the Euler numbers $E_n$. The graph $\mathcal{E}_5$
is pictured in Figure \ref{fig:plabic_transform} and $\mathcal{E}_7$ is pictured in Figure \ref{fig:dimer_graphs}.

Doing this procedure for the networks obtained from the matrix product in Theorem \ref{thm:matrix_formula_zigzag} in the case of the standard labeling
gives a family of planar bipartite graphs $\mathcal{C}_n$ whose number of perfect matchings are the Catalan numbers $C_n$. The graph $\mathcal{C}_6$ 
is pictured in Figure \ref{fig:dimer_graphs}. These graphs have appeared in the literature before (for example, in \cite{doslic}, 
in \cite{matchable} where they were called \emph{prolate triangles}, and in Exercise 206 from \cite{stanley_catalan})

\begin {rmk}
    Zhang, Yang, and Yao called a distributive lattice \emph{matchable} if it is isomorphic to the lattice of perfect matchings of some planar graph \cite{matchable},
    and asked the interesting question of which lattices are matchable. By the preceding discussion, we see that the lattices of mixed dimer covers on any straight or zigzag snake graph (with any vertex labeling) is a matchable lattice.
\end {rmk}

\subsection {Rhombus Tilings and Plane Partitions}

There is a standard well-known correspondence between perfect matchings of a hexagonal grid and rhombus tilings of an associated region.
For each edge in the matching, draw a rhombus surrounding it, whose vertices are the centers of the neighboring hexagons. These can also be viewed
as 3-dimensional pictures depicting plane partitions (stacks of cubes).

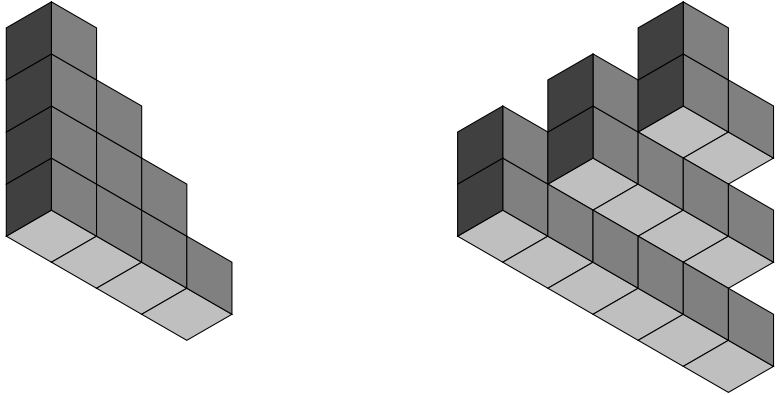
\begin {figure}
\centering
\begin {tikzpicture} 
    \begin {scope}[scale=0.4,rotate=-30]
        \foreach \x/\y in {0/0,1/0,2/0,3/0} {
            \def\xx{sqrt(3)*(\x+0.5*\y)}
            \def\yy{1.5*\y}
            \draw[fill=lightgray] ({\xx},\yy) -- ({\xx + sqrt(3)},\yy) -- ({\xx + 1.5*sqrt(3)}, \yy+1.5) -- ({\xx + 0.5*sqrt(3)}, \yy+1.5) -- cycle;
        }
    
        \foreach \x/\y in {1/1, 2/1, 3/1, 4/1, 1/2, 2/2, 3/2, 1/3, 2/3, 1/4} {
            \def\xx{sqrt(3)*(\x-0.5*\y)}
            \def\yy{1.5*\y}
            \draw[fill=gray] ({\xx}, {\yy}) -- ({\xx - 0.5*sqrt(3)}, {\yy+1.5}) -- ({\xx + 0.5*sqrt(3)}, {\yy+1.5}) -- ({\xx + sqrt(3)}, {\yy}) -- cycle;
        }
    
        \foreach \x/\y in {0/0, 0/1, 0/2, 0/3} {
            \def\xx{0.5*sqrt(3)*(\x-\y)}
            \def\yy{1.5*(\x+\y)}
            \draw[fill=darkgray] ({\xx}, {\yy}) -- ({\xx - 0.5*sqrt(3)}, {\yy+1.5}) -- ({\xx}, {\yy+3}) -- ({\xx + 0.5*sqrt(3)}, {\yy+1.5}) -- cycle;
        }
    \end {scope}
    
    \begin {scope}[shift={(6,0)}, scale=0.4,rotate=-30]
        \foreach \x/\y in {0/0,1/0,2/0,3/0,4/0,5/0,0/2,1/2,2/2,3/2,0/4,1/4} {
            \def\xx{sqrt(3)*(\x+0.5*\y)}
            \def\yy{1.5*\y}
            \draw[fill=lightgray] ({\xx},\yy) -- ({\xx + sqrt(3)},\yy) -- ({\xx + 1.5*sqrt(3)}, \yy+1.5) -- ({\xx + 0.5*sqrt(3)}, \yy+1.5) -- cycle;
        }
    
        \foreach \x/\y in {1/1,2/1,3/1,4/1,5/1,6/1,1/2,3/3,4/3,5/3,6/3,3/4,5/5,6/5,5/6} {
            \def\xx{sqrt(3)*(\x-0.5*\y)}
            \def\yy{1.5*\y}
            \draw[fill=gray] ({\xx}, {\yy}) -- ({\xx - 0.5*sqrt(3)}, {\yy+1.5}) -- ({\xx + 0.5*sqrt(3)}, {\yy+1.5}) -- ({\xx + sqrt(3)}, {\yy}) -- cycle;
        }
    
        \foreach \x/\y in {0/0,0/1,2/0,2/1,4/0,4/1} {
            \def\xx{0.5*sqrt(3)*(\x-\y)}
            \def\yy{1.5*(\x+\y)}
            \draw[fill=darkgray] ({\xx}, {\yy}) -- ({\xx - 0.5*sqrt(3)}, {\yy+1.5}) -- ({\xx}, {\yy+3}) -- ({\xx + 0.5*sqrt(3)}, {\yy+1.5}) -- cycle;
        }
    \end {scope}
\end {tikzpicture}
\caption {Rhombus tilings corresponding to $\mathcal{C}_n$ (left) and $\mathcal{E}_n$ (right).}
\label {fig:rhombus_tilings}
\end {figure}

Since the graphs obtained in the previous section are regions in the hexagonal grid, we obtain for any labeled straight or zigzag snake graph a 
corresponding region of the plane whose rhomubs tilings are equinumerous with the mixed dimer covers of the original snake.
See Figure \ref{fig:rhombus_tilings} for examples corresponding to the minimal matchings of $\mathcal{C}_n$ and $\mathcal{E}_n$.

\subsection {Order Ideals}

Birkhoff's representation theorem says that for any finite distributive lattice $\Omega$, there exists a poset $P_{\Omega}$ whose set of order ideals form a lattice (under containment) isomorphic to $\Omega$ \cite[Thm 3.4.1]{stanley_EC1}. By Theorem \ref{thm:lattice} (to be proved below in Section \ref{sec:dimers_paths}), the mixed dimer covers on a fixed snake graph $\mathcal{G}$ in fact form a distributive lattice $\Omega$. Thus, there is an underlying poset $P_{\Omega}$ whose lattice of order ideals is isomorphic to the lattice of mixed dimer covers on $\mathcal{G}$. We now describe these underlying posets $P_{\Omega}$ for the two special cases of the mixed dimer lattices on $\mathcal{G}_{n}^{s}$ and $\mathcal{G}_{n}^{z}$ (each equipped with their standard labeling). We will use the notations $P_{n}^s$ and $P_n^z$ for the posets underlying the lattices of mixed dimer covers on $\mathcal{G}_{n}^s$ and $\mathcal{G}_n^z$.

\begin{figure}
\centering
\tikzset{every picture/.style={line width=0.75pt}} 
\begin{tikzpicture}[x=0.75pt,y=0.75pt,yscale=-.7,xscale=.7, every node/.style={scale=0.8}]

\draw    (377.1,108.16) -- (322.47,53.53) ;
\draw [shift={(320.81,51.87)}, rotate = 225] [color={rgb, 255:red, 0; green, 0; blue, 0 }  ][line width=0.75]      (0, 0) circle [x radius= 3.35, y radius= 3.35]   ;
\draw [shift={(354.38,85.44)}, rotate = 225] [fill={rgb, 255:red, 0; green, 0; blue, 0 }  ][line width=0.08]  [draw opacity=0] (8.93,-4.29) -- (0,0) -- (8.93,4.29) -- cycle    ;
\draw [shift={(378.76,109.82)}, rotate = 225] [color={rgb, 255:red, 0; green, 0; blue, 0 }  ][line width=0.75]      (0, 0) circle [x radius= 3.35, y radius= 3.35]   ;
\draw    (376.41,109.82) -- (323.29,109.82) ;
\draw [shift={(320.94,109.82)}, rotate = 180] [color={rgb, 255:red, 0; green, 0; blue, 0 }  ][line width=0.75]      (0, 0) circle [x radius= 3.35, y radius= 3.35]   ;
\draw [shift={(356.35,109.82)}, rotate = 180] [fill={rgb, 255:red, 0; green, 0; blue, 0 }  ][line width=0.08]  [draw opacity=0] (8.93,-4.29) -- (0,0) -- (8.93,4.29) -- cycle    ;
\draw [shift={(378.76,109.82)}, rotate = 180] [color={rgb, 255:red, 0; green, 0; blue, 0 }  ][line width=0.75]      (0, 0) circle [x radius= 3.35, y radius= 3.35]   ;
\draw    (377.1,111.48) -- (322.47,166.11) ;
\draw [shift={(320.81,167.77)}, rotate = 135] [color={rgb, 255:red, 0; green, 0; blue, 0 }  ][line width=0.75]      (0, 0) circle [x radius= 3.35, y radius= 3.35]   ;
\draw [shift={(346.25,142.33)}, rotate = 315] [fill={rgb, 255:red, 0; green, 0; blue, 0 }  ][line width=0.08]  [draw opacity=0] (8.93,-4.29) -- (0,0) -- (8.93,4.29) -- cycle    ;
\draw [shift={(378.76,109.82)}, rotate = 135] [color={rgb, 255:red, 0; green, 0; blue, 0 }  ][line width=0.75]      (0, 0) circle [x radius= 3.35, y radius= 3.35]   ;
\draw    (318.46,167.77) -- (265.34,167.77) ;
\draw [shift={(262.99,167.77)}, rotate = 180] [color={rgb, 255:red, 0; green, 0; blue, 0 }  ][line width=0.75]      (0, 0) circle [x radius= 3.35, y radius= 3.35]   ;
\draw [shift={(286.9,167.77)}, rotate = 360] [fill={rgb, 255:red, 0; green, 0; blue, 0 }  ][line width=0.08]  [draw opacity=0] (8.93,-4.29) -- (0,0) -- (8.93,4.29) -- cycle    ;
\draw [shift={(320.81,167.77)}, rotate = 180] [color={rgb, 255:red, 0; green, 0; blue, 0 }  ][line width=0.75]      (0, 0) circle [x radius= 3.35, y radius= 3.35]   ;
\draw    (319.28,111.48) -- (264.65,166.11) ;
\draw [shift={(262.99,167.77)}, rotate = 135] [color={rgb, 255:red, 0; green, 0; blue, 0 }  ][line width=0.75]      (0, 0) circle [x radius= 3.35, y radius= 3.35]   ;
\draw [shift={(288.43,142.33)}, rotate = 315] [fill={rgb, 255:red, 0; green, 0; blue, 0 }  ][line width=0.08]  [draw opacity=0] (8.93,-4.29) -- (0,0) -- (8.93,4.29) -- cycle    ;
\draw [shift={(320.94,109.82)}, rotate = 135] [color={rgb, 255:red, 0; green, 0; blue, 0 }  ][line width=0.75]      (0, 0) circle [x radius= 3.35, y radius= 3.35]   ;
\draw    (319.28,224.06) -- (264.65,169.43) ;
\draw [shift={(262.99,167.77)}, rotate = 225] [color={rgb, 255:red, 0; green, 0; blue, 0 }  ][line width=0.75]      (0, 0) circle [x radius= 3.35, y radius= 3.35]   ;
\draw [shift={(296.56,201.34)}, rotate = 225] [fill={rgb, 255:red, 0; green, 0; blue, 0 }  ][line width=0.08]  [draw opacity=0] (8.93,-4.29) -- (0,0) -- (8.93,4.29) -- cycle    ;
\draw [shift={(320.94,225.72)}, rotate = 225] [color={rgb, 255:red, 0; green, 0; blue, 0 }  ][line width=0.75]      (0, 0) circle [x radius= 3.35, y radius= 3.35]   ;
\draw    (377.1,224.06) -- (322.47,169.43) ;
\draw [shift={(320.81,167.77)}, rotate = 225] [color={rgb, 255:red, 0; green, 0; blue, 0 }  ][line width=0.75]      (0, 0) circle [x radius= 3.35, y radius= 3.35]   ;
\draw [shift={(354.38,201.34)}, rotate = 225] [fill={rgb, 255:red, 0; green, 0; blue, 0 }  ][line width=0.08]  [draw opacity=0] (8.93,-4.29) -- (0,0) -- (8.93,4.29) -- cycle    ;
\draw [shift={(378.76,225.72)}, rotate = 225] [color={rgb, 255:red, 0; green, 0; blue, 0 }  ][line width=0.75]      (0, 0) circle [x radius= 3.35, y radius= 3.35]   ;
\draw    (376.28,167.77) -- (323.16,167.77) ;
\draw [shift={(320.81,167.77)}, rotate = 180] [color={rgb, 255:red, 0; green, 0; blue, 0 }  ][line width=0.75]      (0, 0) circle [x radius= 3.35, y radius= 3.35]   ;
\draw [shift={(344.72,167.77)}, rotate = 360] [fill={rgb, 255:red, 0; green, 0; blue, 0 }  ][line width=0.08]  [draw opacity=0] (8.93,-4.29) -- (0,0) -- (8.93,4.29) -- cycle    ;
\draw [shift={(378.63,167.77)}, rotate = 180] [color={rgb, 255:red, 0; green, 0; blue, 0 }  ][line width=0.75]      (0, 0) circle [x radius= 3.35, y radius= 3.35]   ;
\draw    (434.92,224.06) -- (380.29,169.43) ;
\draw [shift={(378.63,167.77)}, rotate = 225] [color={rgb, 255:red, 0; green, 0; blue, 0 }  ][line width=0.75]      (0, 0) circle [x radius= 3.35, y radius= 3.35]   ;
\draw [shift={(412.2,201.34)}, rotate = 225] [fill={rgb, 255:red, 0; green, 0; blue, 0 }  ][line width=0.08]  [draw opacity=0] (8.93,-4.29) -- (0,0) -- (8.93,4.29) -- cycle    ;
\draw [shift={(436.58,225.72)}, rotate = 225] [color={rgb, 255:red, 0; green, 0; blue, 0 }  ][line width=0.75]      (0, 0) circle [x radius= 3.35, y radius= 3.35]   ;
\draw    (376.41,225.72) -- (323.29,225.72) ;
\draw [shift={(320.94,225.72)}, rotate = 180] [color={rgb, 255:red, 0; green, 0; blue, 0 }  ][line width=0.75]      (0, 0) circle [x radius= 3.35, y radius= 3.35]   ;
\draw [shift={(356.35,225.72)}, rotate = 180] [fill={rgb, 255:red, 0; green, 0; blue, 0 }  ][line width=0.08]  [draw opacity=0] (8.93,-4.29) -- (0,0) -- (8.93,4.29) -- cycle    ;
\draw [shift={(378.76,225.72)}, rotate = 180] [color={rgb, 255:red, 0; green, 0; blue, 0 }  ][line width=0.75]      (0, 0) circle [x radius= 3.35, y radius= 3.35]   ;
\draw    (434.23,225.72) -- (381.11,225.72) ;
\draw [shift={(378.76,225.72)}, rotate = 180] [color={rgb, 255:red, 0; green, 0; blue, 0 }  ][line width=0.75]      (0, 0) circle [x radius= 3.35, y radius= 3.35]   ;
\draw [shift={(414.17,225.72)}, rotate = 180] [fill={rgb, 255:red, 0; green, 0; blue, 0 }  ][line width=0.08]  [draw opacity=0] (8.93,-4.29) -- (0,0) -- (8.93,4.29) -- cycle    ;
\draw [shift={(436.58,225.72)}, rotate = 180] [color={rgb, 255:red, 0; green, 0; blue, 0 }  ][line width=0.75]      (0, 0) circle [x radius= 3.35, y radius= 3.35]   ;
\draw    (318.59,225.72) -- (265.47,225.72) ;
\draw [shift={(263.12,225.72)}, rotate = 180] [color={rgb, 255:red, 0; green, 0; blue, 0 }  ][line width=0.75]      (0, 0) circle [x radius= 3.35, y radius= 3.35]   ;
\draw [shift={(298.53,225.72)}, rotate = 180] [fill={rgb, 255:red, 0; green, 0; blue, 0 }  ][line width=0.08]  [draw opacity=0] (8.93,-4.29) -- (0,0) -- (8.93,4.29) -- cycle    ;
\draw [shift={(320.94,225.72)}, rotate = 180] [color={rgb, 255:red, 0; green, 0; blue, 0 }  ][line width=0.75]      (0, 0) circle [x radius= 3.35, y radius= 3.35]   ;
\draw    (261.46,227.38) -- (206.83,282.01) ;
\draw [shift={(205.17,283.67)}, rotate = 135] [color={rgb, 255:red, 0; green, 0; blue, 0 }  ][line width=0.75]      (0, 0) circle [x radius= 3.35, y radius= 3.35]   ;
\draw [shift={(230.61,258.23)}, rotate = 315] [fill={rgb, 255:red, 0; green, 0; blue, 0 }  ][line width=0.08]  [draw opacity=0] (8.93,-4.29) -- (0,0) -- (8.93,4.29) -- cycle    ;
\draw [shift={(263.12,225.72)}, rotate = 135] [color={rgb, 255:red, 0; green, 0; blue, 0 }  ][line width=0.75]      (0, 0) circle [x radius= 3.35, y radius= 3.35]   ;
\draw    (434.92,227.38) -- (380.29,282.01) ;
\draw [shift={(378.63,283.67)}, rotate = 135] [color={rgb, 255:red, 0; green, 0; blue, 0 }  ][line width=0.75]      (0, 0) circle [x radius= 3.35, y radius= 3.35]   ;
\draw [shift={(404.07,258.23)}, rotate = 315] [fill={rgb, 255:red, 0; green, 0; blue, 0 }  ][line width=0.08]  [draw opacity=0] (8.93,-4.29) -- (0,0) -- (8.93,4.29) -- cycle    ;
\draw [shift={(436.58,225.72)}, rotate = 135] [color={rgb, 255:red, 0; green, 0; blue, 0 }  ][line width=0.75]      (0, 0) circle [x radius= 3.35, y radius= 3.35]   ;
\draw    (377.1,227.38) -- (322.47,282.01) ;
\draw [shift={(320.81,283.67)}, rotate = 135] [color={rgb, 255:red, 0; green, 0; blue, 0 }  ][line width=0.75]      (0, 0) circle [x radius= 3.35, y radius= 3.35]   ;
\draw [shift={(346.25,258.23)}, rotate = 315] [fill={rgb, 255:red, 0; green, 0; blue, 0 }  ][line width=0.08]  [draw opacity=0] (8.93,-4.29) -- (0,0) -- (8.93,4.29) -- cycle    ;
\draw [shift={(378.76,225.72)}, rotate = 135] [color={rgb, 255:red, 0; green, 0; blue, 0 }  ][line width=0.75]      (0, 0) circle [x radius= 3.35, y radius= 3.35]   ;
\draw    (319.28,227.38) -- (264.65,282.01) ;
\draw [shift={(262.99,283.67)}, rotate = 135] [color={rgb, 255:red, 0; green, 0; blue, 0 }  ][line width=0.75]      (0, 0) circle [x radius= 3.35, y radius= 3.35]   ;
\draw [shift={(288.43,258.23)}, rotate = 315] [fill={rgb, 255:red, 0; green, 0; blue, 0 }  ][line width=0.08]  [draw opacity=0] (8.93,-4.29) -- (0,0) -- (8.93,4.29) -- cycle    ;
\draw [shift={(320.94,225.72)}, rotate = 135] [color={rgb, 255:red, 0; green, 0; blue, 0 }  ][line width=0.75]      (0, 0) circle [x radius= 3.35, y radius= 3.35]   ;
\draw    (376.28,283.67) -- (323.16,283.67) ;
\draw [shift={(320.81,283.67)}, rotate = 180] [color={rgb, 255:red, 0; green, 0; blue, 0 }  ][line width=0.75]      (0, 0) circle [x radius= 3.35, y radius= 3.35]   ;
\draw [shift={(344.72,283.67)}, rotate = 360] [fill={rgb, 255:red, 0; green, 0; blue, 0 }  ][line width=0.08]  [draw opacity=0] (8.93,-4.29) -- (0,0) -- (8.93,4.29) -- cycle    ;
\draw [shift={(378.63,283.67)}, rotate = 180] [color={rgb, 255:red, 0; green, 0; blue, 0 }  ][line width=0.75]      (0, 0) circle [x radius= 3.35, y radius= 3.35]   ;
\draw    (318.46,283.67) -- (265.34,283.67) ;
\draw [shift={(262.99,283.67)}, rotate = 180] [color={rgb, 255:red, 0; green, 0; blue, 0 }  ][line width=0.75]      (0, 0) circle [x radius= 3.35, y radius= 3.35]   ;
\draw [shift={(286.9,283.67)}, rotate = 360] [fill={rgb, 255:red, 0; green, 0; blue, 0 }  ][line width=0.08]  [draw opacity=0] (8.93,-4.29) -- (0,0) -- (8.93,4.29) -- cycle    ;
\draw [shift={(320.81,283.67)}, rotate = 180] [color={rgb, 255:red, 0; green, 0; blue, 0 }  ][line width=0.75]      (0, 0) circle [x radius= 3.35, y radius= 3.35]   ;
\draw    (434.1,283.67) -- (380.98,283.67) ;
\draw [shift={(378.63,283.67)}, rotate = 180] [color={rgb, 255:red, 0; green, 0; blue, 0 }  ][line width=0.75]      (0, 0) circle [x radius= 3.35, y radius= 3.35]   ;
\draw [shift={(402.54,283.67)}, rotate = 360] [fill={rgb, 255:red, 0; green, 0; blue, 0 }  ][line width=0.08]  [draw opacity=0] (8.93,-4.29) -- (0,0) -- (8.93,4.29) -- cycle    ;
\draw [shift={(436.45,283.67)}, rotate = 180] [color={rgb, 255:red, 0; green, 0; blue, 0 }  ][line width=0.75]      (0, 0) circle [x radius= 3.35, y radius= 3.35]   ;
\draw    (260.64,283.67) -- (207.52,283.67) ;
\draw [shift={(205.17,283.67)}, rotate = 180] [color={rgb, 255:red, 0; green, 0; blue, 0 }  ][line width=0.75]      (0, 0) circle [x radius= 3.35, y radius= 3.35]   ;
\draw [shift={(229.08,283.67)}, rotate = 360] [fill={rgb, 255:red, 0; green, 0; blue, 0 }  ][line width=0.08]  [draw opacity=0] (8.93,-4.29) -- (0,0) -- (8.93,4.29) -- cycle    ;

\draw [shift={(262.99,283.67)}, rotate = 180] [color={rgb, 255:red, 0; green, 0; blue, 0 }  ][line width=0.75]      (0, 0) circle [x radius= 3.35, y radius= 3.35]   ;
\draw (309.77,34.98) node [anchor=north west][inner sep=0.75pt]    {$1$};
\draw (382.37,94.83) node [anchor=north west][inner sep=0.75pt]    {$1$};
\draw (309.27,93.62) node [anchor=north west][inner sep=0.75pt]    {$0$};
\draw (378.48,149.54) node [anchor=north west][inner sep=0.75pt]    {$0$};
\draw (313.86,149.03) node [anchor=north west][inner sep=0.75pt]    {$1$};
\draw (250.4,151.74) node [anchor=north west][inner sep=0.75pt]    {$1$};
\draw (319.02,206.56) node [anchor=north west][inner sep=0.75pt]    {$1$};
\draw (251.7,209.47) node [anchor=north west][inner sep=0.75pt]    {$0$};
\draw (380.27,207.57) node [anchor=north west][inner sep=0.75pt]    {$2$};
\draw (439.94,208.49) node [anchor=north west][inner sep=0.75pt]    {$2$};
\draw (437.34,266.1) node [anchor=north west][inner sep=0.75pt]    {$0$};
\draw (371.52,264.4) node [anchor=north west][inner sep=0.75pt]    {$2$};
\draw (314.77,264.3) node [anchor=north west][inner sep=0.75pt]    {$4$};
\draw (255.89,265.5) node [anchor=north west][inner sep=0.75pt]    {$5$};
\draw (193.52,268.1) node [anchor=north west][inner sep=0.75pt]    {$5$};

\end{tikzpicture}
\caption {The \emph{boustrophedon} number triangle for the Entringer numbers.}
\label {fig:boustrophedon}
\end{figure}
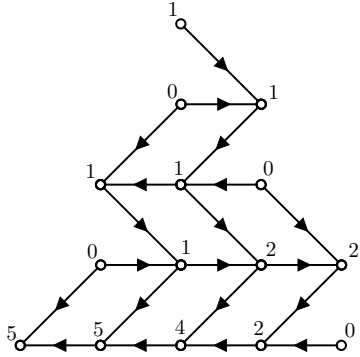

One way to construct $P_{n}^{s}$ is as follows. Consider the directed graph in Figure \ref{fig:boustrophedon}, which encodes the recurrence for the Entringer numbers.
One starts by labeling the sources of the graph from top to bottom (note that there is only one source in each row) by $1,0,0,...$. 
Then, the labels for the rest of the nodes are computed by summing
all incoming labels. Finally, reading from top to bottom the ``ends" of each horizontal row, 
we obtain another sequence, called the \emph{boustrophedon transform} of the original. In this case, we see that the boustrophedon transform of the sequence $(1,0,0,0,0,...)$ is the sequence of Euler numbers $(1,1,2,5,16,...)$.
See \cite{boustrophedon} for more on the boustrophedon transform. 

\begin{figure}
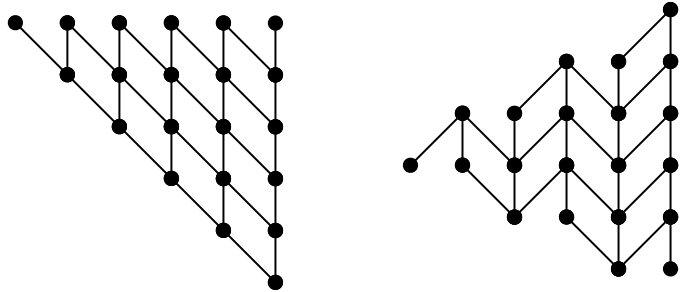

\centering
\tikzset{every picture/.style={line width=0.75pt}} 

\caption {The posets $P_n^z$ (left) and $P_n^s$ (right) such that $\Omega(\mathcal{G}_n^z) \cong J(P_n^z)$ and $\Omega(\mathcal{G}_n^s) \cong J(P_n^s)$.}
\label {fig:underlying_posets}
\end{figure}

The Hasse diagram of the poset $P_{n}^{s}$ is obtained by rotating the underlying undirected graph in Figure \ref{fig:boustrophedon} counterclockwise by $90^{\circ}$ (see the poset on the right side of Figure \ref{fig:underlying_posets}). To see that this is true, notice that the poset $P_{n}^{s}$ naturally fits into the plane partition picture shown in Figure \ref{fig:rhombus_tilings} in the following way. Each node in the poset is mapped to the center of its respective cube in the maximal plane partition, and the covering relations between nodes in $P_{n}^{s}$ point between adjacent cube centers in the three positive directions in the first octant.

The poset $P_{n}^{z}$ is isomorphic to the type $A$ root poset (see e.g. \cite{matchable}). See Figure \ref{fig:underlying_posets} (left). There are similar relationships between the triangle of ballot numbers, the rhombus tilings corresponding to $\mathcal{C}_n$, and the poset $P_{n}^{z}$. 

\section {General Snake Graphs}

\subsection {Standard Labelings}

So far we have considered only straight or zig-zag snake graphs. 
We will now consider the general case, and define the \emph{standard labeling} of any snake graph, 
which generalizes the standard labelings we have seen for straight and zig-zag shapes. 
Some constructions and notation will be introduced in preparation of the general definition.
Then we study mixed dimer configurations on general snake graphs equipped with their standard labeling. 

\begin {defn} 
    A \emph{snake graph} $\mathcal{G}$ is a connected planar graph formed by recursively gluing finitely many squares, each either above or to the right of the previous one. Each square in this construction is referred to as a \emph{tile}. We will denote the tiles by $T_i$, in order from bottom-left to top-right. 
\end {defn}

See \cite{claussen} for more precise definitions of snake graphs. Snake graphs have come to play an important role 
in the theory of cluster algebras \cite{ms_09} \cite{msw_13}, 
and have connections to number theory \cite{snake_cluster} \cite{snake_markov}, knot theory \cite{jones_poly}, and other fields.

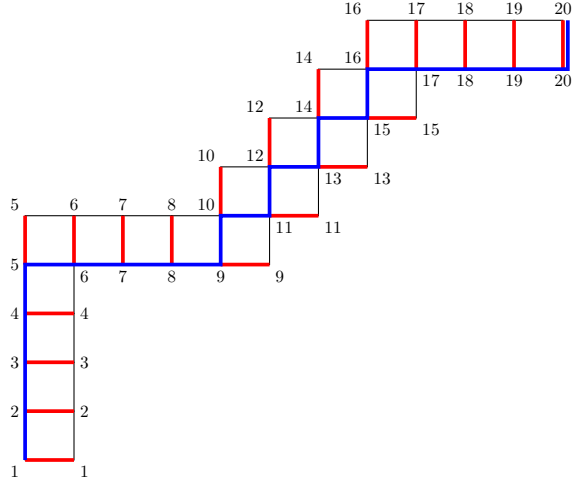
\begin {figure}
\centering
\begin {tikzpicture} [scale=0.65, every node/.style={scale=0.65}]
    \foreach \x/\y in {0/0, 0/1, 0/2, 0/3, 0/4, 1/4, 2/4, 3/4, 4/4, 4/5, 5/5, 5/6, 6/6, 6/7, 7/7, 7/8, 8/8, 9/8, 10/8} {
        \draw [shift={(\x,\y)}] (0,0) -- (1,0) -- (1,1) -- (0,1) -- cycle;
    }

    \foreach \x/\y in {0/0, 0/1, 0/2, 0/3, 4/4, 5/5, 6/6, 7/7} {
        \draw[shift={(\x,\y)}, red, line width=1.4] (0,0) -- (1,0);
    }

    \foreach \x/\y in {0/4, 1/4, 2/4, 3/4, 4/5, 5/6, 6/7, 7/8, 8/8, 9/8, 10/8, 11/8} {
        \draw[shift={(\x,\y)}, red, line width=1.4] (0,0) -- (0,1);
    }

    \draw[blue, line width=1.4] (0,0) --++ (0,4) --++ (4,0) --++ (0,1) --++ (1,0) --++ (0,1) --++ (1,0) --++ (0,1) --++ (1,0) --++ (0,1) --++ (4.1,0) --++ (0,1);

    \draw (0,0) node[below left]  {$1$};
    \draw (1,0) node[below right] {$1$};
    \draw (0,1) node[left]        {$2$};
    \draw (1,1) node[right]       {$2$};
    \draw (0,2) node[left]        {$3$};
    \draw (1,2) node[right]       {$3$};
    \draw (0,3) node[left]        {$4$};
    \draw (1,3) node[right]       {$4$};
    \draw (0,4) node[left]        {$5$};
    \draw (0,5) node[above left]  {$5$};
    \draw (1,4) node[below right] {$6$};
    \draw (1,5) node[above]       {$6$};
    \draw (2,4) node[below]       {$7$};
    \draw (2,5) node[above]       {$7$};
    \draw (3,4) node[below]       {$8$};
    \draw (3,5) node[above]       {$8$};
    \draw (4,4) node[below]       {$9$};
    \draw (5,4) node[below right] {$9$};
    
    \draw (4,5) node[above left]  {$10$};
    \draw (4,6) node[above left]  {$10$};
    \draw (5,6) node[above left]  {$12$};
    \draw (5,7) node[above left]  {$12$};
    \draw (6,7) node[above left]  {$14$};
    \draw (6,8) node[above left]  {$14$};
    \draw (7,8) node[above left]  {$16$};
    \draw (7,9) node[above left]  {$16$};
    
    \draw (5,5) node[below right] {$11$};
    \draw (6,5) node[below right] {$11$};
    \draw (6,6) node[below right] {$13$};
    \draw (7,6) node[below right] {$13$};
    \draw (7,7) node[below right] {$15$};
    \draw (8,7) node[below right] {$15$};
    
    \draw (8,8) node[below right] {$17$};
    \draw (8,9) node[above]       {$17$};
    \draw (9,8) node[below]       {$18$};
    \draw (9,9) node[above]       {$18$};
    \draw (10,8) node[below]      {$19$};
    \draw (10,9) node[above]      {$19$};
    \draw (11,8) node[below]      {$20$};
    \draw (11,9) node[above]      {$20$};
\end {tikzpicture}
\caption {A snake graph $\mathcal{G}^{w}$ with word $w = U^4R^4URURURUR^3$, and its standard labeling. The canonical dimer cover $D_0$ is pictured in red. The canonical lattice path $P_0$ is in blue.}
\label {fig:canonical_dimer}
\end {figure}

The shape of any snake graph can be encoded by its underlying \emph{word} $w \in \{ R, U \}^{*}$, where $\{ R, U \}^{*}$ is the set of all finite words $w$ 
in the alphabet $\{R,U\}$ (including the the empty word). We will write $\mathcal{G}^{w}_{n}$ for the snake graph with $n$ tiles and word $w$ of length $n-1$. 
We will sometimes just write $\mathcal{G}^w$ instead of $\mathcal{G}^w_n$ (omitting the subscript).
For instance, letting $s = RRRR ...$ and $z = URUR ...$ recovers the straight and zig-zag snake graphs $\mathcal{G}^{s}_{n}$ and $\mathcal{G}^{z}_{n}$, respectively. We let the empty word correspond to the unique snake graph with one tile. An example is pictured in Figure \ref{fig:canonical_dimer}.

\begin {defn}
Let $\mathcal{G} = (V,E)$ be a finite graph, and let $D \in \Omega_1(\mathcal{G})$ be a dimer cover (perfect matching) of $\mathcal{G}$. 
If $\n_{D} : V \longrightarrow \Bbb{N}$ is such that 
$\n_{D}(v) = \n_{D}(u)$ whenever $v$ and $u$ are the endpoints of an edge in $D$, then we say $\n_{D}$ is a \emph{D-labeling} of $\mathcal{G}$.  
\end {defn}

Clearly any $D$-labeling admits mixed dimer covers (that is, $\Omega_{\n}(\mathcal{G})$ is non-empty if $\n$ is a $D$-labeling).
We will define the \emph{standard labeling} of any snake graph as a certain $D$-labeling. To do so, we must make a canonical choice
of dimer cover for any snake graph.

\begin {defn} 
    Let $\mathcal{G}$ be a snake graph with word $w$. We define a dimer cover $D_0$, which we call the \emph{canonical dimer cover}, as follows. 
    Whenever $w$ contains the substring $UR$,
    then $D_0$ contains the left vertical edge of the corresponding corner tile. Whenever $w$ contains the substring $RU$, $D_0$ contains the bottom
    horizontal edge of the corner tile. Otherwise, $D_0$ contains all internal edges (i.e. those not adjacent to the outer face) of the straight
    segments of $\mathcal{G}$. Figure \ref{fig:canonical_dimer} shows the canonical dimer cover in red.
\end {defn}

\begin {rmk}
    This definition leaves some ambiguity; namely whether the last tile is considered a corner tile or as part of a straight segment. 
    We resolve this ambiguity with the following rule. Look at the last two letters of the word $w$ (equivalently the last three tiles of the graph).
    If the last two letters are the same (i.e. the last 3 tiles form a straight segment), then treat the final tile as if it were part of a straight segment.
    If the last two letters are different (i.e. the last 3 tiles form a zigzag), then treat the final tile as if it were a corner tile.
\end {rmk}

\begin {defn}
    Let $D_0$ be the canonical dimer cover of a snake graph $\mathcal{G}^w_n$. Define a linear order on the edges of $D_0$ as follows. 
    If $e_1$ is an edge of tile $T_i$ and $e_2$ is an edge of $T_j$ with $i < j$, then order $e_i < e_j$.
    If $e_i$ and $e_j$ are both edges of the same tile $T_i$, then order them left-to-right or bottom-to-top.
    With the edges of $D_0$ ordered $e_1 < e_2 < \cdots < e_{n+1}$, we define the \emph{standard labeling} of $\mathcal{G}^w_n$ as the $D_0$-labeling 
    where $\n(u) = \n(v) = k$ for edge $e_k = (u,v)$. An example is pictured in Figure \ref{fig:canonical_dimer}.
\end {defn}

\bigskip

\subsection {The Dimers-Paths Duality} \label{sec:dimers_paths}

There is an involution on the set of all finite snake graphs, first defined in \cite{propp_20}. 
See \cite{claussen} for more on this involution and its manifestations on some of the objects associated with snake graphs, 
such as triangulations and continued fractions. 
In this section, we will define the duality, and explain how it is related to the enumeration of lattice paths, and finally extend this
to the current case of mixed dimer covers.

Denote the nontrivial involution on $\{ R, U \}$ by $\widetilde{R} = U$ and $\widetilde{U} = R$. Given any word $w = w_1 w_2 \cdots w_n \in \{ R, U \}^*,$ 
let $\widetilde{w} = \widetilde{w}_1 w_{2} \widetilde{w}_{3} \cdots w_{n-1} \widetilde{w}_{n}$ if $n$ is odd, 
or $\widetilde{w} = \widetilde{w}_{1} w_{2} \widetilde{w}_{3} \cdots \widetilde{w}_{n-1} w_{n}$ if $n$ is even.

\begin {defn} \label {def:duality}
    If $\mathcal{G} = \mathcal{G}^w$ is the snake graph with word $w$, we define the \emph{dual snake graph} to be 
    $\widetilde{\mathcal{G}} = \mathcal{G}^{\widetilde{w}}.$ 
\end {defn}

\begin {ex}
    The dual of the straight snake graph $\mathcal{G}^{s}_n$ with word $w = RRRR\dots$ is the zig-zag snake graph $\mathcal{G}^z_n$, where $\widetilde{w} = URUR \dots$
\end {ex}

\begin {figure}
\centering
\begin{tikzpicture}[x=0.75pt,y=0.75pt,yscale=-1,xscale=1]
 
\draw    (240.56,219.62) -- (290.95,219.62) ;
\draw    (240.56,169.23) -- (240.56,219.62) ;
\draw    (240.56,169.23) -- (290.95,169.23) ; 
\draw    (290.95,169.23) -- (290.95,219.62) ;
\draw    (290.95,118.84) -- (290.95,169.23) ; 
\draw    (240.56,118.84) -- (240.56,169.23) ;
\draw    (240.56,118.84) -- (290.95,118.84) ; 
\draw    (240.56,68.45) -- (240.56,118.84) ; 
\draw    (290.95,68.45) -- (290.95,118.84) ; 
\draw    (240.56,68.45) -- (290.95,68.45) ; 
\draw    (29.54,219.62) -- (79.93,219.62) ; 
\draw    (29.54,169.23) -- (29.54,219.62) ; 
\draw    (29.54,169.23) -- (79.93,169.23) ; 
\draw    (79.93,169.23) -- (79.93,219.62) ; 
\draw    (79.93,219.62) -- (130.31,219.62) ; 
\draw    (79.93,169.23) -- (130.31,169.23) ;
\draw    (130.31,169.23) -- (130.31,219.62) ;
\draw    (130.31,169.23) -- (180.7,169.23) ;
\draw    (130.31,219.62) -- (180.7,219.62) ;
\draw    (180.7,169.23) -- (180.7,219.62) ; 
\draw    (359.99,219.62) -- (410.38,219.62) ; 
\draw    (359.99,169.23) -- (359.99,219.62) ;
\draw    (359.99,169.23) -- (410.38,169.23) ;
\draw    (410.38,169.23) -- (410.38,219.62) ;
\draw    (359.99,118.84) -- (410.38,118.84) ;
\draw    (410.38,169.23) -- (460.77,169.23) ;
\draw    (359.99,118.84) -- (359.99,169.23) ;
\draw    (410.38,118.84) -- (410.38,169.23) ;
\draw    (460.77,118.84) -- (460.77,169.23) ;
\draw    (410.38,118.84) -- (460.77,118.84) ;

\draw  [dash pattern={on 0.84pt off 2.51pt}]  (222.76,101.43) -- (306.74,185.41) ;
\draw  [dash pattern={on 0.84pt off 2.51pt}]  (14.13,154.43) -- (98.11,238.41) ;
\draw  [dash pattern={on 0.84pt off 2.51pt}]  (393.59,102.05) -- (477.57,186.03) ;

\draw    (152,242.25) .. controls (176.16,259.37) and (209.08,261.3) .. (236.78,243.09) ;
\draw [shift={(238.05,242.25)}, rotate = 145.66] [color={rgb, 255:red, 0; green, 0; blue, 0 }  ][line width=0.75]    (10.93,-3.29) .. controls (6.95,-1.4) and (3.31,-0.3) .. (0,0) .. controls (3.31,0.3) and (6.95,1.4) .. (10.93,3.29)   ;

\draw    (282.29,242.25) .. controls (306.45,259.37) and (339.37,261.3) .. (367.08,243.09) ;
\draw [shift={(368.34,242.25)}, rotate = 145.66] [color={rgb, 255:red, 0; green, 0; blue, 0 }  ][line width=0.75]    (10.93,-3.29) .. controls (6.95,-1.4) and (3.31,-0.3) .. (0,0) .. controls (3.31,0.3) and (6.95,1.4) .. (10.93,3.29)   ;

\draw    (410.99,242.2) .. controls (435.15,259.32) and (468.07,261.25) .. (495.77,243.05) ;
\draw [shift={(497.04,242.2)}, rotate = 145.66] [color={rgb, 255:red, 0; green, 0; blue, 0 }  ][line width=0.75]    (10.93,-3.29) .. controls (6.95,-1.4) and (3.31,-0.3) .. (0,0) .. controls (3.31,0.3) and (6.95,1.4) .. (10.93,3.29)   ;

\draw    (521.33,219.57) -- (571.72,219.57) ;
\draw    (521.33,169.18) -- (521.33,219.57) ; 
\draw    (521.33,169.18) -- (571.72,169.18) ; 
\draw    (571.72,169.18) -- (571.72,219.57) ;
\draw    (521.33,118.8) -- (571.72,118.8) ;
\draw    (571.72,169.18) -- (622.1,169.18) ; 
\draw    (521.33,118.8) -- (521.33,169.18) ;
\draw    (571.72,118.8) -- (571.72,169.18) ;
\draw    (622.1,118.8) -- (622.1,169.18) ; 
\draw    (571.72,118.8) -- (622.1,118.8) ;


\draw (79.75,162) node  [font=\scriptsize]  {$2$};
\draw (29.54,225) node  [font=\scriptsize]  {$1$};
\draw (130.1,162) node  [font=\scriptsize]  {$3$};
\draw (179.1,162) node  [font=\scriptsize]  {$4$};
\draw (79.78,225.02) node  [font=\scriptsize, red]  {$2$};
\draw (29.5,162) node  [font=\scriptsize, red]  {$1$};
\draw (130.1,224.97) node  [font=\scriptsize]  {$3$};
\draw (179.14,225) node  [font=\scriptsize]  {$4$};
\draw (54.66,226) node  {$p$};
\draw (104.94,226) node  {$q$};
\draw (155.22,226) node  {$r$};
\draw (54.66,162) node  {$x$};
\draw (104.94,162) node  {$y$};
\draw (155.22,162) node  {$z$};
\draw (24,193) node {$a$};
\draw (74.28,193) node {$b$};
\draw (124.56,193) node {$c$};
\draw (174.84,193) node {$d$};

\draw (238,225) node  [font=\scriptsize]  {$1$};
\draw (293.9,225) node  [font=\scriptsize]  {$1$};
\draw (295.75,169.2) node  [font=\scriptsize, red]  {$2$};
\draw (236,169.23) node  [font=\scriptsize]  {$2$};
\draw (236,117.34) node  [font=\scriptsize, red]  {$3$};
\draw (295.06,117.3) node  [font=\scriptsize]  {$3$};
\draw (236,67.8) node  [font=\scriptsize]  {$4$};
\draw (295.31,67.84) node  [font=\scriptsize]  {$4$};
\draw (265.95,226) node {$p$};
\draw (232.5,193) node {$a$};
\draw (232.5,142.72) node {$q$};
\draw (232.5,92.44) node {$r$};
\draw (265.95,175.72) node {$b$};
\draw (265.95,125.44) node {$c$};
\draw (265.95,75.16) node {$d$};
\draw (297,193) node {$x$};
\draw (297,142.72) node {$y$};
\draw (297,92.44) node {$z$};

\draw (357.2,225) node  [font=\scriptsize]  {$1$};
\draw (413.2,225) node  [font=\scriptsize]  {$1$};
\draw (355.03,169.2) node  [font=\scriptsize]  {$2$};
\draw (355,117.3) node  [font=\scriptsize]  {$2$};
\draw (414.36,174.5) node  [font=\scriptsize]  {$3$};
\draw (409.56,111.5) node  [font=\scriptsize, red]  {$3$};
\draw (465.16,115.04) node  [font=\scriptsize]  {$4$};
\draw (465.2,169.2) node  [font=\scriptsize, red]  {$4$};
\draw (385.2,226) node {$p$};
\draw (351.75,193) node {$a$};
\draw (351.75,142.72) node {$q$};
\draw (416.25,193) node {$x$};
\draw (385.2,175.72) node {$b$};
\draw (385.2,112) node {$y$};
\draw (435.48,112) node {$z$};
\draw (435.48,175.72) node {$r$};
\draw (402.03,142.72) node {$c$};
\draw (466.53,142.72) node {$d$};

\draw (518.53,224.95) node  [font=\scriptsize]  {$1$};
\draw (574.53,224.95) node  [font=\scriptsize]  {$1$};
\draw (516.36,169.15) node  [font=\scriptsize]  {$2$};
\draw (516.33,117.25) node  [font=\scriptsize]  {$2$};
\draw (575.69,174.45) node  [font=\scriptsize]  {$3$};
\draw (570.89,112.45) node  [font=\scriptsize, red]  {$4$};
\draw (626.49,115) node  [font=\scriptsize]  {$4$};
\draw (626.53,169.15) node  [font=\scriptsize, red]  {$3$};
\begin {scope}[shift={(161.33,0)}]
\draw (385.2,226) node {$p$};
\draw (351.75,193) node {$a$};
\draw (351.75,142.72) node {$q$};
\draw (416.25,193) node {$x$};
\draw (385.2,175.72) node {$b$};
\draw (385.2,112) node {$y$};
\draw (435.48,112) node {$d$};
\draw (435.48,175.72) node {$r$};
\draw (402.03,142.72) node {$c$};
\draw (466.53,142.72) node {$z$};
\end {scope}

\draw (186.67,260.07) node [anchor=north west][inner sep=0.75pt] {$\tau_1$};
\draw (317,260.07)    node [anchor=north west][inner sep=0.75pt] {$\tau_2$};
\draw (448,260.1)     node [anchor=north west][inner sep=0.75pt] {$\tau_3$};
\end{tikzpicture}
\caption {The snake graph duality $\mathcal{G} \mapsto \widetilde{\mathcal{G}}$, realized as the composition of tile maps $\tau_3 \circ \tau_2 \circ \tau_1$.}
\label {fig:tile_maps}
\end {figure}
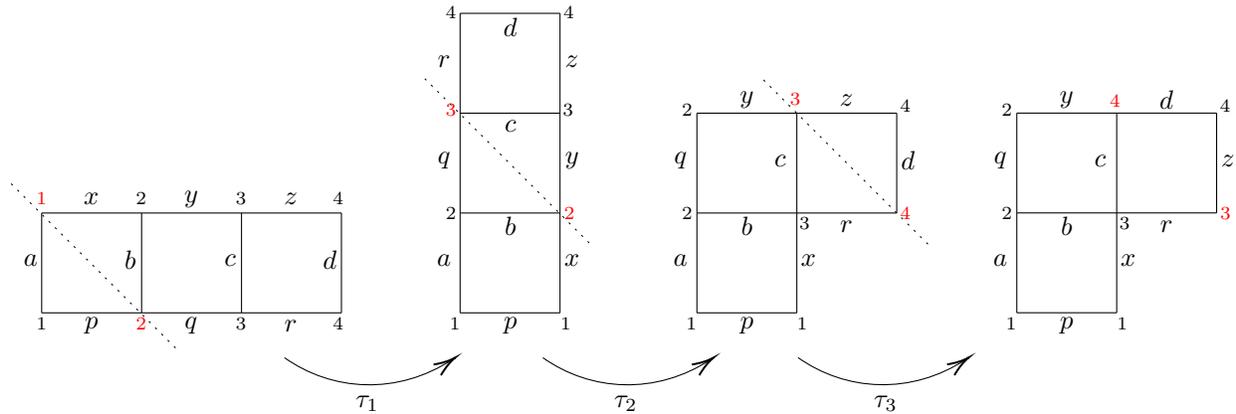

There is a more refined version of the duality $\mathcal{G} \mapsto \widetilde{\mathcal{G}}$
which gives a bijection between the edges of $\mathcal{G}$ and the edges of $\widetilde{\mathcal{G}}$. So if $\mathcal{G}$ has edge weights,
it induces a weighting of the edges of its dual $\widetilde{\mathcal{G}}$.
For a snake graph with $n$ tiles, we will describe the dual map as a composition of $n$ \emph{tile maps} (as in \cite{claussen}). 
For $1 \leq i \leq n$, let $\ell_-$ and $\ell_+$ be the diagonals of tile $T_i$ with slope $-1$ and $+1$.
Define the tile map $\tau_i$ (mapping the set of edge-labeled snake graphs to itself) by the following process:
\begin {itemize}
    \item Cut the snake graph into two pieces along the line $\ell_-$ 
    \item Reflect the right half over the line $\ell_+$ (leaving the left half intact) 
    \item Glue the two resulting pieces back together 
\end {itemize}
The dual map from Definition \ref{def:duality} can then be written as the composition $\tau = \tau_n \circ \tau_{n-1} \circ \cdots \circ \tau_2 \circ \tau_1.$ 
Following this procedure maps each edge of $\mathcal{G}$ bijectively to an edge of $\widetilde{\mathcal{G}}$.
See Figure \ref{fig:tile_maps} for an illustration. Each $\tau_i$ is clearly an involution, and it is not hard to see that the composition $\tau$ is as well.

Roughly, a \emph{mixed lattice path} on a snake graph $\mathcal{G}_{n}^{w}$ with its standard labeling is a collection of $n+1$ north-east paths $L_0,L_1, L_2, ... , L_n$ on $\mathcal{G}_{n}^{w}$ with distinct lengths in $\{1, 2, 3, ... , n+1 \}$ that satisfy certain restrictions. Mixed lattice paths are in a sense dual to mixed dimer covers, which we will now explain. 

We generalize the dual map to any snake graph equipped with a vertex labeling by declaring that the tile map $\tau_{i}$ switches the two labels  
on the upper-left and lower-right vertices of tile $T_{i}$. All other labels unambiguously remain attached to their corresponding vertices. 
See Figure \ref{fig:tile_maps} for an example of how the vertex labels change under the tile maps.

\begin {defn}
    Let $\mathcal{G} = \mathcal{G}^w$ be a snake graph with word $w$. 
    The \emph{canonical lattice path} $P_0$ is a north-east lattice path on $\mathcal{G}$ from the bottom-left to the top-right
    vertex whose sequence of up and right steps is given by the word $w$. As with the definition of the canonical dimer cover,
    this leaves some ambiguity on the last tile, which we resolve in the same way: treat the last tile as a corner or part of a straight segment based on the
    pattern formed by the final three squares. Equivalently, it is the unique lattice path which has no edges in common with the canonical dimer cover (except on the final square). See Figure \ref{fig:canonical_dimer} for an example.
\end {defn}

\begin {defn}
    Let $\mathcal{G}^w_n$ be a snake graph with its standard labeling. A \emph{mixed lattice path} is a multiset of edges which can be realized
    (not necessarily uniquely) as the union of $n+1$ lattice paths $L_0, L_1, L_2, \dots, L_n$ such that
    \begin {itemize}
        \item All paths $L_i$ begin at one of the vertices along the canonical path $P_0$, and end at the top-right vertex of $\mathcal{G}$.
        \item The path $L_i$ has length $i+1$.
    \end {itemize}
\end {defn}
 
Note that these conditions uniquely determine the starting vertex for each path. If the vertices along the canonical path are $v_0,v_1,v_2,\dots$,
then path $L_i$ begins at vertex $v_{n-i}$.

Much like for ordinary single lattice paths, mixed lattice paths form a poset under flips 
(i.e. if we represent a lattice path by a word from $\{ R,U \}^{*}$ in the obvious way, then a flip replaces an instance or $RU$ with $UR$). 

In \cite{propp_02} and \cite{claussen}, it was shown that the duality given by the composition of tile maps induces a bijection between
dimer covers of $\mathcal{G}$ and north-east lattice paths on $\widetilde{\mathcal{G}}$. We now give a generalization of this fact
to the current setting.

\begin {thm} \label{thm:duality}
    Let $\mathcal{G} = \mathcal{G}^w_n$ be a snake graph with word $w$, and $\widetilde{\mathcal{G}} = \mathcal{G}^{\tilde{w}}_n$ its dual.
    \begin {enumerate}
        \item[(a)] The duality map $\tau$ sends the standard labeling of $\mathcal{G}$ to the standard labeling of $\widetilde{\mathcal{G}}$.
        \item[(b)] The duality map $\tau$ sends the canonical dimer cover of $\mathcal{G}$ to the canonical lattice path of $\widetilde{\mathcal{G}}$.
        \item[(c)] The duality induces a bijection between mixed dimer covers of $\mathcal{G}$ and mixed lattice paths on $\widetilde{\mathcal{G}}$.
        \item[(d)] The bijection from part $(c)$ is a poset isomorphism. 
    \end {enumerate}
\end {thm}
\begin {proof}
    Note that parts $(a)$ and $(b)$ are equivalent, since the canonical dimer cover and the canonical lattice path are deterimined
    by the standard labeling, and vice versa.
    So we will just prove $(b)$, and $(a)$ will follow. 
    
    Part $(b)$ can be proved by looking locally at three consecutive tiles of $\mathcal{G}$, and
    breaking into several different cases. For example, if three tiles of $\mathcal{G}$ go right and then up, then $D_0$ has an edge on the bottom 
    of the corner tile. The corresponding three tiles of $\widetilde{\mathcal{G}}$ form a straight segment. 
    If the corner tile is even (i.e. $T_i = T_{2k}$), these three tiles in the dual are a vertical straight segment, and the bottom edge
    of $T_i$ maps to the left edge in the dual. Similarly, if the corner tile is odd then these three tiles in the dual form a horizontal segment,
    and the bottom edge of $T_i$ maps to the bottom edge of the dual tile. One easily sees from the definition of the canonical path $P_0$ that 
    in both cases the bottom edge of tile $T_i$ maps to an edge of $P_0$.

    There are several other cases to consider, namely if tile $T_i$ is a northwest corner (rather than southeast), and if $T_i$ is part of a straight segment.
    All these cases will similarly depend on the parity of the tile in question. The analysis of these cases is similar (and straightforward), and we omit the details.
    
    $(c)$ A mixed dimer cover may be realized as the union of partial dimer covers $M_0, M_1, \dots M_n$, where $M_i$ is a dimer cover of the subgraph
    obtained by removing the first $n-i$ edges of the canonical dimer cover from $\mathcal{G}$. Under the duality (from \cite{propp_02} and \cite{claussen}), 
    these partial dimer covers are mapped to
    partial lattice paths. By parts $(a)$ and $(b)$, the edges of the canonical dimer cover are mapped to the edges along the canonical lattice path,
    and so the earliest edge in the partial dimer cover $M_i$ maps to the first edge in a partial lattice path beginning at the appropriate vertex
    along the canonical path.

    $(d)$ We need to see that covering relations on mixed dimer covers of $\mathcal{G}$ are mapped via the duality $\tau$ to covering relations
    of mixed lattice paths on $\widetilde{\mathcal{G}}$. This follows immediately from the case of single dimer covers, since the covering relations
    are the same in the mixed case.
\end {proof}

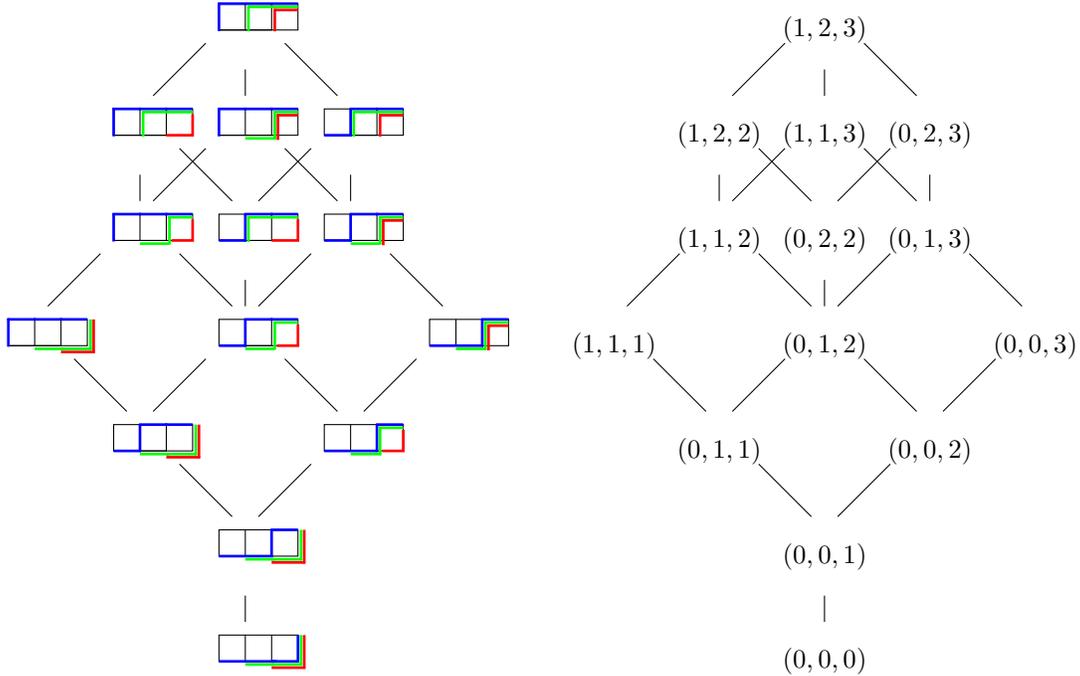
\begin {figure}
\begin{center}
\begin {tikzpicture}[scale=0.35]
    \foreach \x/\y in {0/0, 0/1, -1/2, 1/2, -2/3, 0/3, 2/3, -1/4, 0/4, 1/4, -1/5, 0/5, 1/5, 0/6} {
        \begin {scope}[shift={(4*\x,4*\y)}]
            \draw (0,1) -- (1,1) -- (2,1) -- (3,1) -- (3,2) -- (2,2) -- (1,2) -- (0,2) -- cycle;
            \draw (0,1) -- (1,1) -- (1,2);
            \draw (1,1) -- (2,1) -- (2,2);
        \end {scope}
    }

    \foreach \x/\y in {0/0, 0/3, -1/4, 1/4, 0/5} {
        \begin {scope}[shift={(4*\x,4*\y)}]
            \draw (1,2.5) -- (1,3.5);
        \end {scope}
    }

    \foreach \x/\y in {0/1, -1/2, 1/2, -2/3, 0/3, -1/4, 0/4, -1/5} {
        \begin {scope}[shift={(4*\x,4*\y)}]
            \draw (1.5,2.5) -- (3.5,4.5);
        \end {scope}
    }

    \foreach \x/\y in {0/1, -1/2, 1/2, 0/3, 2/3, 0/4, 1/4, 1/5} {
        \begin {scope}[shift={(4*\x,4*\y)}]
            \draw (0.5,2.5) -- (-1.5,4.5);
        \end {scope}
    }

    \begin {scope}[shift={(4*0, 4*0)}]
        \draw[blue, line width=1] (0,1) -- (1,1);
        \draw[blue, line width=1] (1,1) -- (2.2,1);
        \draw[blue, line width=1] (2,1) -- (3.04,1);
        \draw[blue, line width=1] (3,.96) -- (3,2);
        \draw[green, line width=1] (1,.88) -- (2.1,.88);
        \draw[green, line width=1] (2,0.88) -- (3.16,0.88);
        \draw[green, line width=1] (3.12,.84) -- (3.12,2);
        \draw[red, line width=1] (3.24,.72) -- (3.24,2);
        \draw[red, line width=1] (2,.76) -- (3.28,.76);

    \end {scope}

    \begin {scope}[shift={(4*0, 4*1)}]
        \draw[blue, line width=1] (0,1) -- (1,1);
        \draw[blue, line width=1] (1,1) -- (2.04,1);
        \draw[blue, line width=1] (1.96,2) -- (3,2);
        \draw[blue, line width=1] (2,.96) -- (2,2);
        \draw[green, line width=1] (1,.88) -- (2.1,.88);
        \draw[green, line width=1] (2,0.88) -- (3.16,0.88);
        \draw[green, line width=1] (3.12,.84) -- (3.12,2);
        \draw[red, line width=1] (3.24,.72) -- (3.24,2);
        \draw[red, line width=1] (2,.76) -- (3.28,.76);

    \end {scope}

    \begin {scope}[shift={(4*-1, 4*2)}]
         \draw[blue, line width=1] (0,1) -- (1.04,1);
        \draw[blue, line width=1] (.96,2) -- (2.04,2);
        \draw[blue, line width=1] (1.96,2) -- (3,2);
        \draw[blue, line width=1] (1,.96) -- (1,2);
        \draw[green, line width=1] (1,.88) -- (2.1,.88);
        \draw[green, line width=1] (2,0.88) -- (3.16,0.88);
        \draw[green, line width=1] (3.12,.84) -- (3.12,2);
        \draw[red, line width=1] (3.24,.72) -- (3.24,2);
        \draw[red, line width=1] (2,.76) -- (3.28,.76);

    \end {scope}

    \begin {scope}[shift={(4*1, 4*2)}]
        \draw[blue, line width=1] (0,1) -- (1,1);
        \draw[blue, line width=1] (1,1) -- (2.04,1);
        \draw[blue, line width=1] (1.96,2) -- (3,2);
        \draw[blue, line width=1] (2,.96) -- (2,2);
        \draw[green, line width=1] (1,.88) -- (2.16,.88);
        \draw[green, line width=1] (2.08,1.88) -- (3,1.88);
        \draw[green, line width=1] (2.12,.84) -- (2.12,1.92);
        \draw[red, line width=1] (2.2,1) --++ (0.8,0) --++ (0,0.8);

    \end {scope}

    \begin {scope}[shift={(4*-2, 4*3)}]
     \draw[blue, line width=1] (-.04,2) -- (1.04,2);
        \draw[blue, line width=1] (.96,2) -- (2.04,2);
        \draw[blue, line width=1] (1.96,2) -- (3,2);
        \draw[blue, line width=1] (0,.96) -- (0,2.04);
        \draw[green, line width=1] (1,.88) -- (2.1,.88);
        \draw[green, line width=1] (2,0.88) -- (3.16,0.88);
        \draw[green, line width=1] (3.12,.84) -- (3.12,2);
        \draw[red, line width=1] (3.24,.72) -- (3.24,2);
        \draw[red, line width=1] (2,.76) -- (3.28,.76);

    \end {scope}

    \begin {scope}[shift={(4*0, 4*3)}]
        \draw[blue, line width=1] (0,1) -- (1.04,1);
        \draw[blue, line width=1] (.96,2) -- (2.04,2);
        \draw[blue, line width=1] (1.96,2) -- (3,2);
        \draw[blue, line width=1] (1,.96) -- (1,2);
        \draw[green, line width=1] (1,.88) -- (2.1,.88);
        \draw[green, line width=1] (2.08,1.88) -- (3,1.88);
        \draw[green, line width=1] (2.12,.84) -- (2.12,1.92);
        \draw[red, line width=1] (2.2,1) --++ (0.8,0) --++ (0,0.8);

    \end {scope}

    \begin {scope}[shift={(4*2, 4*3)}]
        \draw[blue, line width=1] (0,1) -- (1,1);
        \draw[blue, line width=1] (1,1) -- (2.04,1);
        \draw[blue, line width=1] (1.96,2) -- (3,2);
        \draw[blue, line width=1] (2,.96) -- (2,2);
        \draw[green, line width=1] (1,.88) -- (2.16,.88);
        \draw[green, line width=1] (2.08,1.88) -- (3,1.88);
        \draw[green, line width=1] (2.12,.84) -- (2.12,1.92);
        \draw[red, line width=1] (2.24,.82) -- (2.24,1.80);
        \draw[red, line width=1] (2.20,1.76) -- (3,1.76);

    \end {scope}

    \begin {scope}[shift={(4*-1, 4*4)}]
        \draw[blue, line width=1] (-.04,2) -- (1.04,2);
        \draw[blue, line width=1] (.96,2) -- (2.04,2);
        \draw[blue, line width=1] (1.96,2) -- (3,2);
        \draw[blue, line width=1] (0,.96) -- (0,2.04);
        \draw[green, line width=1] (1,.88) -- (2.1,.88);
        \draw[green, line width=1] (2.08,1.88) -- (3,1.88);
        \draw[green, line width=1] (2.12,.84) -- (2.12,1.92);;
        \draw[red, line width=1] (2.2,1) --++ (0.8,0) --++ (0,0.8);

    \end {scope}

    \begin {scope}[shift={(4*0, 4*4)}]
        \draw[blue, line width=1] (0,1) -- (1.04,1);
        \draw[blue, line width=1] (.96,2) -- (2.04,2);
        \draw[blue, line width=1] (1.96,2) -- (3,2);
        \draw[blue, line width=1] (1,.96) -- (1,2);
        \draw[green, line width=1] (1.12,1.88) -- (2.1,1.88);
        \draw[green, line width=1] (2.08,1.88) -- (3,1.88);
        \draw[green, line width=1] (1.12,.95) -- (1.12,1.924);
        \draw[red, line width=1] (2,1) --++ (1,0) --++ (0,0.8);

    \end {scope}

    \begin {scope}[shift={(4*1, 4*4)}]
        \draw[blue, line width=1] (0,1) -- (1.04,1);
        \draw[blue, line width=1] (.96,2) -- (2.04,2);
        \draw[blue, line width=1] (1.96,2) -- (3,2);
        \draw[blue, line width=1] (1,.96) -- (1,2);
        \draw[green, line width=1] (1,.88) -- (2.1,.88);
        \draw[green, line width=1] (2.08,1.88) -- (3,1.88);
        \draw[green, line width=1] (2.12,.84) -- (2.12,1.92);
        \draw[red, line width=1] (2.24,.82) -- (2.24,1.80);
        \draw[red, line width=1] (2.20,1.76) -- (3,1.76);

    \end {scope}

    \begin {scope}[shift={(4*-1, 4*5)}]
        \draw[blue, line width=1] (-.04,2) -- (1.04,2);
        \draw[blue, line width=1] (.96,2) -- (2.04,2);
        \draw[blue, line width=1] (1.96,2) -- (3,2);
        \draw[blue, line width=1] (0,.96) -- (0,2.04);
        \draw[green, line width=1] (1.12,1.88) -- (2.1,1.88);
        \draw[green, line width=1] (2.08,1.88) -- (3,1.88);
        \draw[green, line width=1] (1.12,.95) -- (1.12,1.924);
        \draw[red, line width=1] (2,1) --++ (1,0) --++ (0,0.8);

    \end {scope}

    \begin {scope}[shift={(4*0, 4*5)}]
        \draw[blue, line width=1] (-.04,2) -- (1.04,2);
        \draw[blue, line width=1] (.96,2) -- (2.04,2);
        \draw[blue, line width=1] (1.96,2) -- (3,2);
        \draw[blue, line width=1] (0,.96) -- (0,2.04);
        \draw[green, line width=1] (1,.88) -- (2.16,.88);
        \draw[green, line width=1] (2.08,1.88) -- (3,1.88);
        \draw[green, line width=1] (2.12,.84) -- (2.12,1.92);
        \draw[red, line width=1] (2.24,.82) -- (2.24,1.80);
        \draw[red, line width=1] (2.20,1.76) -- (3,1.76);

    \end {scope}
    
    \begin {scope}[shift={(4*1, 4*5)}]
        \draw[blue, line width=1] (0,1) -- (1.04,1);
        \draw[blue, line width=1] (.96,2) -- (2.04,2);
        \draw[blue, line width=1] (1.96,2) -- (3,2);
        \draw[blue, line width=1] (1,.96) -- (1,2);
        \draw[green, line width=1] (1.12,1.88) -- (2.1,1.88);
        \draw[green, line width=1] (2.08,1.88) -- (3,1.88);
        \draw[green, line width=1] (1.12,.95) -- (1.12,1.924);
        \draw[red, line width=1] (2.12,.95) -- (2.12,1.80);
        \draw[red, line width=1] (2.08,1.76) -- (3,1.76);

    \end {scope}

    \begin {scope}[shift={(4*0, 4*6)}]
        \draw[blue, line width=1] (-.04,2) -- (1.04,2);
        \draw[blue, line width=1] (.96,2) -- (2.04,2);
        \draw[blue, line width=1] (1.96,2) -- (3,2);
        \draw[blue, line width=1] (0,.96) -- (0,2.04);
        \draw[green, line width=1] (1.12,1.88) -- (2.1,1.88);
        \draw[green, line width=1] (2.08,1.88) -- (3,1.88);
        \draw[green, line width=1] (1.12,.95) -- (1.12,1.924);
        \draw[red, line width=1] (2.12,.95) -- (2.12,1.80);
        \draw[red, line width=1] (2.08,1.76) -- (3,1.76);

    \end {scope}


    \begin {scope}[shift={(22,0)}]
        \foreach \x/\y in {0/0, 0/3, -1/4, 1/4, 0/5} {
            \begin {scope}[shift={(4*\x,4*\y)}]
                \draw (1,2.5) -- (1,3.5);
            \end {scope}
        }

        \foreach \x/\y in {0/1, -1/2, 1/2, -2/3, 0/3, -1/4, 0/4, -1/5} {
            \begin {scope}[shift={(4*\x,4*\y)}]
                \draw (1.5,2.5) -- (3.5,4.5);
            \end {scope}
        }

        \foreach \x/\y in {0/1, -1/2, 1/2, 0/3, 2/3, 0/4, 1/4, 1/5} {
            \begin {scope}[shift={(4*\x,4*\y)}]
                \draw (0.5,2.5) -- (-1.5,4.5);
            \end {scope}
        }

        \begin {scope}[shift={(4*0, 4*0)}]
            \draw (1,1) node {$(0,0,0)$};
        \end {scope}

        \begin {scope}[shift={(4*0, 4*1)}]
            \draw (1,1) node {$(0,0,1)$};
        \end {scope}

        \begin {scope}[shift={(4*-1, 4*2)}]
            \draw (1,1) node {$(0,1,1)$};
        \end {scope}

        \begin {scope}[shift={(4*1, 4*2)}]
            \draw (1,1) node {$(0,0,2)$};
        \end {scope}

        \begin {scope}[shift={(4*-2, 4*3)}]
            \draw (1,1) node {$(1,1,1)$};
        \end {scope}

        \begin {scope}[shift={(4*0, 4*3)}]
            \draw (1,1) node {$(0,1,2)$};
        \end {scope}

        \begin {scope}[shift={(4*2, 4*3)}]
            \draw (1,1) node {$(0,0,3)$};
        \end {scope}

        \begin {scope}[shift={(4*-1, 4*4)}]
            \draw (1,1) node {$(1,1,2)$};
        \end {scope}

        \begin {scope}[shift={(4*0, 4*4)}]
            \draw (1,1) node {$(0,2,2)$};
        \end {scope}

        \begin {scope}[shift={(4*1, 4*4)}]
            \draw (1,1) node {$(0,1,3)$};
        \end {scope}

        \begin {scope}[shift={(4*-1, 4*5)}]
            \draw (1,1) node {$(1,2,2)$};
        \end {scope}

        \begin {scope}[shift={(4*0, 4*5)}]
            \draw (1,1) node {$(1,1,3)$};
        \end {scope}

        \begin {scope}[shift={(4*1, 4*5)}]
            \draw (1,1) node {$(0,2,3)$};
        \end {scope}

        \begin {scope}[shift={(4*0, 4*6)}]
            \draw (1,1) node {$(1,2,3)$};
        \end {scope}
    \end {scope}
\end {tikzpicture}
\end{center}
\caption {(Left) The mixed lattice paths on $\mathcal{G}^s_3$. (Right) The tuples $(a_1,a_2,a_3)$, where $a_i$ is the number of paths above tile $T_i$.
The reverse sequences $(a_3,a_2,a_1,0)$ are the Lehmer codes of the 132-avoiding permutations in Figure \ref{fig:catalan_posets}.}
\label {fig:catalan_posets2}
\end {figure}

\begin {ex}
    For the straight snake graph $\mathcal{G}^s_n$, the standard labeling considered in Section \ref{sec:euler_straight_snakes} agrees with the
    one just defined. The canonical path $P_0$ is the minimal lattice path, which follows the bottom boundary of $\mathcal{G}^s_n$. Therefore the
    unique minimal mixed lattice path is the one where all edges are contained in the canonical (i.e. minimal) path. Although a mixed lattice path
    is by definition a multiset of edges and in general there is not a unique way to write it as a union of single lattice paths, in this case there
    is a particularly simple decomposition. We can always write a mixed lattice path as the union of partial paths $L_0,L_1,L_2,\dots,L_n$ which are 
    weakly increasing (i.e. if $i<j$, then $L_i \leq L_j$ in the usual partial order on lattice paths). 
    
    This is illustrated in Figure \ref{fig:catalan_posets2}. In the figure, path $L_1$ is in red, $L_2$ in green, and $L_3$ in blue. In the picture, the
    path $L_0$ (which is just a single edge) is omitted.
\end{ex}

Now we will finally prove Theorem \ref{thm:lattice}, which says that the partial order on mixed dimer covers is a distributive lattice.

\begin {proof}[Proof of Theorem \ref{thm:lattice}]
    Let $\mathcal{L}(\mathcal{G})$ be the set of mixed lattice paths on $\mathcal{G}$ (with its standard labeling), and let $\mathcal{L}_n(\mathcal{G})$
    be the set of \emph{$n$-lattice paths} (i.e. for the vertex labeling where every vertex has label $n$). 
    By Theorem \ref{thm:duality}, we have the isomorphism of posets $\Omega(\mathcal{G}) \cong \mathcal{L}(\widetilde{\mathcal{G}})$. So it suffices to
    show that $\mathcal{L}(\mathcal{G})$ is a distributive lattice for any snake graph $\mathcal{G}$.

    Let $\mathcal{G}$ be a snake graph with $n$ squares. There is a natural injection $\mathcal{L}(\mathcal{G}) \hookrightarrow \mathcal{L}_{n+1}(\mathcal{G})$
    defined as follows. Decompose a mixed lattice path into the union of $n+1$ paths. Each of these paths begins at a different point along the canonical path.
    Extend each of these paths by adding the beginning part of the canonical path to it. This way we obtain $n+1$ paths $L_0,L_1,\dots,L_n$ which all begin 
    at the bottom-left vertex of $\mathcal{G}$, where $L_i$ shares (at least) the first $i$ steps with the canonical lattice path.
    
    It was observed in \cite{mosz} and \cite{bosz_24} that $\mathcal{L}_{n+1}(\mathcal{G})$ is a distributive lattice, isomorphic to the lattice of
    $P$-partitions (with parts at most $n+1$) on the poset whose elements are the boxes of $\mathcal{G}$. A $P$-partition is a labeling of the boxes of $\mathcal{G}$
    by non-negative integers, such that labels weakly increase in both directions (left-to-right and top-to-bottom). Given an $(n+1)$-lattice path, decomposed
    as a union $L_0 \cup L_1 \cup \cdots \cup L_n$, the labels in the boxes of the corresponding $P$-partition are simply the number of $L_i$'s which go above that box.

    Since $\mathcal{L}(\mathcal{G})$ can be seen as a subset of a distributive lattice, it suffices to check that this subset is closed under the
    meet and join operations of the lattice. Note that 
    \[ \mathcal{L}_{n+1}(\mathcal{G}) \cong \mathcal{L}_1(\mathcal{G})^{n+1} / \sim \; = \; \mathcal{L}_1(\mathcal{G}) \times \cdots \times \mathcal{L}_1(\mathcal{G}) / \sim \]
    That is, $\mathcal{L}_{n+1}(\mathcal{G})$ is a quotient of the product of $n+1$ copies of the lattice of \emph{single} lattice paths by the equivalence relation
    $(L_0,\dots,L_n) \sim (L'_0,\dots,L'_n)$ if $\bigcup_i L_i = \bigcup_i L'_i$ as multisets. For two $(n+1)$-lattice paths $(L_0,\dots,L_n)$ and $(L_0',\dots,L_n')$,
    the meet and join are simply component-wise:
    \[ 
        (L_0,\dots,L_n) \wedge (L_0',\dots,L_n') = (L_0 \wedge L_0', \dots, L_n \wedge L_n') \quad \text{and} \quad 
        (L_0,\dots,L_n) \vee (L_0',\dots,L_n') = (L_0 \vee L_0', \dots, L_n \vee L_n')
    \]
    But since $L_i$ and $L_i'$ both share the first $i$ steps in common (both agree with the canonical path for the first $i$ steps), then both $L_i \wedge L_i'$
    and $L_i \vee L_i'$ also share the first $i$ steps in common with the canonical path. This shows that $\mathcal{L}(\mathcal{G})$ is closed under the meet
    and join operations of the enveloping distributive lattice $\mathcal{L}_{n+1}(\mathcal{G})$, and hence $\mathcal{L}(\mathcal{G})$ is itself a distributive lattice.
\end {proof}

\subsection {Dimers, Paths, and Permutations}

Recall from Section \ref{sec:lattices} that the covering relations for the face twist order on $\Omega(\mathcal{G}^s_n)$ and $\Omega(\mathcal{G}^z_n)$
are mapped, under the respective bijections, to the covering relations in the left middle order on permutations. More specifically, performing
a face twist (or dually on lattice paths, performing a bump $RU \mapsto UR$) corresponds to incrementing one entry of the Lehmer code of a permutation.
This is illustrated in the right side of Figure \ref{fig:catalan_posets2}, where the tuples are the Lehmer codes written in reverse order.

We can extend this idea to the case of more general snake graphs. Let $\mathcal{G} = \mathcal{G}^w$ be a snake graph with its standard labeling.
In light of Theorem \ref{thm:duality}, we can consider mixed lattice paths on $\widetilde{\mathcal{G}}$ rather than mixed dimer covers of $\mathcal{G}$.
Define a mapping from the set of mixed lattice paths to the product of intervals $[0,1] \times [0,2] \times \cdots \times [0,n]$ as follows.
For a mixed lattice path $p$, we let $p \mapsto (a_1,a_2,\dots,a_n)$, where in any maximal chain from the minimal mixed lattice path to $p$,
a total of $a_i$ flips are performed in tile $T_i$. It is easy to see that this is independent of the maximal chain.

Note that by definition of mixed lattice path, there are at most $i$ paths beginning at vertices before tile $T_i$, and
so there can be at most $i$ flips in tile $T_i$. This shows that the mapping is well-defined. Reversing the order of this tuple
$(a_1,\dots,a_n) \mapsto (a_n, \dots a_1)$ gives an element of $[0,n] \times [0,n-1] \times \cdots \times [0,1]$. We may then append
this sequence with an extra $0$ to get an element of $[0,n] \times \cdots \times [0,1] \times [0,0]$, which can be interpreted as
the Lehmer code of a permutation $\sigma \in S_{n+1}$. By construction, the partial order on mixed lattice paths on $\widetilde{\mathcal{G}}$
(or equivalently mixed dimer covers of $\mathcal{G}$) will be isomorphic to the restriction of the left middle order
on this set of permutations.
To summarize,
the preceding discussion establishes the following.

\begin {thm} \label{thm:dimer_permutations}
    Let $\mathcal{G}$ be a snake graph with $n-1$ tiles, endowed with its standard labeling. Then there is a set of permutations $S_{\mathcal{G}} \subseteq S_n$
    such that the partial order on $\Omega(\mathcal{G})$ (and by Theorem \ref{thm:duality}, also $\mathcal{L}(\widetilde{G})$) is isomorphic to the 
    left middle order on $S_{\mathcal{G}}$.
\end {thm}

\begin {rmk}
    While the theorem is only stated for the standard labeling, the construction still makes sense whenever tile $T_i$ admits at most $i$ flips.
    In particular, any vertex labeling where every vertex label is less than or equal to the standard one will produce sequences which can be
    interpreted as Lehmer codes.
\end {rmk}

When $\widetilde{G}$ is a straight snake graph, this gives precisely the middle order on 132-avoiding permutations (as illustrated in Figures \ref{fig:catalan_posets}
and \ref{fig:catalan_posets2}). When $\widetilde{G}$ is a zigzag snake graph, this construction doesn't quite give the set of alternating permutations,
but it is easy to say how it is related. The minimal alternating permutation has Lehmer code $(1,0,1,0,\dots)$ with alternating 1's and 0's,
and so the Lehmer codes of the alternating permutations differ from the ones given by this construction by subtracting the vector $(1,0,1,0,\dots)$.

In these two extreme cases of straight and zigzag snakes, we obtain interesting classes of permutations. So it seems to be an interesting
problem to characterize for more general snake graphs exactly which sets of permutations we obtain from this construction. 

\begin {prob}
    For a general snake graph $\mathcal{G}$ (with the standard labeling), characterize the corresponding set $S_\mathcal{G}$ of permutations.
\end {prob}

Another interesting example to consider is that of ordinary single dimer covers. \c{C}anak\c{c}i and Schroll showed in \cite{cs_21} that
the lattice of single dimer covers of a snake graph is isomorphic to an order ideal in the weak order on $S_n$ generated by
a Coxeter element (whose reduced decomposition into adjacent transpositions is determined by the shape of the snake).

While similar in spirit, Theorem \ref{thm:dimer_permutations} above gives something different in this case.
In the construction from \cite{cs_21}, a face twist in tile $T_i$ (which are the covering relations) corresponds to composing with
the adjacent transposition $s_i = (i,i+1)$. In our setup described above, the face twist in tile $T_i$ corresponds to incrementing
the $i^\mathrm{th}$ entry of the Lehmer code. So there is a poset isomorphism from the \c{C}anak\c{c}i-Schroll poset to ours, where a reduced word
$\sigma = s_{i_1} s_{i_2} \cdots s_{i_k}$ is mapped to a permutation whose Lehmer code is the 0-1 vector with $L_i = 1$ if and only if $s_i$
appears in the reduced word for $\sigma$.

Note that if a permutation's Lehmer code consists of only 0's and 1's, then it can be written in cycle notation with the numbers $1,2,3,\dots,n$ appearing
in order from left-to-right. Indeed, a substring $111 \cdots 10$ of the Lehmer code, starting at position $i$, and having $k$ $1$'s followed by a $0$,
implies that part of the permutation's one-line notation is $i+1, i+2, \dots, i+k, i$, which means it contains the cycle $(i,i+1,\dots,i+k)$.
For example, the permutation $\sigma = (1)(234)(56)(789)$ has Lehmer code $011010110$.

\begin{figure}
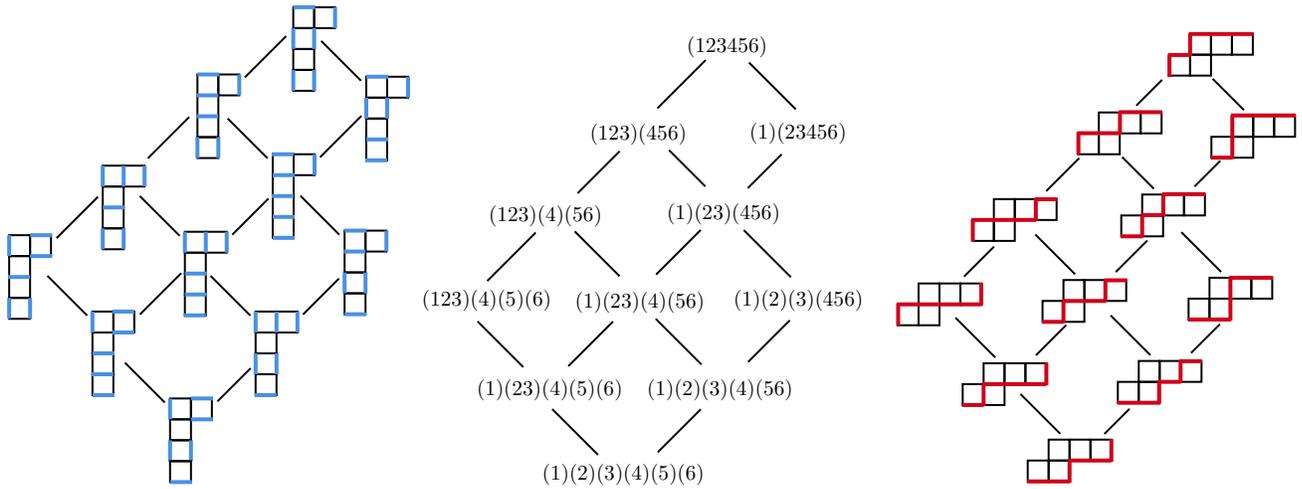

    \centering
    
\tikzset{every picture/.style={line width=0.75pt}} 

\tikzset{every picture/.style={line width=0.75pt}} 


\caption {(left) the poset $\Omega_1(\mathcal{G})$, (middle) the poset $S_{\mathcal{G}}$, (right) the poset $\mathcal{L}(\widetilde{\mathcal{G}})$}
\label {fig:single_dimer_poset}
\end{figure}

For single dimer covers, the poset $S_{\mathcal{G}}$  has bottom element the identity,
top element the long cycle $(123 \cdots n)$, and in between, the covering relations are given by joining adjacent cycles. The restrictions of when cycles
can be joined is dictated by the shape of the snake graph. Specifically, performing a face twist of a dimer cover at tile $T_i$ corresponds
to joining the cycle of $\sigma$ ending with $i$ and the cycle beginning with $i+1$. An example is illustrated in Figure \ref{fig:single_dimer_poset},
which shows both the dimer covers of $\mathcal{G} = \mathcal{G}^{UUUR}$ and the lattice paths on $\widetilde{\mathcal{G}}$.

\section {Acknowledgments}

Nicholas Ovenhouse was partially supported by the Simons Foundation grant 327929. We would like to thank the following people for
helpful comments and enlightening discussions: Seok Hyun Byun, Luca Ferrari, Richard Kenyon, Gregg Musiker, James Propp, Bridget Tenner, and Sylvester W. Zhang.

\vfill

\bibliographystyle{alpha}
\bibliography{refs.bib}

\end {document}